\documentclass[a4paper,12pt]{amsart}
\usepackage{amsmath}
\usepackage{amsthm}
\usepackage{amscd}
\usepackage{amssymb}
\usepackage{amsfonts}
\usepackage{latexsym}
\usepackage[mathscr]{eucal}

\newtheorem{thm}{Theorem}[section]
\newtheorem{prop}[thm]{Proposition}
\newtheorem{lemma}[thm]{Lemma}
\newtheorem{cor}[thm]{Corollary}

\def\graph{\mathop{\rm {graph}}\nolimits}

\def\Im{\mathop{\rm {Im}}\nolimits}
\def\codim{\mathop{\rm {codim}}\nolimits}
\def\Index{\mathop{\rm {Index}}\nolimits}
\def\Ker{\mathop{\rm {Ker}}\nolimits}
\def\cograph{\mathop{\rm {cograph}}\nolimits}

\begin{document}
\title{Relative position of four subspaces
in a Hilbert space}
\author{Masatoshi Enomoto}
\address[Masatoshi Enomoto]{College of Business Administration 
and Inforation Science, 
Koshien University, Takarazuka, Hyogo 665, Japan}      

\author{Yasuo Watatani}
\address[Yasuo Watatani]{Department of Mathematical Sciences, 
Kyushu University, Hakozaki, 
Fukuoka, 812-8581,  Japan}
\maketitle
\begin{abstract}
The relative position of one subfactor of a factor 
has been proved quite rich since the work of Jones. 
We shall show that the theory of relative position of
{\it several subspaces } of a separable infinite-dimensional Hilbert
 space is also rich. 
In finite-dimensonal case, 
Gelfand and Ponomarev  gave a complete classification of 
indecomposable systems of four subspaces.  
We construct exotic examples of indecomposable 
systems of four subspaces in infinite-dimensional Hilbert spaces.  
We extend  their Coxeter functors  and defect using Fredholm index. 
There exist close connections with strongly 
irreducible operators and transitive lattices.


\medskip\par\noindent
AMS SUBJECT CLASSIFICATION: 46C07, 47A15, 15A21, 16G20, 16G60.

\end{abstract}

\section{Introduction}
The relative position of one subfactor of a factor 
has been proved quite rich since the work of Jones [J] .  
On the other hand, the relative position of one
subspace of a Hilbert space is extremely simple and determined 
by the dimension and the co-dimension of the subspace.
But we shall show that the theory of relative position of
{\it several subspaces } of a Hilbert space is rich 
as subfactor theory. 

It is a well known fact that the relative position of  two subspaces
$E$ and $F$ in a Hilbert space $H$  can be described completely up to
unitary equivalence  as in Araki \cite{Ar}  Dixmier \cite{D} and 
Halmos \cite{Ha1}.  The Hilbert
space is the direct sum of five subspaces:
$$
H = (E \cap F) \oplus (\text{the rest}) \oplus (E \cap F^{\perp})
\oplus (E^{\perp} \cap F) \oplus (E^{\perp} \cap F^{\perp}).
$$
In the rest part,  $E$ and $F$ are in generic position and the
relative position is described only by \lq\lq the angles" between them.

We disregard \lq\lq the angles" and study the still-remaining fundamental 
feature of the relative position of $n$ subspaces.   
As it is important to study {\it irreducible} subfactors 
in subfactor theory, we should study an {\it indecomposable} system 
of $n$ subspaces in the sense that the system can not be isomorphic 
to a direct sum of two  non-zero systems.

On the other hand, many problems of linear algebra can be reduced 
to the classification of the systems of subpaces in a 
finite-dimensional vector space.  In a finite-dimensional space,  
the classification of indecomposable systems of $n$ subspaces 
for $n = 1,2$ and $3$  was simple. Jordan blocks give 
indecomposable systems of 4 subspaces.  But there exist 
many other kinds of indecomposable systems of 4 subspaces. 
Therefore it was surprising that  
Gelfand and Ponomarev \cite{GP} gave a complete classification of 
indecomposable systems of four subspaces in a finite-dimensional 
space over an algebraically closed field. 

In this note we study relative position of $n$ subspaces
in a separable  infinite-dimensional Hilbert space. The fact 
that the sum of closed subspaces is not necessary closed causes some 
troubles in several arguments in Gelfand-Ponomarev \cite{GP}.  
Let $H$ be a Hilbert space and $E_1, \dots E_n$ be $n$ subspaces 
in $H$.  Then we say that  ${\mathcal S} = (H;E_1, \dots , E_n)$  
is a system of $n$ subspaces in $H$ or a $n$-subspace system in $H$.
A system ${\mathcal S}$ is called indecomposable if ${\mathcal S}$ 
can not be decomposed into a nontrivial direct sum.  
For any bounded linear operator $A$ on a Hilbert space $K$, we can 
associate a system ${\mathcal S}_A$ of four subspaces in 
$H = K \oplus K$ by 
\[
{\mathcal S}_A = (H;K\oplus 0,0\oplus K,\graph A, \{(x,x) ; x \in K\}).
\]
Two such systems ${\mathcal S}_A$ and ${\mathcal S}_B$  are isomorphic 
if and only if the two operators $A$ and $B$ are similar. The 
direct sum of such  systems corresponds to the direct sum of 
the operators. In this sense the theory of operators is 
included into the theory of relative positions of four subspaces. 
In particular on a finite dimesional space, Jordan blocks correspond 
to indecomposable systems. Moreover on an infinite dimensional 
Hilbert space, the above system ${\mathcal S}_A$ is indecomposable 
if and only if $A$ is strongly irreducible, which is an 
infinite-dimensional analog of a Jordan block, see, for example, 
a monograph by Jiang and Wang \cite{JW}. 
Therefore there exist uncountably many indecomposable 
systems of four subspaces. But it is rather difficult to know 
whether there exists another kind of indecomposable system of four 
subspaces. One of the main result of the paper is 
to give uncountably many, exotic,  indecomposable systems 
of four subspaces on an infinite-dimensional separable Hilbert space. 
The $\ell ^2$-boundedness is crucially used.    

Gelfand and Ponomarev introduced an integer valued  invariant 
$\rho ({\mathcal S})$, 
called  
{\it defect},  for a system 
${\mathcal S} = (H ; E_1, E_2, E_3, E_4)$ of four subspaces by 
\[
\rho ({\mathcal S}) = \sum _{i=1} ^4 \dim E_i - 2\dim H.  
\]
We extend the defect to a certain class of 
systems of four subspaces on an infinite dimesional Hilbert space 
using  Fredholm index.  
We believe that there exists an analogy between a classification of 
systems of subspaces and a classification of subfactors, and 
the defect by Gelfand and Ponomarev seems to correspond  
 to the index by Jones \cite{J}. Therefore the 
determination of possible value of defect is also important.
If a pair $N \subset M$ of factor-subfactor is finite-dimensional, 
then Jones index $[M:N]$ is an integer.  But if $N \subset M$ is infinite-dimensional, then Jones index $[M:N]$ is  a  non-integer in general. 
One of the amazing fact was that the possible value of Jones index
is  in $\{4\cos ^2\frac{\pi}{n} \ | \ n = 3,4,...\} \cup [4,\infty]$.
We show that a similar situation occurs for the possible value 
of defect. If a system $\mathcal S = (H;E_1,E_2,E_3,E_4)$ of 
four subspaces is finite-dimensional,
then the defect $\rho (\mathcal S)$ is an integer.  Gelfand and Ponomarev 
showed that the possible value of defect $\rho (\mathcal S)$ is exactly in
$\{-2, -1, 0, 1, 2\}$.   We show that the set of values of defect for 
indecomposable systems of four subspaces in an infinite-dimesional 
Hilbert spaces is exactly  $\{ \frac{n}{3} ; n  \in \mathbb Z  \}$. 

We extend Coxeter functors after Gelfand-Ponomarev and show that the Coxeter 
functors preserve the defect and indecomposability under certain conditions. 

Halmos initiated the study of transitive lattices and gave an example of 
transitive lattice consisting of seven subspaces in \cite{Ha2}. 
Harison-Radjavi-Rosenthal \cite{HRR} constructed a transitive lattice  
consisting of six subspaces using the graph of an unbounded closed operator. 
Hadwin-Longstaff-Rosenthal found a transitive lattice of five  non-closed 
linear subspaces in \cite{HLR}.     
Any finite transitive lattice which consists of  $n$ subspaces of a Hilbert 
space $H$  
gives an indecomposable system of $n-2$ subspaces by withdrawing 
$0$ and $H$,  but the converse is not true. It is still unknown whether or not 
there exists a transitive lattice 
consisting of five subspaces.  Therefore it is also an interesting  problem to know whether there exists an indecomposable system of three subspaces in an infinite-dimensional Hilbert space.  
 
Throughout the paper a projection means an operator 
$e$ with $e^2 = e = e^*$ and 
an idempotent  means an operator $p$ with $p^2 = p$. 
 
Sunder also considered $n$ subspaces in \cite{S}.  
But his interest is extremely 
opposite to ours.  In fact he studied the decomposable case such that the 
Hilbert space $H$ is an algebraic sum of the $n$ subspaces. 
He solved the statistical problem of computing the canonical 
partial correlation coefficients between three sets of random 
variables. 

When we announced some part of our result in US-Japan seminar 
at Fukuoka in 1999, we had not yet known the notion and interesting  
works on strong irreducible operators which are summarized in a monograph 
by Jiang and Wang \cite{JW}.

There seems to be interesting relations with the study of 
representations of $*$-algebras generated by idempotents by 
S. Kruglyak and Y. Samoilenko \cite{KS} and the study on sums of 
projections by S. Kruglyak, V. Rabanovich and Y. Samoilenko 
\cite{KRS}. But we do not know the exact implication, because 
their objects are different with ours. 

In finite dimensional case, the classification of four subspaces is 
described as  the classification of the representations  of the 
extended Dynkin diagram $D_4^{(1)}$. Recall that Gabriel \cite{G} 
listed Dynkin diagrams $A_n, D_n, E_6, E_7,E_8$ in his theory on 
finiteness of indecomposable representations of quivers. We will 
discuss on indecomposable representations of 
quivers on {\it infinite-dimensinal Hilbert spaces} 
somewhere else \cite{EW} as a continuation of this paper.  

In purely algebraic setting, it is known that if a finite-dimensional 
algebra $R$ is not of representation-finite type, then there exist 
indecomposable 
$R$-modules of infinite length as in M. Auslander \cite{Au}. 
Since  we consider 
representations on Hilbert spaces, the result in \cite{Au} cannot be applied 
directly.  We need several techniques in functional analysis. See a book 
\cite{KR} for infinite length modules. 

The authors are supported by the Grant-in-Aid for Scientific Research of 
JSPS.

\section{systems of $n$ subspaces}
We study the relative position of $n$ subspaces
in a separable  Hilbert space. 
Let $H$ be a Hilbert space and $E_1, \dots E_n$ be $n$ subspaces 
in $H$.  Then we say that  ${\mathcal S} = (H;E_1, \dots , E_n)$  
is a system of $n$-subspaces in $H$ or a $n$-subspace system in $H$. 
Let ${\mathcal T} = (K;F_1, \dots , F_n)$  
be  another system of $n$-subspaces in a Hilbert space $K$. Then  
$\varphi : {\mathcal S} \rightarrow {\mathcal T}$ is called a 
homomorphism if $\varphi : H \rightarrow K$ is a bounded linear 
operator satisfying that  
$\varphi(E_i) \subset F_i$ for $i = 1,\dots ,n$. And 
$\varphi : {\mathcal S} \rightarrow {\mathcal T}$
is called an isomorphism if $\varphi : H \rightarrow K$ is 
an invertible (i.e., bounded  bijective) linear 
operator satisfying that  
$\varphi(E_i) = F_i$ for $i = 1,\dots ,n$. 
We say that systems ${\mathcal S}$ and ${\mathcal T}$ are 
{\it isomorphic} if there is an isomorphism  
$\varphi : {\mathcal S} \rightarrow {\mathcal T}$. This means 
that the relative positions of $n$ subspaces $(E_1, \dots , E_n)$ in $H$ 
and   $(F_1, \dots , F_n)$ in $K$ are same under disregarding angles. 
We say that systems ${\mathcal S}$ and ${\mathcal T}$ are 
{\it unitarily equivalent } if the above isomorphism 
$\varphi : H \rightarrow K$ can be chosen to be a unitary. 
This means that the relative positions of 
$n$ subspaces $(E_1, \dots , E_n)$ in $H$ 
and  $(F_1, \dots , F_n)$ in $K$ are same with preserving the angles
between the subspaces. We are interested in the relative position of 
subspaces up to isomorphims to study the still-remaining fundamental 
feature of the relative position after disregarding 
 \lq\lq the angles" .  

We denote by 
$Hom(\mathcal S, \mathcal T)$ the set of homomorphims of 
$\mathcal S$ to $\mathcal T$ and  
$End(\mathcal S) := Hom(\mathcal S, \mathcal S)$ 
the set of endomorphisms on $\mathcal S$. 

Let $G_2 = \mathbb Z/ 2\mathbb Z * \mathbb Z/ 2\mathbb Z =  
\langle a_1, a_2 \rangle$
be the free product of the cyclic groups of order two with generators 
$a_1$ and $a_2$. 
For two subspaces $E_1$ and $E_2$ of a Hilbert space $H$, 
let $e_1$ and $e_2$ be the projections onto $E_1$ and $E_2$.  
Then  $u_1 = 2e_1 -I$ and $u_2 = 2e_2 - I$ are self-adjoint unitaries.   
Thus there is a bijective correspondence between the set 
${\mathcal Sys}^2(H)$ of systems 
${\mathcal S} = (H;E_1,E_2)$ of two subspaces in a Hilbert space $H$ 
and the set $Rep(G_2,H)$ of unitary representations $\pi$ of $G_2$ on $H$ 
such that  $\pi(a_1) = u_1$ and $\pi(a_2) = u_2$.
Similarly let $G_n = \mathbb Z/ 2\mathbb Z * ... * \mathbb Z/ 2\mathbb Z$
be the $n$-times free product of the cyclic groups of order two. Then
there is a bijective correspondence between the set 
${\mathcal Sys}^n(H)$ of systems of $n$
subspaces in a Hilbert space $H$  and the set 
$Rep(G_n,H)$ of unitary representations on $H$ of $G_n$ on $H$.  
It is well known that if $n \geq 3$,
then the group $G_n$ is non-amenable. We should be careful that 
even if two systems of 
$n$ subspaces are isomorphic, the corresponding unitary representations 
are {\it not} necessary to be similar, although the converse is always 
true.   

\noindent
{\bf Example 1} Let $H = \mathbb C ^2$. Fix an angle $\theta$
with $0 < \theta < \pi /2$.  Put $E_1 = \mathbb C(1,0)$ and
$E_2 = \mathbb C(cos\theta, sin\theta)$.  Then
${\mathcal S}_1 = (H;E_1,E_2)$ is isomorphic to 
${\mathcal S}_2 = 
({\mathbb C}^2 ; {\mathbb C} \oplus 0, 0 \oplus {\mathbb C})$.
But the corresponding two unitary representations $\pi_1$ and 
$\pi_2$ are not similar, because 
$\frac{1}{2}(\pi_1(a_1) + 1) \frac{1}{2}(\pi_1(a_2) + 1)
 \not= 0$
and 
$\frac{1}{2}(\pi_2(a_1) + 1) \frac{1}{2}(\pi_2(a_2) + 1) =  0$.

We start with a known fact to recall some notation.  
   
\begin{lemma}  Let $H$ be a Hilbert space and $H_1$ and $H_2$ be 
two subspaces of $H$.  Then the following are equivalent:
\begin{enumerate}
\item $H = H_1 + H_2$ \ and  \ $H_1 \cap H_2 = 0$.
\item There exists a closed subspace $M \subset H$ such that
$(H;H_1,H_2)$ is isomorphic to  $(H;M,M^{\perp})$
\item There exists an idempotnet $P \in B(H)$ such that 
$H_1 = \Im P$ and $H_2 = \Im (1-P)$. 
\end{enumerate}
\label{lemma:decompose}
\end{lemma}

\begin{proof}The equivalence between (1) and (3) is trivial and 
it is immediate that (2)$\Rightarrow$(1) . We show that 
(1)$\Rightarrow$(2).     
Assume (1) and put $M = H_1$. 
Let $e_1$ be the (orthogonal) projection onto $H_1$. Let $P$ be 
the idempotent onto $H_1$ along $H_2$, so that 
$P\xi = \xi_1$ for $\xi = \xi_1 + \xi_2$, 
$(\xi_1 \in H_1, \xi_2 \in H_2)$.  Define an operator 
$T : H \rightarrow H$ by $T\xi = P\xi + (I-e_1)(I-P)\xi$ for $\xi \in H$. 
The operator $P$, $T$ and $T^{-1}$  are also writen as operator matrices 
\[
P = 
\begin{pmatrix}
     I & B \\
     0 & 0 
    \end{pmatrix} 
\text{ , }
T = 
\begin{pmatrix}
     I & B \\
     0 & I 
\end{pmatrix} 
\text{ and }
T^{-1} = 
\begin{pmatrix}
     I & -B \\
     0 & I 
    \end{pmatrix} 
\]
under the decomposion $H = H_1 \oplus H_1^{\perp}$. 
Thus $T$ is an invertible bounded linear operator satisfying 
$TH_1 = H_1$ and $TH_2 = H_1^{\perp}$. Hence $T$ gives an isomorphism.   
\end{proof}

\begin{lemma} Let $H$ and $K$ be Hilbert spaces and 
$E \subset H$ and $F \subset K$ be closed subspaces of $H$ and $K$.  
Let $e \in B(H)$ and $f \in B(K)$ be the projections onto $E$ and $F$. 
Then the following are equivalent:
\begin{enumerate}
\item There exists an invertible operator $T: H \rightarrow K$
such that $T(E)=F$. 
\item There exists an invertible operator $T: H \rightarrow K$
such that $e =(T^{-1}fT)e$ and $f=(TeT^{-1})f$. 
\end{enumerate}
\end{lemma}

\begin{proof}
(1)$\Rightarrow$(2):Assume there exists an invertible operator 
$T: H \rightarrow K$ such that $T(E)=F$. Then for any $\xi \in H$, 
$Te(\xi) \in T(E) = F$.  Hence $f(Te(\xi)) = Te(\xi)$.  
Thus $T^{-1}fTe = e$.  Similarly we have $f=TeT^{-1}f$. \\
$(2)\Rightarrow$ (1):Assume (2). For $\xi \in E$, 
$T(\xi) = Te(\xi)=fTe(\xi) \in F$. Thus $T(E) \subset F$.  
Similarly $T^{-1}(F) \subset E$. Hence $F \subset T(E)$. 
Therefore $T(E) = F$.  
\end{proof}

Using the above lemma, we can describe an isomorphism between 
two systems of $n$ suspaces in terms of operators only as follows:

\begin{cor}
Let $\mathcal S=(H;E_{1},\cdots,E_{n})$
and $\mathcal S^{\prime}
 =(H^{\prime};E_{1}^{\prime},\cdots ,E_{n}^{\prime})$
be two systems of n-subspaces. Let $e_i$ 
(resp. $e_{i}^{\prime}$) be the projection onto 
$E_i$ (resp. $E_{i}^{\prime}$) .
Then two systems $\mathcal S$ and $\mathcal S^{\prime}$ are 
isomorphic if and only if there exists an invertible operator
$T: H \rightarrow H^{\prime}$
such that $e_i =(T^{-1}e_{i}^{\prime}T)e_i$ and 
$e_{i}^{\prime} =(Te_iT^{-1})e_{i}^{\prime}$ for 
$i = 1, \dots, n$.  
\end{cor}

\noindent
{\bf Remark.} If there exists an invertible operator
$T: H \rightarrow H^{\prime}$
such that $e_{i}^{\prime}=Te_iT^{-1}$ for 
$i = 1, \dots, n$, then 
two systems $\mathcal S$ and $\mathcal S^{\prime}$ are 
isomorphic.  But the converse is not true as in example 1.  

We often want to disregard the order of the subspaces. 

\bigskip
\noindent
{\bf Definition} 
Let $\mathcal S=(H;E_{1},\cdots,E_{n})$
and $\mathcal S^{\prime}
 =(H^{\prime};E_{1}^{\prime},\cdots ,E_{n}^{\prime})$
be two systems of $n$-subspaces. Then we say that $\mathcal S$ and 
$\mathcal S^{\prime}$ are isomorphic up to a permutation of subspaces 
if there exists a permutation $\sigma$ on $\{1,2,\dots,n\}$ such that 
$\sigma (\mathcal S) := (H;E_{\sigma(1)},\cdots,E_{\sigma(n)})$ and 
$\mathcal S^{\prime}=(H^{\prime};E_{1}^{\prime},\cdots ,E_{n}^{\prime})$
are isomorphic, i.e., there exists a bounded invertible operator 
$\varphi : H \rightarrow H^{\prime}$ satisfying that  
$\varphi(E_{\sigma(i)}) = E_{i}^{\prime}$ for $i = 1,\dots ,n$.

\section{indecomposable systems}
In this section we  shall introduce a notion of indecomposable system,  
that is, a system which cannot be decomposed into a direct 
sum of smaller systems anymore. 

\bigskip
\noindent
{\bf Definition} (direct sum)
Let ${\mathcal S} = (H;E_1, \dots , E_n)$ and 
 $\mathcal S^{\prime}=
 (H^{\prime};E_{1}^{\prime},\cdots,E_{n}^{\prime})$ be 
systems of $n$ subspaces in Hilbert spaces 
 $H$ and $H^{\prime}$.  Then their direct sum 
$\mathcal {S} \oplus\mathcal {S}^{\prime}$ is defined by 
\[
\mathcal {S} \oplus\mathcal {S}^{\prime}
:= (H\oplus H^{\prime};
E_{1}\oplus E_{1}^{\prime},\dots,E_{n}\oplus E_{n}^{\prime}).
\]

\bigskip
\noindent
{\bf Definition}(indecomposable system)
A system $\mathcal S=(H;E_{1},\dots,E_{n})$
of $n$ subspaces is called {\it decomposable} 
if the system $\mathcal S$ is isomorphic to 
a direct sum of two non-zero systems.  
A system $\mathcal S=(H;E_{1},\cdots,E_{n})$ is said to be 
{\it indecomposable} if it is not decomposable. 

\bigskip
\noindent
{\bf Example 2.} Let $H = \mathbb C ^2$. Fix an angle $\theta$
with $0 < \theta < \pi /2$.  Put $E_1 = \mathbb C(1,0)$ and
$E_2 = \mathbb C(cos\theta, sin\theta)$.  Then
$(H;E_1,E_2)$ is isomorphic to 
\[ 
({\mathbb C}^2 ; {\mathbb C}\oplus 0, 0 \oplus {\mathbb C}) 
\cong (\mathbb C; \mathbb C, 0) \oplus (\mathbb C; 0, \mathbb C).
\]
Hence $(H;E_1,E_2)$ is decomposable. 
 
\bigskip
\noindent
{\bf Remark.} Let $e_1$ and $e_2$ be the
projections onto $E_1$ and $E_2$ in the example 2 above.  
Then the $C^*$-algebra 
$C^*(\{e_1,e_2\})$ generated by $e_1$ and $e_2$ is exactly $B(H)
\cong M_2(\mathbb C)$. Therefore 
the irreducibility of 
$C^*(\{e_1,e_2\})$ does {\it not} imply the 
indecomposability of $(H;E_1,E_2)$. Thus seeking an 
indecomposable system of subspaces is much more difficult and 
fundamental task than showing irreducibility of 
the $C^*$-algebra generated by the corresponding projectios for 
the subspaces.   
\bigskip

We can characterize decomposability of systems inside the 
ambient Hilbert space. 

\begin{lemma}  Let $H$ be a Hilbert space and 
${\mathcal S} = (H;E_1, \dots , E_n)$  a system of $n$ subspaces.  
Then the following condition are equivalent:
\begin{enumerate}
\item ${\mathcal S}$ is decomposable.
\item there exist non-zero closed
subspaces $H_1$ and $H_2$ of $H$ such that $H_1 + H_2 = H$,
$H_1 \cap H_2 = 0$ and $E_i = E_i \cap H_1 + E_i \cap H_2$
for $i= 1,\ldots , n$.
\end{enumerate}
\end{lemma}
\begin{proof}(1)$\Rightarrow$ (2): It is trivial. 
(2)$\Rightarrow$ (1): Assume (2). By \ref{lemma:decompose}, 
there exist a closed subspace $M \subset H$ (in fact we can 
choose $M = H_1$)  and an invertible operator $T \in B(H)$ 
such that $T(H_1) = M$ and $T(H_2) = M^{\perp}$.  Then 
${\mathcal S}$ is isomorphic to a direct sum 
\[
(M; T(E_1 \cap H_1), \dots , T(E_n \cap H_1)) \oplus 
 (M^{\perp}; T(E_1 \cap H_2), \dots , T(E_n \cap H_2)). 
\]
\end{proof}

We give a condition of decomposability  in terms of 
endomorphism algebras for the systems.

\begin{lemma} Let $H$ be a Hilbert space and  
${\mathcal S} = (H;E_1, \dots , E_n)$  a system of 
$n$ subspaces in $H$.  Let $e_i$ be the projection onto $E_i$.
Then the following are equivalent:
\begin{enumerate}
\item There exist non-zero closed subspaces 
$H_{1},H_{2}\subset H$
such that
$H=H_{1}+H_{2}, \ 
H_{1}\cap H_{2}=(0)$ and 
$E_i=E_i\cap H_{1}+ E_i\cap H_{2}$,  $(i = 1, \dots , n)$. 
\item There exists a non-trivial idempotent 
$R\in B(H)$ such that $R(E_i)\subset E_i$,  $(i = 1, \dots , n)$. 
\item There exists a non-trivial idempotent 
$R\in B(H)$ such that $e_iRe_i=Re_i$, $(i = 1, \dots , n)$ . 
\end{enumerate}
\label{lemma:endo-decompose}
\end{lemma}
\begin{proof}
$(1)\Rightarrow(2)$: Assume (1). Let $R$ be the idempotent onto 
$H_1$ along $H_2$. For any $\xi \in E_i$, there exist
$\xi_{1} \in E_i \cap H_{1}$ and
$\xi_{2} \in E_i \cap H_{2}$ 
such that $\xi= \xi_{1}+ \xi_{2}$. Then 
$R(\xi)=\xi_{1} \in E_i$. Thus $R(E_i)\subset E_i.$ \\
$(2)\Rightarrow$(1): Assume (2). We put $H_{1}= \Im R $ and 
$H_{2}=\Im(I-R)$. For $\xi \in E_i$,  we have 
$\xi=R(\xi)+(I-R)(\xi)$. 
Since $R(E_i) \subset E_i$, $R(\xi) \in E_i$. Then 
$(I-R)(\xi) = \xi-R(\xi)\in E_i$. 
Thus $E_i \subset E_i \cap H_{1}+E_i \cap H_{2}.$ The other 
inclusion \lq\lq $\supset$" is trivial. 
$(2)\Leftrightarrow(3):$ It is trivial. 
\end{proof}

We put $Idem(\mathcal S) := \{T \in End(\mathcal S) ; T = T^2 \}$. 

\begin{cor}
Let $\mathcal S = (H;E_1,\ldots ,E_n)$ be a system
of $n$ subspaces in a Hilbert space $H$. 
Then $\mathcal S$ is indecomposable if and only if 
$Idem(\mathcal S) = \{0,I\}$. 
\end{cor}

\begin{cor}
Let $\mathcal S = (H;E_1,\ldots ,E_n)$ be a system
of $n$ subspaces in a Hilbert space $H$. 
Let $e_i$ be the projection of $H$ onto $E_i$ for $i = 1,\ldots ,n$.  
If $\mathcal S = (H;E_1,\ldots ,E_n)$
is indecomposable, then the $C^*(\{e_1,\ldots, e_n\})$ generated by
$e_1, \ldots, e_n$ is irreducible.  But the converse is not true.
\end{cor}

\begin{cor}
Let $\mathcal S = (H;E_1,\ldots ,E_n)$ be a system
of $n$ subspaces in a Hilbert space $H$. 
Let $e_i$ be the projection of $H$ onto $E_i$ for $i = 1,\ldots ,n$.  
Let $P$ be a closed subspace of $H$ and $p$ the projection of $H$ 
onto $P$. If $p$ commutes with any $e_i$, then 
\[
E_i = E_i \cap P + E_i \cap P^{\perp}
\]
 
\end{cor}
\begin{proof}
The projection $R$ of $H$ onto $P$ satisfies the condition (3) 
in \ref{lemma:endo-decompose}. 
\end{proof}

\noindent
{\bf Definition.} Let $\mathcal S = (H;E_1,\ldots ,E_n)$ 
be a system of $n$ subspaces in a Hilbert space $H$. 
Let $e_i$ be the projection of $H$ onto $E_i$ for $i = 1,\ldots ,n$.  
We say that $\mathcal S$ is a {\it commutative} system 
if the $C^*(\{e_1,\ldots, e_n\})$ generated by
$e_1, \ldots, e_n$ is commutative. Be carefull that commutativity 
is {\it not} an isomorphic invariant as shown in Example 1. But 
it is meaningful that a system is isomorphic to a commutative 
system.
 
\begin{prop}Let $\mathcal S = (H;E_1,\ldots ,E_n)$ 
be a system of $n$ subspaces in a Hilbert space $H$.
Assume that $\mathcal S$ is a commutative system. 
Then $\mathcal S$ is indecomposable if and only if 
$\dim H = 1$. Moreover each subset $\Lambda \subset \{1,\dots,n\}$ 
corresponds to a commutative system satisfying $\dim E_i = 1$ for 
$i \in \Lambda$ and $\dim E_i = 0$ for $i \notin \Lambda$. 
\label{prop:commutative}
\end{prop}
\begin{proof}Let $e_i$ be the projection of $H$ onto 
$E_i$ for $i = 1,\ldots ,n$.  
If $\mathcal S$ is a commutative, indecomposable system, then  
the $C^*(\{e_1,\ldots, e_n\}) \subset B(H)$ is commutative 
and irreducible.  Thus $\dim H = 1$.  The converse and the rest 
is clear.
\end{proof}

\noindent{\bf Example 3}. Let $H = \mathbb C ^2$.
Put $E_1 = \mathbb C (1,0)$,
$E_2 = \mathbb C (0,1)$ and $E_3 = \mathbb C (1,1)$.
Then $\mathcal S = (H;E_1,E_2,E_3)$ is indecomposable.
The system $\mathcal S$ is the lowest dimensional one 
among non-commutative indecomposable systems.

\bigskip
\noindent
{\bf Example 4}. Let $H = \mathbb C ^3$ and 
$\{a_1,a_2,a_3\}$
be a linearly independent subset of $H$.  Put $E_1 = \mathbb C a_1$,
$E_2 = \mathbb C a_2$ and $E_3 = \mathbb C a_3$.
Then $\mathcal S = (H;E_1,E_2,E_3)$ is decomposable.  In fact, let
$H_1 = E_1 \vee E_2 \not= 0$ and $H_2 = E_3 \not= 0$.  Then
$H_1 + H_2 = H$, $H_1 \cap H_2 = 0$ and
$E_i = E_i \cap H_1 + E_i \cap H_2$,
for $i=1,2,3$.

\bigskip
\noindent{\bf Example 5}.  Let $H = \mathbb C ^3$ and
$\{b_1,b_2,b_3,b_4\}$
be a  subset of $H$.
Put $E_i = \mathbb C b_i$ for $ i=1,\ldots, 4$.
Consider a system $\mathcal S = (H;E_1,E_2,E_3,E_4)$
of four subspaces.  Then the following are equivalent:
\begin{enumerate}
\item $\mathcal S$ is indecomposable.
\item Any three vectors of $\{b_1,b_2,b_3,b_4\}$ is linearly independent.
\item The set  $\{b_1,b_2,b_3\}$ is linearly independent and
$b_4 = \lambda _1 b_1 + \lambda _2 b_2 + \lambda _3 b_3$ for some
scalars $\lambda _i \not= 0 \ (i = 1,2,3)$.
\end{enumerate} 
Assume that $\{u_1,u_2,u_3,u_4\} \subset H$ and
$\{v_1,v_2,v_3,v_4\} \subset H$ satisfy the above condition (2).
Then $\mathcal S = (H;\mathbb C u_1,\mathbb C u_2,\mathbb C u_3,\mathbb C u_4)$
and  $\mathcal T = (H;\mathbb C v_1,\mathbb C v_2,\mathbb C v_3,\mathbb C v_4)$
are isomorphic. 

\bigskip
\noindent{\bf Example 6}. Let $H = \mathbb C ^3$.  Put
$E_1 = \mathbb C \oplus \mathbb C \oplus 0$, $E_2 = \mathbb C(1,1,1)$ and
$E_3 = \mathbb C(1,2,3)$.  Then a system $\mathcal S = (H;E_1,E_2,E_3)$ is
decomposable.  In fact, let $E'_1 = (E_2 \vee E_3) \cap E_1$ and
$H_1 = E_1 \cap (E'_1)^{\perp} \not= 0$.  Let $H_2 = E_2 \vee E_3 \not= 0$.
Then $H_1 + H_2 = H$, $H_1 \cap H_2 = 0$ and 
$E_i = E_i \cap H_1 + E_i \cap H_2$ for $i=1,2,3$.

\bigskip
\noindent{\bf Example 7}. Let $H = \mathbb C ^3$.  Put
$E_1 = \mathbb C \oplus \mathbb C \oplus 0$, $E_2 = \mathbb C(0,0,1)$,
$E_3 = \mathbb C(0,1,1)$ and
$E_4 = \mathbb C(1,0,1)$.  
Then a system $\mathcal S_7 = (H;E_1,E_2,E_3,E_4)$
of four subspaces is indecomposable.

\bigskip
\noindent{\bf Example 8}. Let $H = \mathbb C ^3$.  Put
$E_1 = \mathbb C \oplus \mathbb C \oplus 0$, $E_2 = \mathbb C(0,0,1)$,
$E_3 =  \mathbb C (1,0,0) + \mathbb C(0,1,1)$ and 
$E_4 = \mathbb C(1,0,1)$.  
Then a system $\mathcal S_8 = (H;E_1,E_2,E_3,E_4)$
of four subspaces is indecomposable.

\bigskip
\noindent{\bf Example 9}. Let $H = \mathbb C ^3$.  Put
$E_1 = \mathbb C \oplus \mathbb C \oplus 0$, $E_2 = \mathbb C(0,0,1)$,
$E_3 =  \mathbb C (1,0,0) + \mathbb C(0,1,1)$ and
$E_4 = \mathbb C(1,0,1) + \mathbb C(0,1,0)$.  Then a system
$\mathcal S_9 = (H;E_1,E_2,E_3,E_4)$ of
four subspaces is indecomposable.

\bigskip
\noindent{\bf Example 10}. Let $H = \mathbb C ^3$.  Put
$E_1 = \mathbb C(1,0,0) + \mathbb C(0,1,0)$,
$E_2 = \mathbb C(0,1,0) + \mathbb C(0,0,1)$
$E_3 =  \mathbb C (1,0,0) + \mathbb C(0,1,1)$ and
$E_4 = \mathbb C(0,0,1) + \mathbb C(1,1,0)$.  Then a system
$\mathcal S_{10} = (H;E_1,E_2,E_3,E_4)$ of
four subspaces is indecomposable.

\bigskip
\noindent{\bf Remark} 
Any two of the above indecomposable systems
$\mathcal S_7, \ldots, \mathcal S_{10} $ of four subspaces 
are not isomorphic each other.

\bigskip
\noindent{\bf Example 11}. 
Let $K = \ell^2(\mathbb N)$ and $H = K \oplus K$.
Consider a unilateral shift $S : K \rightarrow K$.  Let $E_1 = K \oplus 0$,
$E_2 = 0 \oplus K$, $E_3 = \{(x,Sx) \in H ; x \in K \}$ and
$E_4 = \{(x,x) \in H ; x \in K \}$.
Then a system $\mathcal S_{11} = (H;E_1,E_2,E_3,E_4)$ of 
four subspaces in $H$ is indecomposable. In fact, let 
$R$ be an idempotent which commutes with $S$. Then  
$R$ is a lower triangular Toeplitz matrix.  Since 
$R$ is an idempotent, $R = 0$ or $R = I$.   

\bigskip

Recall that Halmos initiated the study of transitive lattices.  
A complete lattice of closed subspaces of a Hilbert space $H$
containing $0$ and $H$ is called 
{\it transitive} if every bounded operator on $H$ leaving each subspace 
invariant is a scalar multiple of the identity. Halmos gave an example of 
transitive lattice consisting of seven subspaces in \cite{Ha2}. 
Harison-Radjavi-Rosenthal \cite{HRR} constructed a transitive lattice  
consisting of six subspaces using the graph of an unbounded operator.   
Any finite transitive lattice which consists of  $n$ subspaces  
gives an indecomposable system of $n$-$2$ subspaces  but the converse is 
not true. Following the study of transitive lattices, we shall introduce 
the notion of  transitive system. 

\bigskip
\noindent{\bf Definition.} Let $\mathcal S = (H;E_1,\ldots ,E_n)$ 
be a system of $n$ subspaces in a Hilbert space $H$. Then we say that 
$\mathcal S$ is {\it transitive} if 
$End(\mathcal S) = {\mathbb C}I_H$.  Recall that  $\mathcal S$ is 
indecomposable if and only if $Idem(\mathcal S) = \{0,I\}$. Hence 
if $\mathcal S$ is transitive, then $\mathcal S$ is indecomposable.
But the converse is not true. In fact the system $\mathcal S_{11}$ 
as above is indecomposable but is not transitve, because 
$End(\mathcal S)$ contains $S \oplus S$.

\bigskip
\noindent{\bf Example 12}.(Harrison-Radjavi-Rosenthal \cite{HRR})
Let $K = \ell^2(\mathbb Z)$ and $H = K \oplus K$.
Consider a sequence $(\alpha _n)_n$ given by $\alpha _n = 1$ for
$n \le 0$ and $\alpha _n = exp((-1)^nn!)$ for $n > 1$.
Consider a bilateral weighted shift $S : \mathcal D_T \rightarrow K$ such
that $T(x_n)_n = (\alpha _{n-1}x_{n-1})_n$ with the domain
$\mathcal D_T = \{(x_n)_n \in  \ell^2(\mathbb Z) ; 
\sum_n |\alpha _nx_n|^2 < \infty \}$.
Let $E_1 = K \oplus 0$,
$E_2 = 0 \oplus K$, $
E_3 = \{(x,Tx) \in H ; x \in \mathcal D_T \}$ and
$E_4 = \{(x,x) \in H ; x \in K \}$.
Harrison, Radjavi and Rosental showed that 
$\{0, H, E_1,E_2,E_3,E_4\}$ is a transitive lattice. 
Hence the  system $\mathcal S = (H;E_1,E_2,E_3,E_4)$
of four subspaces in H is transitive and in particular 
indecomposable. 

\bigskip
Let $\mathcal S = (H;E_1,\ldots ,E_n)$ be a system
of $n$ subspaces  in a finite-dimensional vector  space $H$. 
Gelfand and Ponomarev \cite{GP} introduced the conjugate system 
$\mathcal S^* = (H^*;E_1',\ldots ,E_n')$, 
where $E_i' = \{f \in H^* ; f(x) = 0 
{\text \ for \ all \  } x \in E_i \}$.  In our setting of Hilbert 
spaces,  their conjugate system $\mathcal S^*$ could  be 
replaced by the following orthogonal complement.

\bigskip
\noindent{\bf Definition}. Let $\mathcal S = (H;E_1,\ldots ,E_n)$ be a system
of $n$ subspaces  in a Hilbert space $H$. Then the orthogonal complement of
$\mathcal S$, denoted by $\mathcal S^{\perp}$, is defined by
$\mathcal S^{\perp} = (H;E_1 ^{\perp},\ldots ,E_n ^{\perp})$.
Let ${\mathcal T} = (K;F_1, \dots , F_n)$  
be  another system of $n$ subspaces in a Hilbert space $K$ and  
$\varphi : {\mathcal S} \rightarrow {\mathcal T}$ be a 
homomorphism.  We define a homomorphism 
$\varphi^* : {\mathcal T}^{\perp} \rightarrow 
{\mathcal S}^{\perp}$ by $\varphi^*: K \rightarrow H$. In fact,  
$\varphi^*(F_i^{\perp}) \subset E_i ^{\perp}$, because 
$\varphi (E_i) \subset F_i$. 

We denote by ${\mathcal Sys}^n$ 
the category of the systems of $n$ subspaces in Hilbert spaces
and homomorphisms. Then we can introduce a contravariant 
functor $\Phi ^{\perp}: {\mathcal Sys}^n \rightarrow 
{\mathcal Sys}^n$ by 
\[
\Phi ^{\perp}({\mathcal S}) = {\mathcal S}^{\perp} 
\text{ and } \Phi ^{\perp}(\varphi) =  \varphi^* . 
\]

\begin{prop} Let $H$ be a Hilbert space and $\mathcal S
= (H;E_1,\ldots ,E_n)$ a system of $n$ subspaces in $H$.  
Then $\mathcal S$
is indecomposable if and only if $\mathcal S^{\perp}$ is indecomposable.
\end{prop}
\begin{proof}
If $\mathcal S$ is decomposable, then there exists an idempotent 
$R \in End(\mathcal S)$ with $R \not= 0$ and $R \not= I_H$. Since 
$R(E_i) \subset E_i$, we have 
$R^*(E_i^{\perp}) \subset E_i^{\perp}$. Thus 
$R^* \in End(\mathcal S^{\perp})$ is an idempotent with 
$R^* \not= 0$ and $R^* \not= I_H$, that is, $\mathcal S^{\perp}$ 
is decomposable. 
This implies the desired 
conclusion. 
\end{proof}

Similarly we have a same fact for transitive systems.

\begin{prop} Let $H$ be a Hilbert space and $\mathcal S
= (H;E_1,\ldots ,E_n)$ a system of $n$ subspaces in $H$.  
Then $\mathcal S$
is transitive  if and only if $\mathcal S^{\perp}$ is transitive .
\end{prop}

\section{indecomposable systems of one subspace}
It is easy to see the case of indecomposable systems 
of one subspace even in an infinite-dimensional Hilbert space. 

\begin{prop}
Let H be a Hilbert space and 
$\mathcal {S}=(H;E)$ a system of one subspace. 
Then $\mathcal{S}=(H;E)$ is indecomposable 
if and only if 
$\mathcal {S} \cong (\mathbb C;0)$ or 
$\mathcal {S} \cong (\mathbb C;\mathbb C)$.
\end{prop}
\begin{proof} If $E \not=0$ and  $E \not = H$, then 
$\mathcal{S} = (E;E) \oplus (E^{\perp};0)$ gives a 
non-trivial decomposition. Assume that $\mathcal{S}$ 
is indecomposable.  Then $E = 0$ or $E = H$.  
Suppose we had $dim H\geq 2$, then there exist 
non-zero closed subspaces
$H_{1}$ and $ H_{2}$ such that $H=H_{1}+H_{2}$ and 
$H_{1}\cap H_{2}= 0$. This gives a non-trivial decompositon 
of $\mathcal{S}$. The contradiction shows that $dim H = 1$. 
Hence $\mathcal {S} \cong (\mathbb C;0)$ or 
$\mathcal {S} \cong (\mathbb C;\mathbb C)$. The 
converse is trivial.
\end{proof}

Let $\mathcal {S}=(H;E)$ and $\mathcal {S}'=(H';E')$ be 
two systems of one subspace.  Then $\mathcal {S}$ and 
$\mathcal {S}'$ are isomorphic if and only if 
$\dim E = \dim E'$ and $\codim E = \codim E'$.

\section{indecomposable systems of two subspaces}
It is a well known fact that the relative position of two subspaces
$E_1$ and $E_2$ in a Hilbert space $H$  can be described 
completely up to unitary equivalence  as in Araki \cite{Ar}, Dixmier \cite{D} 
and Halmos \cite{Ha1}.  The Hilbert space $H$ is the direct sum of five
subspaces:
$$
H = (E_1\cap E_2) \oplus (\text{the rest}) \oplus (E_1\cap E_2^{\perp})
\oplus (E_1^{\perp} \cap E_2) \oplus (E_1^{\perp} \cap E_2^{\perp}).
$$
In the rest part,  $E_1$ and $E_2$ are in generic position and the
relative position is described only by \lq\lq the angles" between them.
In fact the rest part is written as $K \oplus K$ for some subspace
$K$ and there exist two positive operators $c,s \in B(K)$ 
with null kernels with $c^2 + s^2 = 1$ such that 
\[
E_{1}=(E_{1}\cap E_{2})\oplus 
\Im \left( 
\begin{array}
{@{\,}cccc@{\,}} 1&0\\ 0&0
\end{array}
\right) \oplus 
(E_{1}\cap E_{2}^{\perp})\oplus 0 \oplus 0,
\]
and
\[
E_{2}=(E_{1}\cap E_{2}) \oplus 
\Im\left(
\begin{array}
{@{\,}cccc@{\,}}c^{2}&cs\\ cs&s^{2}
\end{array}
\right)
\oplus 0 \oplus (E_{1}^{\perp}\cap E_{2}) \oplus 0.
\]
By the functional calculus, there exists a unique positive operator
$\theta$, called the angle operator, such that
$c = \cos \theta \ \  \mbox{and} \ \ s = \sin \theta$ 
with $0 \leq \theta \leq \frac{\pi}{2}$.

\begin{prop} Let $\mathcal S = (H;E_1,E_2)$ be a system of two subspaces in 
a Hilbert space $H$.
Then $\mathcal S$ is indecomposable if and only if $\mathcal S$ is 
isomorphic to one of the
following four commutative systems:
$$
\mathcal S_1 = (\mathbb C; \mathbb C, 0), \ \  \mathcal S_2 = (\mathbb C; 0, \mathbb C), \\
\mathcal S_3 = (\mathbb C; \mathbb C, \mathbb C), \ \  \mathcal S_4 = (\mathbb C; 0, 0).
$$
\label{prop:two-subspaces}
\end{prop}
\begin{proof} Let $e_i \in B(H)$ be the projection of $H$ onto 
$E_i$, $i = 1,2$ with the canonical decomposition as above. 
Suppose that $\dim K \geq 2$. Then there exists 
a projection $p \in B(K)$ with $0 \not= p \not= I_K$ satisfying 
$p$ commutes with $c$ and $s$. Let 
$H_1 := \Im (p\oplus p) \subset K \oplus K$ and 
$H_2 := H_1^{\perp} \cap H$. Let $p_1 \in B(H)$ be the projection 
of $H$ onto $H_1$. 
Since non-trivial projection $p_1$ commute with $e_1$ and $e_2$, 
$\mathcal S$ is decomposable by Lemma \ref{lemma:endo-decompose}.  
Therefore if $\mathcal S$ is indecomposable, 
then $\dim K \leq 1$ and only one of the five direct summands is non-zero. 
If the rest component were non-zero, then it is  isomorphic to a decomposable 
one as in Example 2. Thus the rest component does not appear. One of the 
other part is commutative.  Since $\mathcal S$ is indecomposable, 
$\mathcal S$ is one of $\mathcal S _1, \dots , {\mathcal S}_4$ by 
Proposition \ref{prop:commutative}. The converse is clear. 
\end{proof}

\section{some properties of indecomposable systems of $n$-subspaces}

Let ${\mathcal S} = (H;E_1, \dots , E_n)$ be a system of 
$n$ subspaces 
in a Hilbert space.  We denote by $\vee _{i = 1}^n E_i$ 
the closed subspace spanned by $E_1,\dots,E_n$. 
If ${\mathcal S}$ is indecomposable and 
$\dim H \geq 2$, then it is easy to see that 
\[
\bigcap _{i = 1}^n E_i = 0  \text{ and  } 
\bigvee _{i = 1}^n E_i = H.  
\]
In fact, on the contrary suppose that 
$M:= \cap _{i = 1}^n E_i \not= 0$. 
We choose a one-dimensional subspace $F \subset M$.
Since $\dim H \geq 2$, the orthogonal decomposition 
$H = F \oplus F^{\perp}$  of the Hilbert space $H$ gives 
a non-trivial decomposition of the system 
${\mathcal S}$.   This contradicts to that ${\mathcal S}$ is indecomposable. 
Hence we have $\cap _{i = 1}^n E_i = 0$. Since the orthogonal 
complement ${\mathcal S}^{\perp}$ is also indecomposable,  we have 
$\vee _{i = 1}^n E_i = (\cap _{i = 1}^n E_i^{\perp})^{\perp} = H$. 
But we can say more as follows:

\begin{prop}
Let ${\mathcal S} = (H;E_1, \dots , E_n)$ be a system of $n$ subspaces 
in a Hilbert space.  If ${\mathcal S}$ is indecomposable and 
$\dim H \geq 2$, then for any distinct $n$-$1$ subspaces 
$E_{i_1}, \dots, E_{i_{n-1}}$, we have that
\[
\bigcap _{k = 1}^{n-1} E_{i_k} = 0  \text{ and  } 
\bigvee _{k = 1}^{n-1} E_{i_k} = H.  
\]
\label{prop:n-1}
\end{prop}
\begin{proof}

We may and do assume that $E_{i_1} = E_1, E_{i_2}= E_2,
\dots, E_{i_{n-1}} = E_{n-1}$.  On the contrary suppose that 
$M:= \cap _{i = 1}^{n-1} E_i \not= 0$. 
Since $\dim H \geq 2$, we can choose a one-dimensional subspace 
$F \subset M$. Consider two subspaces $F$ and $E_n$ in $H$. We 
have the following canonical decomposition into five parts: 
\[
F=(F\cap E_{n})\oplus 
\Im
\left(
\begin{array}{@{\,}cccc@{\,}}
1&0\\
0&0
\end{array}
\right)
\oplus 
(F\cap E_{n}^{\perp})
\oplus 0 \oplus 0,
\]
\[
E_{n}=(F\cap E_{n})
\oplus 
\Im \left(
\begin{array}{@{\,}cccc@{\,}}
c^{2} & cs \\
 cs&s^{2}
\end{array}
\right)
\oplus 0
\oplus 
(F^{\perp}\cap E_{n}) \oplus 0.
\]
We denote by  $K\oplus K$
the underlying subspace of the part in generic position. 

\noindent
(i)(the case that $K = 0)$: Since $F\cap E_{n} = \cap _{i = 1}^n E_i = 0$, 
we have $F = F\cap E_{n}^{\perp}$, so that $F \subset E_n^{\perp}$. 
Let $e_i$ and $f$ be the projections of $H$  onto $E_i$ and $F$ 
respectively. Then 
$f$ commutes with each $e_i$.  Therefore 
the orthogonal decomposition $H = F \oplus F^{\perp}$ of  $H$ 
gives a non-trivial decomposition of the system 
${\mathcal S}$.  This contradicts to that ${\mathcal S}$ is indecomposable.  
Hence $M = \cap _{i = 1}^{n-1} E_i = 0$. 

\noindent
(ii)(the case that $K \not= 0$): Since $F$ is one-dimensional,
\[
K\oplus 0 +
\Im \left(
\begin{array}{@{\,}cccc@{\,}}
c^{2} & cs \\
 cs&s^{2}
\end{array}
\right)
= K \oplus K
\]
and 
\[
(K\oplus 0)
\cap  
\Im \left(
\begin{array}{@{\,}cccc@{\,}}
c^{2} & cs \\
 cs&s^{2}
\end{array}
\right)
=0.
\]
Then there exists an invertible
operator $T\in B(K\oplus K)$
such that
$T(K\oplus 0)=K\oplus 0,$
and
$T(\Im
\left(
\begin{array}{@{\,}cccc@{\,}}
c^{2} & cs \\
 cs&s^{2}
\end{array}
\right))
=0\oplus K$.

We define an invertible operator 
$\varphi:=I\oplus T\oplus I\oplus I\oplus I \in B(H)$.
Let $E_{i}^{\prime}:= \varphi(E_i)$ for $i = 1,\dots ,n$. 
Since $\mathcal S$ is indecomposabe, a new system 
 $\mathcal S^{\prime}:=(H;E_{1}^{\prime},\dots, E_{n}^{\prime})$
is indecomposable. 
Since $F = \varphi(F)$,  $F \subset \cap _{i = 1}^{n-1} E_i^{\prime}$ 
and $F$ is orthogonal to $E_n^{\prime}$. 
Let $e_i^{\prime}$ and $f$ be the projections of $H$  onto 
$E_i^{\prime}$ and $F$. Then 
$f$ commutes with each $e_i^{\prime}$.  Therefore 
the orthogonal decomposition $H = F \oplus F^{\perp}$  of $H$ 
gives a non-trivial decomposition of the system  
${\mathcal S}^{\prime}$.  This contradicts to that 
${\mathcal S}^{\prime}$ is indecomposable.  
Hence $M = \cap _{i = 1}^{n-1} E_i = 0$

Since  the orthogonal complement ${\mathcal S}^{\perp}$ is also 
indecomposable, we also have $\vee _{k = 1}^{n-1} E_{i_k} = H$. 
\end{proof}

\begin{cor}
Let ${\mathcal S} = (H;E_1, \dots , E_n)$  a system of $n$ subspaces 
in a Hilbert space.  If ${\mathcal S}$ is indecomposable and 
$H$ is infinite-dimensional,  then 
$\ ^\#\{i ; E_i \text{ is finite dimensional } \} \leq n-2$. 
\end{cor}
\begin{proof} 
On the contrary, suppose that there were 
distinct $n$-$1$  finite-dimensional subspaces 
$E_{i_1}, \dots, E_{i_{n-1}}$. Then 
$H = \bigvee _{k = 1}^{n-1} E_{i_k}$ is also 
finite-dimensional. This is a contradiction.
\end{proof}

\section{indecomposable systems of three subspaces}

Gelfand and Ponomarev (\cite{GP}) claimed that 
there exist only nine, finite-dimensional, indecomposable 
systems of three subspaces.  We shall include a direct proof of it. 
We do not know whether there exists an infinite-dimensional 
transitive systems of three subspaces. In fact it is still 
an unsolved problem whether there exists a transitive lattice 
consisting of five elements in an infinite-dimensional Hilbert space. 
Therefore it is worth while investigating the existence of 
infinite-dimensional indecomposable systems of three subspaces. 
 
\begin{prop}
Let $\mathcal{S}=
(H;E_{1},E_{2},E_{3})$
be an indecomposable system of three subspaces. 
If $H$ is infinite dimensional, then 
$E_{i} \not= 0$  and $E_{i}\not=H$ for $i=1,2,3.$
\end{prop}
\begin{proof}
On the contrary suppose that $E_{1}=0$.
Then $\mathcal S^{\prime}=(H;E_{2},E_{3})$
is an indecomposable system of two subspaces. 
Hence by Proposition \ref{prop:two-subspaces}, $H$ is finite 
dimensional.  This is a contradiction. Hence $E_{1}\not=0$. 
Similary $E_{i} \not= 0$  and $E_{i}\not=H$ for $i=1,2,3.$ 
\end{proof}

\begin{thm}
Let $\mathcal{S}=(H;E_{1},E_{2},E_{3})$
be an indecomposable system of three subspaces in a Hilbert space $H$.
Then the following hold. 

\noindent
(1)If $H$ is infinite-dimensional, then for any $i \not= j$, 
$E_i \cap E_j = 0$ and $E_i + E_j$ is a non-closed dense subspace of $H$. 
In particular each $E_i$ is infinite-dimensional. 

\noindent
(2)\cite{GP} If $H$ is finite-dimensional, 
then $\mathcal S$  is isomorphic to 
one of the following eight commutaitve systems 
${\mathcal S}_1, \dots, {\mathcal S}_8$  
 and one non-commutative system ${\mathcal S}_9$:
\[
{\mathcal S}_1 = (\mathbb C;0,0,0),
\ \ {\mathcal S}_2 = (\mathbb C;\mathbb C,0,0),
\ \ {\mathcal S}_3 = (\mathbb C; 0,\mathbb C,0),
\]
\[
{\mathcal S}_4 = (\mathbb C; 0,0, \mathbb C), 
\ \ {\mathcal S}_5 = (\mathbb C;\mathbb C,\mathbb C,0),
\ \ {\mathcal S}_6 = (\mathbb C;\mathbb C,0,\mathbb C),
\]
\[
{\mathcal S}_7 = (\mathbb C; 0,\mathbb C,\mathbb C),
\ \ {\mathcal S}_8 = (\mathbb C; \mathbb C , \mathbb C,\mathbb C), 
\ \ {\mathcal S}_9 = (\mathbb C^2; \mathbb C (1,0), 
\mathbb C (0,1),\mathbb C (1,1)).
\]
\end{thm}
\begin{proof}
If $\dim H = 1$, then ${\mathcal S}$ is commutative. Hence 
if ${\mathcal S}$ is isomorphic to one of 
${\mathcal S}_1, \dots, {\mathcal S}_8$. Therefore  
we may assume that ${\mathcal S}$ is indecomposable and 
$\dim H \geq 2$. Then, by Proposition \ref{prop:n-1}, 
for any $i \not= j$, 
$E_i \cap E_j = 0$ and $E_i + E_j$ is a dense subspace of $H$. 
We claim  that if $E_1 + E_2 = H$, then 
$H$ is finite-dimensional and ${\mathcal S}$ is isomorphic to 
${\mathcal S}_9$. It is enough to show the claim to prove the 
theorem. In fact, assume that the claim holds.  
(1)If $H$ is infinite-dimensional, then $E_1 + E_2$ is not closed. 
Similarly  for any $i \not= j$, $E_i + E_j$ is not closed. 
(2)If $H$ is finite-dimensional, then $E_1 + E_2 = H$. Thus 
${\mathcal S}$ is isomorphic to ${\mathcal S}_9$. 
We shall show the claim.  Since $E_1 \cap E_2 = 0$ and 
$E_1 + E_2 = H$, there exists $T\in B(H)^{-1}$
such that $T(E_{1})=E_{1}$ and $T(E_{2})=E_{1}^{\perp}$. 
Therefore we may assume that $E_{2} =E_{1}^{\perp}$ to show 
the claim.  Considering the canonical decomposition for  
two subspaces $E_1$ and $E_3$, we have the following descripton 
of three subspaces: 
\[
E_{1}=(E_{1}\cap E_{3} )
\oplus 
\Im\left(
\begin{array}{@{\,}cccc@{\,}}
1&0\\
0&0
\end{array}
\right)
\oplus 
(E_{1}\cap E_{3}^{\perp})\oplus 0 \oplus 0  ,
\]
\[
E_{3}=(E_{1}\cap E_{3})\oplus 
\Im
\left(
\begin{array}{@{\,}cccc@{\,}}
c^{2} & cs \\
 cs&s^{2}
\end{array}
\right)
\oplus 0
\oplus (E_{1}^{\perp}\cap E_{3}) \oplus 0 , 
\]
\[
E_2 = E_{1}^{\perp}
= 0 \oplus \Im
\left(
\begin{array}{@{\,}cccc@{\,}}
1&0\\
0&0
\end{array}
\right)
\oplus 0 \oplus 
(E_{1}^{\perp}\cap E_{3})\oplus (E_{1}^{\perp}\cap E_{3}^{\perp}), 
\]
where the underlying Hilbert space $H$ is decomposed into 
five parts 
\[
H = (E_{1}\cap E_{3}) \oplus (K \oplus K) \oplus 
(E_{1}\cap E_{3}^{\perp}) \oplus 
(E_{1}^{\perp}\cap E_{3}) \oplus (E_{1}^{\perp}\cap E_{3}^{\perp}) .
\]
If two parts of the above five parts were non-zero, 
then ${\mathcal S}$ can be decomposed non-trivially. 
This contradicts to that ${\mathcal S}$ is indecomposable. 
Hence only one of the above five parts is non-zero. 
If the part $K \oplus K = 0$, then ${\mathcal S}$ is commutative. 
Since ${\mathcal S}$ is indecomposable, $\dim H = 1$.  This contradicts 
to that $\dim H \geq 2$.  Hence the only the part $K \oplus K \not= 0$. 
If $\dim K = 1$, then it is clear that ${\mathcal S}$ is 
isomorphic to ${\mathcal S}_9$.  If $\dim K  \geq 2$, then there exists 
a projection $p \in B(K)$ with $0 \not= p \not= I_K$ satisfying 
$p$ commute with $c$ and $s$. Let 
$H_1 := \Im (p\oplus p) \subset K \oplus K = H$ and 
$H_2 := H_1^{\perp} \cap H$. Let $p_1, e_1, e_2, e_3 \in B(H)$ 
be the projections of $H$ onto $H_1, E_1, E_2, E_3$ respectively. 
Since non-trivial projection $p_1$ commute with $e_1$, $e_2$ and $e_3$, 
$\mathcal S$ is decomposable by Lemma \ref{lemma:endo-decompose}. 
This is a contradiciton. Hence the case that $\dim K  \geq 2$ does not 
occur.  We have shown the claim.    
\end{proof}

\section{operator systems }
We can associate a system of four subspaces for any operator.

\noindent
{\bf Definition.} (bounded operator system) 
We say that a system $\mathcal{S}=(H;E_1,E_2,E_3,E_4)$ of 
four subspaces is a {\it bounded operator system} if  
there exist a Hilbert space $K_1,K_2$ and bounded 
operators $T: K_1 \rightarrow K_2 $, 
$S:K_2 \rightarrow K_1$ such that $H=K_1\oplus K_2$ and 
\[
E_{1}=K_1\oplus 0, \ 
E_{2}=0\oplus K_2, 
\]
\[
E_{3}=\{(x,Tx); x\in K_1\}, \ 
E_{4}=\{(Sy,y); y\in K_2\}. 
\]
We denote by $\mathcal{S}_{T,S}$ the above operator system $\mathcal{S}$. 
We often identify $E_1$ with $K_1$ and $E_2$ with $K_2$.   
In particular we associate an operator system 
$\mathcal{S}_T := \mathcal{S}_{T,I} = (H;E_1,E_2,E_3,E_4)$ 
for any single operator $T \in B(K)$ such that  $H = K \oplus K$ and 
\[
E_{1}=K\oplus 0,
E_{2}=0\oplus K, 
E_{3}=\{(x,Tx); x\in K\},
E_{4}=\{(y,y); y\in K\}. 
\]

We shall study a relation between the system $\mathcal{S}_T$ 
of four subspaces and a single operator $T$. 
  
\begin{prop} Let $\mathcal{S}_{T,S} = $ $(H;E_1,E_2,E_3,E_4)$ be a 
bounded operator system associated with 
$T: K_1 \rightarrow K_2 $ and $S:K_2 \rightarrow K_1$.
Then

\begin{align*}
End (\mathcal{S}_{T,S}) & = \{ A_1 \oplus A_2 \in B(H) ; 
A_1 \in B(K_1), A_2 \in B(K_2), \\
& \ A_1S = SA_2, \ A_2T = TA_1 \}, {\text and}
\end{align*}

\begin{align*}
Idem (\mathcal{S}_{T,S}) & = \{ A_1 \oplus A_2 \in B(H) ; 
A_1 \in B(K_1), A_2 \in B(K_2), \\
& \ A_1S = SA_2, \ A_2T = TA_1, 
A_1^2 = A_1, \ A_2^2 = A_2 \}
\end{align*}

\end{prop}

\begin{proof}Let $A \in End (\mathcal{S})$. Since 
$A(E_1) \subset E_1$ and  $A(E_2) \subset E_2$, 
we have $A = A_1 \oplus A_2$ 
for some $A_1 \in B(K_1), A_2 \in B(K_2)$. Since $A(E_3) \subset E_3$, 
for any $x \in K_1$, $(A_1 \oplus A_2)(x,Tx) \in E_3$. Thus 
$(A_1x,A_2Tx) = (y,Ty)$ for some $y \in K_2$. Therefore 
$A_2Tx = TA_1x$. Thus $A_2T = TA_1$. 
Similarly  $A(E_3) \subset E_3$ implies  $A_1S = SA_2$. 
The converse is clear. We get the equality for 
$Idem(\mathcal{S})$ immediately. 
\end{proof}

\begin{cor}
Let $\mathcal{S}_T = (H;E_1,E_2,E_3,E_4)$ be a  
bounded operator system associated with a single operator
 $T \in B(K)$. Then
\[ 
End (\mathcal{S}_T) = \{ B \oplus B \in B(H) ; 
B \in B(K), \ BT = TB \}, {\text and}
\]  
\[
Idem(\mathcal{S}_T) = \{ B \oplus B \in B(H) ; 
B \in B(K), \ BT = TB, B^2 = B \}. 
\] 
\end{cor}

\noindent
{\bf Definition.} Recall that a bounded operator $T$ on a 
Hilbert space $K$ is called {\it strongly irreducible} 
if there do not exist two non-trivial subspaces $M \subset K$  
and $N \subset K$ such that $T(M) \subset M$, $T(N) \subset N$,  
$M \cap N = 0$ and $M + N = K$. We also see that $T$ is 
strongly irreducible if and only if there does not exist 
any non-trivial idempotent $P$ such that $PT = TP$. See a 
monograph \cite{JW} by Jiang and Wang. 

\begin{cor}
Let $\mathcal{S}_T = (H;E_1,E_2,E_3,E_4)$ be a  
bounded operator system associated with a single operator
 $T \in B(K)$. Then $\mathcal{S}_T$ is indecomposable 
if and only if $T$ is strongly irreducible. 
\end{cor} 

\noindent
{\bf Example. } Let $K = \ell ^2 (\mathbb{N})$ and $S \in B(K)$ be 
the unilateral shift.  Let $P \in B(K)$ be an idempotent 
which commutes with $S$. Then $P$ is a lower triangular  Toeplitz 
matrix. Since $P$ is an idempotent, we have  $P = 0$ or $P = I$ as in 
Lemma \ref{lemma:Toeplitz-idempotent}. Thus 
$S$ is strongly irreducible, as already known, for example, in \cite{JW}, and 
$\mathcal{S}_S$ is indecomposable. 

\begin{cor}
Let $\mathcal{S}_T = (H;E_1,E_2,E_3,E_4)$ be a 
bounded operator system associated with a single operator
 $T \in B(K)$. If $\mathcal{S}_T$ is decomposable,  
then $T$ has a non-trivial invariant subspace. 
\end{cor} 
\begin{proof}Let $\mathcal{S}_T$ be decomposable. Then 
there exists a non-trivial idempotent $P$ such that $PT = TP$.
Then $\Im P$ is a non-trivial invariant subspace.  
\end{proof}

\begin{prop}
Let $\mathcal{S}_T = (H;E_1,E_2,E_3,E_4)$ and 
$\mathcal{S}_{T'}  = (H';E_1',E_2',$ $E_3',E_4')$
be bounded operator systems associated with operators
$T \in B(K)$ and  $T' \in B(K')$.  Then 
$\mathcal{S}_T$ and $\mathcal{S}_{T'}$ are isomorphic 
if and only if $T$ and $T'$ are similar.    
\end{prop}
\begin{proof} 
Assume that $\mathcal{S}_T$ and $\mathcal{S}_{T'}$
are isomorphic.  Then there exists a bounded invertible 
operator $A: H \rightarrow H'$ with $A(E_i) = E_i'$ for 
$i = 1,2,3,4$.  Since $A(E_i) = E_i'$ for $i = 1,2,4$, we have 
$A = B \oplus B$ for some invertible operator 
$B: K \rightarrow K'$.  And $A(E_3) \subset E_3$ implies that 
$BT = T'B$, that is, $T$ and $T'$ are similar. The converse 
is clear.
\end{proof}

\noindent
{\bf Remark.} The above proposition shows that the classification 
of systems of four subspaces contains the classification of 
operators up to similarity in a certain sense. 

\noindent
{\bf Example}.(an uncountable family  of 
indecomposable systems of four subspaces)   
Let $K=\ell^{2}(\mathbb{N})$ and 
 $H=K\oplus K$. Consider a unilateral shift
$S:K\to K$. For a parameter $\alpha\in \mathbb {C}$, 
let $E_{1}=K\oplus 0, E_{2}= 0\oplus  K, 
E_{3}=\{(x,(S+\alpha I)x)\vert x\in K\}$ and 
$E_{4}=\{(x,x)\vert x\in K\}$. 
Then the  system $\mathcal {S}_{\alpha}=(H;E_1,E_2,E_3,E_4)$
of four subspaces  are indecomposable.
If $\alpha\not= \beta$, then $\mathcal {S}_{\alpha}$ and
$\mathcal {S}_{\beta}$ are not isomorphic, because 
the spectra $\sigma (S+\alpha)\not=
\sigma(S+\beta)$ and 
$S+\alpha I$ and $S + \beta I$ are not similar. 
Thus we can easily construct an uncountable family 
$(\mathcal {S}_{\alpha})_{\alpha \in \mathbb {C}}$  of 
indecomposable systems of four subspaces. 

As the single operator case,  we also obtain the following: 
 
\begin{prop}
Let $\mathcal{S}_{T,S} = (H;E_1,E_2,E_3,E_4)$ and 
$\mathcal{S}_{T',S'} = (H';E_1',$ $E_2',E_3',E_4')$
be bounded operator systems associated with operators
$S \in B(K_2,K_1), T \in B(K_1,K_2)$,
$S' \in B(K_2',K_1'), T' \in B(K_1',K_2')$ .  Then 
$\mathcal{S}_{T,S} $ and $\mathcal{S}_{T',S'}$ are isomorphic 
if and only if there exist bounded invertible operators 
$A_1: K_1 \rightarrow K_1'$ and $A_2: K_2 \rightarrow K_2'$ such that    
$A_1S = S'A_2$ and $A_2T = T'A_1$.
\end{prop}

\begin{prop}
Let $\mathcal{S}_{T,S} = (H;E_1,E_2,E_3,E_4)$
be a bounded operator system associated with operators
$S \in B(K_2,K_1), T \in B(K_1,K_2)$.　Then the orthogonal 
complement of the system $\mathcal{S}_{T,S}$ is isomorphic to 
another  bounded operator system up to a permutation of subspaces 
and given by 
\[
\mathcal{S}_{T,S}^{\perp} \cong \sigma_{1,2} \sigma_{3,4} \mathcal{S}_{-S^*,-T^*} ,
\]
where $\sigma_{i,j}$ is a transposition of $i$ and $j$. 
\end{prop}
\begin{proof}It is evident from the fact 
$\{(x,Tx) \in K_1\oplus K_2 ; x \in K_1 \}^{\perp} 
= \{(-T^*y,y) \in K_1\oplus K_2 ; y \in K_2 \}$ and etc. 
\end{proof}

\begin{prop}
Let $\mathcal{S}_{T,S} = (H;E_1,E_2,E_3,E_4)$
be a bounded operator system associated with operators
$S \in B(K_2,K_1), T \in B(K_1,K_2)$. If $T$ is invertible, 
then $\mathcal{S}_{T,S}$ is isomorphic to 
$\mathcal{S}_{I,TS}$. 
If $S$ is invertible, 
then $\mathcal{S}_{T,S}$ is isomorphic to 
$\mathcal{S}_{ST,I}$. 
\label{prop:ST}
\end{prop}
\begin{proof}Let $T$ be invertible.  Define an invertible operator 
$\varphi : K_1 \oplus K_2 \rightarrow K_2 \oplus K_2$ by 
$\varphi (x,y) = (Tx,y)$. Then 
$\varphi(E_1) = \varphi (K_1 \oplus 0) = K_2 \oplus 0$. 
$\varphi(E_2) = \varphi(0 \oplus K_2) = 0 \oplus K_2$. 
Since $\varphi (x,Tx) = (Tx,Tx)$, 
$\varphi(E_3) = \varphi (\graph T) = \{(y,y) ; y \in K_2 \}$. 
Because $\varphi (Sy,y) = (TSy, y)$,   
$\varphi(E_4) = \varphi (\cograph S) = \{(TSy,y) ; y \in K_2 \} 
= \cograph TS$. 
Hence $\mathcal{S}_{T,S}$ is isomorphic to 
$\mathcal{S}_{I,TS}$. If $S$ is invertible, use an invertible 
operator
$\psi : K_1 \oplus K_2 \rightarrow K_1 \oplus K_1$ defined by 
$\psi (x,y) = (x,Sy)$.  
\end{proof}

Bounded operator systems can be extended to (unbounded) 
closed operator systems. 

\noindent 
{\bf Definition.}(closed operator systems)
We say that a system $\mathcal{S}=$ $ (H;E_1,$ $E_2,E_3,E_4)$ of 
four subspaces is a {\it closed operator system} if  
there exist Hilbert spaces  $K_1,K_2$ and closed
operators $T: K_1 \supset D(T) \rightarrow K_2 $, 
$S:K_2 \supset D(S) \rightarrow K_1$ such that 
$H=K_1\oplus K_2$ and $E_{1}=K_1\oplus 0$, \ 
\[
E_{2}=0\oplus K_2, \ 
E_{3}=\{(x,Tx); x\in D(T) \},
E_{4}=\{(Sy,y); y\in D(S) \}. 
\]
We also denote by $\mathcal{S}_{T,S}$ the above operator system $\mathcal{S}$.

We shall give a characterization of (densely defined) 
closed operator systems. 

\begin{prop}
Let $\mathcal{S}=(H;E_1,E_2,E_3,E_4)$ be a system of  
four subspaces in a Hilbert space $H$. Then the following are 
equivalent:
\begin{enumerate}
\item $\mathcal{S}$ is isomorphic to a closed operator system 
$\mathcal{S}_{T,S}$ for some closed operators 
$T: E_1 \supset D(T) \rightarrow E_2$ and 
$S: E_2 \supset D(S) \rightarrow E_1$.  
\item $E_1 + E_2 = H$ and $E_i \cap E_j = 0$ for 
$(i,j) = (1,2), (2,3)$ and $(4,1)$.   
\end{enumerate}
Moreover if these conditions are satisfied, then 
$D(T):= E_1 \cap (E_3 + E_2)$ and 
$D(S):= E_2 \cap (E_4 + E_1)$.
\label{prop:operator-system} 
\end{prop}
\begin{proof}
(1)$\Rightarrow$(2): It is trivial. 
(2)$\Rightarrow$(1): By Lemma \ref{lemma:decompose}, 
we may assume that $E_2 = E_1^{\perp}$. 
Put $K_1 = E_1$ and $K_2 = E_2$.  Then $H = E_1 \oplus E_2$. Since 
$E_3 \cap E_2 = 0$, for any 
$x_1 \in E_1 \cap (E_3 + E_2)$, 
there exist unique $x_3 \in E_3$ and $x_2 \in E_2$ such that 
$x_1 = x_3 - x_2$.  Define a linear operator  
$T:E_1 \supset D(T) \rightarrow 
E_2$ by $Tx_1 = x_2$ with a domain $D(T):= E_1 \cap (E_3 + E_2)$. 
Since $E_1 + E_2 = H$, for any $x_3 \in E_3$ 
there exist $x_1 \in E_1$ and 
$x_2 \in E_2$ with $x_3 = x_1 + x_2$. This implies that 
$\graph T = E_3$.  Hence $T$ is a closed operator. Similarly 
there exists a closed operator 
$S: E_2 \supset D(S) \rightarrow E_1$ 
with a domain $D(S):= E_2 \cap (E_4 + E_1)$. 
\end{proof}

\begin{cor}
Let $\mathcal{S}=(H;E_1,E_2,E_3,E_4)$ be a system of  
four subspaces in a Hilbert space $H$. Then the following are 
equivalent:
\begin{enumerate}
\item $\mathcal{S}$ is isomorphic to a 
closed operator system 
$\mathcal{S}_{T,S}$ for some densely defined closed operators 
$T: E_1 \supset D(T) \rightarrow E_2$ and 
$S: E_2 \supset D(S) \rightarrow E_1$.  
\item $E_1 + E_2 = H$ and $E_i \cap E_j = 0$ for 
$(i,j) = (1,2), (2,3),(4,1)$, \\
$E_1 \cap (E_3 + E_2)$ is dense in $E_1$ , 
$E_2 \cap (E_4 + E_1)$ is dense in $E_2$   
\end{enumerate}
\end{cor}

We immediately have a characterization of bounded operator systems.   

\begin{cor}
Let $\mathcal{S}=(H;E_1,E_2,E_3,E_4)$ be a system of  
four subspaces in a Hilbert space $H$. Then the following are 
equivalent:
\begin{enumerate}
\item $\mathcal{S}$ is isomorphic to a bounded operator system.
\item $E_i + E_j = H$ and $E_i \cap E_j = 0$ for 
$(i,j) = (1,2), (2,3)$ and $(4,1)$  
\end{enumerate} 
\label{cor:bounded operator system}
\end{cor}
\begin{proof}
(1)$\Rightarrow$(2): It is trivial. 
(2)$\Rightarrow$(1):
Since $E_3 + E_2 = H$,  we have 
$D(T)= E_1 \cap (E_3 + E_2) = E_1$. 
Because $\graph T = E_3$ is closed , 
$T$ is bounded by the closed graph theorem. 
Similarly 
$E_4 + E_1 = H$ implies that $D(S) = E_2$ and $S$ is bounded. 
\end{proof}

\begin{cor}
Let $\mathcal{S}=(H;E_1,E_2,E_3,E_4)$ be a system of  
four subspaces in a Hilbert space $H$. Then the following are 
equivalent:
\begin{enumerate}
\item $\mathcal{S}$ is isomorphic to a bounded operator system associated with a single operator. 
\item $E_i + E_j = H$ and $E_i \cap E_j = 0$ for 
$(i,j) = (1,2), (2,3),(4,1)$ and $(2,4)$.   
\end{enumerate} 
\end{cor}
\begin{proof}
(1)$\Rightarrow$(2): It is trivial. 
(2)$\Rightarrow$(1): By the preceding Corollary, 
$\mathcal{S}$ is isomorphic to a bounded operator system 
$\mathcal{S}_{T,S}$. Since $E_2 \cap E_4 = 0$, $S$ is one to one. 
Since $E_2 + E_4 = H$, $S$ is onto.  Therefore 
$\mathcal{S}_{T,S}$ is isomorphic to a bounded operator system 
$\mathcal{S}_{ST,I} = \mathcal{S}_{ST}$ associated with a single operator 
$ST$ by Proposition \ref{prop:ST}. 
\end{proof}

\begin{prop} Let $\mathcal{S}_T = (H;E_1,E_2,E_3,E_4)$ be  
a bounded operator system associated with a single operator 
$T \in B(K)$.
Then  $\mathcal{S}_T$ is transitive if and only if 
$\dim K = 1$. If it is so, then  $\mathcal{S}_T$ is isomorphic 
to 
\[
(\mathbb{C}^2;\mathbb{C}\oplus 0, 0 \oplus \mathbb{C},
\{(x,\lambda x);x \in \mathbb{C}\},\{(x,x);x \in \mathbb{C}\})
\]
for some $\lambda \in \mathbb{C}$. 
\label{prop:trasitive}
\end{prop}
\begin{proof}Recall that  $\mathcal{S}_T$ is transitive if 
\[ 
End (\mathcal{S}_T) = \{ B \oplus B \in B(H) ; 
B \in B(K), \ BT= TB \} = \mathbb{C}I. 
\]  
Hence $\mathcal{S}_T$ is transitive if and only if 
$\{T\}' := \{B \in B(K) ; \ BT = TB \} = \mathbb{C}I$ 
if and only if $\dim K = 1$. 

\end{proof}

But certain unbounded operators on an infnite dimensional 
Hilbert space give transitive systems of four subspaces.
 
{\bf Example}(Harrison-Radjavi-Rosenthal \cite{HRR})        
Let $K=\ell^{2}(\mathbb{Z})$ and $H=K\oplus K$.
Let  $(a_n)_{n\in \mathbb{Z}}$ be a sequence given by
$a_{n}=1$ for $n\leq 0$
and 
$a_{n}=exp((-1)^{n}n!)$ for $n\geq 1$.
Define a bilateral weighted shift
$T:K \supset D(T) \to K$ by $(Tx)_{n} = a_{n-1}x_{n-1}$
with the domain 
$D(T) = \{(x_{n})_{n}\in \ell^{2}(\mathbb{Z}) ;  
\sum_{n}\vert a_{n}x_{n}\vert^{2}<\infty\}.$
Let
$E_{1}=K\oplus 0$, 
$E_{2}= 0\oplus  K$, 
$E_{3}=\{(x,Tx) \in K \oplus K ;  x \in D(T)\}$,
and $E_{4}=\{(x,x) \in K \oplus K ; x \in K\}$.
Harrison, Radjavi and Rosenthal showed that 
$\{H, E_1, E_2, E_3, E_4, 0 \}$ is a transitive lattice 
in \cite{HRR}. Hence $\mathcal {S} = (H;E_1,E_2,E_3,E_4)$
is a transitive system of four subspaces.

We can extend their example to construct uncountably many 
transitive systems.

\begin{lemma}
Let $\mathcal{S}_T = (H;E_1,E_2,E_3,E_4)$ and 
$\mathcal{S}_{T'}  = (H';E_1',E_2',E_3',$ $E_4')$
be closed operator systems associated with operators
$T: D(T) \rightarrow K$, $T': D(T') \rightarrow K'$.  Then 
$\mathcal{S}_T$ and $\mathcal{S}_{T'}$ are isomorphic 
if and only if $T$ and $T'$ are similar.    
\end{lemma}
\begin{proof} The proof is as same as bounded operators if 
we see the domains of the closed operators carefully.   
\end{proof}

\noindent
{\bf  Example}.  
Let $K=\ell^{2}(\mathbb{Z})$ and
$H=K\oplus K$. For a fixed number $\alpha >1$, let 
$(w_n)_{n\in \mathbb{Z}} = (w_n(\alpha))_{n\in \mathbb{Z}}$ be a sequence given by
$w_{n}=1$ for $n\leq 0$
and 
$w_{n}=exp((-\alpha)^{n})$for $(n\geq 1)$.
Define a bilateral weighted shift
$T_{\alpha}:K \supset D_{\alpha} \to K$ by 
$(T_{\alpha}x)_{n} = w_{n-1}x_{n-1}$
with the domain 
$D_{\alpha} = \{(x_{n})_{n}\in \ell^{2}(\mathbb{Z}) ;  
\sum_{n}\vert w_{n}x_{n}\vert^{2}<\infty\}.$
Let
$E_{1}=K\oplus 0$, 
$E_{2}= 0\oplus  K$, 
$E_{3}^{\alpha}=\{(x,T_{\alpha}x) \in K \oplus K ;  x \in D_{\alpha}$,
and $E_{4}=\{(x,x) \in K \oplus K ; x \in K\}$.

\begin{prop}If $\alpha >1$, then the above system 
$\mathcal {S}_{\alpha}= (H;E_1,E_2,E_3^{\alpha},$ $E_4)$
is a transitive system. Furthermore   
if $\alpha\not= \beta,$ then $\mathcal {S}_{\alpha}$
and $\mathcal {S}_{\beta}$ are not isomorphic. 
\end{prop} 
\begin{proof}
Let $V \in Hom(\mathcal {S}_{\alpha}, \mathcal {S}_{\beta})$. 
Since $V(E_i) \subset E_i$ for $i = 1,2,4$, $V = A \oplus A$ 
for some $A=(a_{ij})_{ij} \in B(K)$. Since 
$V(E_{3}^{\alpha}) \subset E_{3}^{\beta}$ and  $e_n \in D_{\alpha}$, 
\[
(A \oplus A)(e_n,T_{\alpha}e_n) = (Ae_n,AT_{\alpha}e_n) \in E_{3}^{\beta}.
\]
Hence $AT_{\alpha}e_n = T_{\beta}Ae_n$.  Comparing $(m+1)$-th component, 
we have 
$w_n(\alpha) a_{m+1,n+1} = w_m(\beta)a_{m,n}$, that is, 
\[
a_{m+1,n+1} = \frac{w_m(\beta)}{w_n(\alpha)}a_{m,n}. 
\]
Therefore for any $k \in {\mathbb N}$, 
\[
a_{m+k,n+k} = \frac{w_m(\beta)\dots w_{m+k-1}(\beta)}
              {w_n(\alpha)\dots w_{n+k-1}(\alpha)}a_{m,n} 
            = exp(c_k(m,n))a_{m,n}, 
\]
where 

\begin{align*}
c_k(m,n) & = ((-\beta)^m + \dots + (-\beta)^{m +k-1}) 
                - ((-\alpha)^n + \dots + (-\alpha)^{n +k-1}) \\
         & = \frac{(-\beta)^m(1-(-\beta)^k)}{1 + \beta}
            - \frac{(-\alpha)^n(1-(-\alpha)^k)}{1 + \alpha} \ .
\end{align*}

\noindent
(i)(the case when $\alpha = \beta$): Putting $n = m$, we have $c_k(m,m) =0$. 
Hence the diagonal of $A$ is constant. If $A$ were not a multiple of the 
identitiy, then there exist distinct $m$ and $n$ with $a_{m,n} \not= 0$. 
According to $m<n$ or $m>n$,  for a sufficient large $k$, 
\[
c_k(m,n) = ((-\alpha)^m + \dots + (-\alpha)^{n}) 
- ((-\alpha)^{m+k-1} + \dots + (-\alpha)^{n +k-1}) 
\] 
or 
\[
c_k(m,n) =  -((-\alpha)^n + \dots + (-\alpha)^{m}) 
+((-\alpha)^{n+k-1} + \dots + (-\alpha)^{m +k-1}) .
\] 
In either case we have $\limsup _k c_k(m,n) = \infty$.  Hence 
$a_{m+k,n+k}$ is not bounded as $k \rightarrow \infty$. This contradicts 
to that $A$ is bounded.  Therefore $A$ is a scalar.  We have shown that 
$\mathcal {S}_{\alpha}$ is a transitive system.

\noindent
(ii)the case when $\alpha \not= \beta$: We may and do assume that 
$1< \alpha < \beta$. If $A$ were not equal to $0$, then there exist 
$m$ and $n$ with $a_{m,n} \not= 0$. 
Since 
\[
c_k(m,n) =  \frac{(-\beta)^m(1-(-\beta)^k)}{1 + \beta}
          \{1 - \frac{(-\alpha)^n(1+\beta)(1-(-\alpha)^k)}
               {(-\beta)^m(1 + \alpha)(1-(-\beta)^k)}\} \ , 
\]
we have $\limsup _k c_k(m,n) = \infty$. This contradicts 
to that $A$ is bounded.  Therefore $A=0$. We have shown that 
$Hom (\mathcal {S}_{\alpha}, \mathcal {S}_{\beta}) = 0$.  Therefore 
$\mathcal {S}_{\alpha}$
and $\mathcal {S}_{\beta}$ are not isomorphic. 
\end{proof}

\begin{prop}
Let $\mathcal{S}_{T,S} = (H;E_1,E_2,E_3,E_4)$
be a bounded operator system associated with operators
$S \in B(K_2,K_1), T \in B(K_1,K_2)$.
Then $\mathcal{S}$ is transitive if and only if 
$\mathcal{S}$ is isomorphic 
to $(\mathbb{C};\mathbb{C},0,\mathbb{C},0)$, 
$(\mathbb{C};0,\mathbb{C},0,\mathbb{C})$, 
$(\mathbb{C}^2;\mathbb{C}\oplus 0, 0 \oplus \mathbb{C},
\{(x,x);x \in \mathbb{C}\},0 \oplus \mathbb{C})$ 
or $(\mathbb{C}^2;\mathbb{C}\oplus 0, 0 \oplus \mathbb{C},
\{(x,\lambda x);x \in \mathbb{C}\},\{(x,x);x \in \mathbb{C}\})$
for some $\lambda \in \mathbb{C}$.
\end{prop}
\begin{proof}
Suppose that $\mathcal{S} = \mathcal{S}_{T,S}$ is transitive.  
If $\dim H = 1$, then 
$\mathcal{S}$ is isomorphic 
to $(\mathbb{C};\mathbb{C},0,\mathbb{C},0)$ or 
$(\mathbb{C};0,\mathbb{C},0,\mathbb{C})$. 
We assume that $\dim H \geq 2$. 
Since $ST \oplus TS \in End (\mathcal{S}_{T,S})$ and 
$\mathcal{S}$ is transitive, there exists 
$\lambda \in \mathbb{C}$ such that $ST = \lambda I_{K_1}$ and 
$TS = \lambda I_{K_2}$. 

In the case that $\lambda \not= 0$, $T$ and $S$ are invertible 
and $S = \lambda T^{-1}$.  By Proposition \ref{prop:ST}, 
$\mathcal{S}_{T,S}$ is isomorphic to 
$\mathcal{S}_{\lambda I_{K_1},I_{K_1}}$. 
Applying  Proposition \ref{prop:trasitive}, 
$\mathcal{S}$ is isomorphic to 
\[
(\mathbb{C}^2;\mathbb{C}\oplus 0, 0 \oplus \mathbb{C},
\{(x,\lambda x);x \in \mathbb{C}\},\{(x,x);x \in \mathbb{C}\})
\]
for some $\lambda \in \mathbb{C}$. 

In the case that $\lambda = 0$, we have $ST = 0$ and $TS = 0$. 
Since $SS^* \oplus S^*S, T^*T \oplus TT^* 
\in End (\mathcal{S}_{T,S})$ and $\mathcal{S}$ is transitive, 
we have $SS^* = \alpha I_{K_1}, S^*S = \alpha I_{K_2}, 
T^*T = \beta I_{K_1}$ and $TT^* = \beta I_{K_2}$.  Because 
$ST = 0$, $\alpha \beta = 0$.  Hence $\alpha = 0$ or 
$\beta = 0$, so that $S = 0$ or $T = 0$.  If $T = 0$, then 
a subsystem $(H;K_1\oplus 0, 0\oplus K_2, \{(Sy,y); y\in K_2\})$
of three subspaces is transitive. Since $\dim H \geq 2$, 
the subsystem is isomorphic to 
$(\mathbb{C}^2;\mathbb{C}\oplus 0, 0 \oplus \mathbb{C},
\{(x,x);x \in \mathbb{C}\})$. Hence $\mathcal{S}$ is isomorphic to 
$(\mathbb{C}^2;\mathbb{C}\oplus 0, 0 \oplus \mathbb{C},
\mathbb{C}\oplus 0,\{(x,x);x \in \mathbb{C}\})$ . Similarly if 
$S= 0$, then $\mathcal{S}$ is isomorphic to 
$(\mathbb{C}^2;\mathbb{C}\oplus 0, 0 \oplus \mathbb{C},
\{(x,x);x \in \mathbb{C}\},0 \oplus \mathbb{C})$. 
The converse is clear. 
\end{proof}

\section{classification theorem by Gelfand-Ponomarev}
One of the main problem to attack is a classification of 
indecomposable systems $\mathcal S = (H;E_1,E_2,E_3,$ $E_4)$ of  
four subspaces in a Hilbert space $H$. 
In the case when $H$ is finite-dimensional, 
Gelfand and Ponomarev  completely
classified indecomposable systems and gave a complete list of 
them in \cite{GP}. 
The important numerical invariants are $dim \ H$ and the 
defect defined by 
$$
\rho (\mathcal S) := \sum_{i=1}^4 dim \ E_i - 2dim \ H.
$$
\begin{thm}[Gelfand-Ponomarev \cite{GP}]
The set of possible values of the defect 
$\rho (\mathcal S)$ for indecomposable systems $\mathcal S$
 of four subspaces in a finite-dimensional space 
is exactly the set $\{-2,-1,0,1,2\}$.
\end{thm}

The defect characterizes  an essential feature of the system.
If $\rho (\mathcal S) = 0$, then ${\mathcal S}$ is 
isomorphic to a bounded operator system up to permutation of subspaces ,
that is, 
there exists a permutation $\sigma$ on $\{1,2,3,4\}$ and    
a pair of linear operators
$A: E \rightarrow F$ and $B: F \rightarrow E$ 
such that 
$H = E \oplus F$, $E_{\sigma(1)} = E \oplus 0$, 
$E_{\sigma(2)} = 0 \oplus F$, 
$E_{\sigma(3)} = \{(x,Ax) \in H ; x \in E\}$ and 
$E_4 =   \{(By,y) \in H ; y \in F \}$.
If $\rho (\mathcal S) = \pm 1$, 
$\mathcal S$ is represented up to permutation by $H = E \oplus F$,
$E_1 = E \oplus 0$, $E_2 = 0 \oplus F$, $E_3$ and $E_4$ are subspaces
of $H$ that do not reduced to the graphs of the operators as in the case that
$\rho (\mathcal S) = 0$.  
A system with $\rho (\mathcal S) = \pm 2$ cannot be
described in the above forms.
\par
Following \cite{GP}, we recall the canonical forms of indecomposable 
systems  
$\mathcal S = (H;E_1,E_2,E_3,E_4)$ of four subspaces in a finite-dimensional
space $H$ up to permutation  in the following: 
\noindent
(A) the case when $dim \ H = 2k$ for some positive integer $k$.

There exist no indecomposable systems $\mathcal S$ with $\rho (\mathcal S) =
\pm 2$.  Let $H$ be a space with a basis 
$\{e_1,\ldots, e_k,f_1,\ldots,f_k\}$.

\noindent
(1)$\mathcal S _3(2k,-1) = (H;E_1,E_2,E_3,E_4)$ with $\rho (\mathcal S) = -1$
\begin{eqnarray*}
&&H = [e_1,\ldots, e_k,f_1,\ldots,f_k], \\
&&E_1 = [e_1,\ldots, e_k], \  E_2 = [f_1,\ldots,f_k],\\
&&E_3 = [(e_2+f_1), \ldots , (e_k+f_{k-1})],\\
&&E_4 = [(e_1 + f_1), \ldots, (e_k + f_k)].
\end{eqnarray*}

\noindent
(2)$S _3(2k,1) = (H;E_1,E_2,E_3,E_4)$ with $\rho (\mathcal S) = 1$ 
\begin{eqnarray*}
&&H = [e_1,\ldots, e_k,f_1,\ldots,f_k], \\
&&E_1 = [e_1,\ldots, e_k], \  E_2 = [f_1,\ldots,f_k],\\
&&E_3 = [e_1, (e_2+f_1), \ldots , (e_k+f_{k-1}), f_k],\\
&&E_4 = [(e_1 + f_1), \ldots, (e_k + f_k)].
\end{eqnarray*}

\noindent
(3)$\mathcal S _{1,3}(2k,0)  = (H;E_1,E_2,E_3,E_4)$ with 
$\rho (\mathcal S) = 0$ 
\begin{eqnarray*}
&&H = [e_1,\ldots, e_k,f_1,\ldots,f_k], \\
&&E_1 = [e_1,\ldots, e_k], \  E_2 = [f_1,\ldots,f_k],\\
&&E_3 = [e_1, (e_2+f_1), \ldots , (e_k+f_{k-1})],\\
&&E_4 = [(e_1 + f_1), \ldots, (e_k + f_k)].
\end{eqnarray*}

\noindent
(4)$\mathcal S (2k,0;\lambda) = (H;E_1,E_2,E_3,E_4)$ with 
$\rho (\mathcal S) = 0$ 
\begin{eqnarray*}
&&H = [e_1,\ldots, e_k,f_1,\ldots,f_k], \\
&&E_1 = [e_1,\ldots, e_k], \  E_2 = [f_1,\ldots,f_k],\\
&&E_3 = [(e_1+\lambda f_1), (e_2+f_1 +\lambda f_2),
\ldots , (e_k+f_{k-1}+ \lambda f_k)],\\
&&E_4 = [(e_1 + f_1), \ldots, (e_k + f_k)].
\end{eqnarray*}

Every other system  $\mathcal S _i(2k,\rho)$, $\mathcal S _{i,j}(2k,0)$ can be
obtained from the systems $\mathcal S _3(2k,\rho)$, $\mathcal S _{i,3}(2k,0)$
by a suitable permutation of the subspaces.  
Let $\sigma_{i,j}$  be the transposition $(i,j)$.
We put $\mathcal S _i(2k,\rho) = \sigma_{3,i}\mathcal S _3(2k,\rho)$ for $\rho =-1, 1$.
We also define $\mathcal S _{i,j}(2k,0) = \sigma_{1,i}\sigma_{3,j}\mathcal S _{1,3}(2k,0)$
for $i,j \in \{1,2,3,4 \}$.

\noindent
(B)the case $dim \ H=2k+1$ is odd for some integer $k \geq 0$ .
Let $H$ be a space with a basis 
$\{e_1,\ldots, e_k,e_{k+1},f_1,\ldots,f_k\}$.

\noindent
(5)$\mathcal S _1(2k+1,-1) = (H;E_1,E_2,E_3,E_4)$ 
with $\rho (\mathcal S) = -1$ 
\begin{eqnarray*}
&&H = [e_1,\ldots, e_k,e_{k+1}, f_1,\ldots,f_k], \\
&&E_1 = [e_1,\ldots, e_k, e_{k+1}], \  E_2 = [f_1,\ldots,f_k],\\
&&E_3 = [(e_2+f_1), \ldots ,(e_{k+1}+f_{k})],\\
&&E_4 = [(e_1 + f_1), \ldots, (e_k + f_k)].
\end{eqnarray*}

\noindent
(6)$\mathcal S _2(2k+1,1) = (H;E_1,E_2,E_3,E_4)$ 
with $\rho (\mathcal S) = 1$ 
\begin{eqnarray*}
&&H = [e_1,\ldots, e_k,e_{k+1}, f_1,\ldots,f_k], \\
&&E_1 = [e_1,\ldots, e_k, e_{k+1}], \  E_2 = [f_1,\ldots,f_k],\\
&&E_3 = [e_1, (e_2+f_1), \ldots , (e_{k+1}+f_{k})],\\
&&E_4 = [(e_1 + f_1), \ldots, (e_k + f_k),e_{k+1}].
\end{eqnarray*}

\noindent
(7)$\mathcal S _{1,3}(2k+1,0) = (H;E_1,E_2,E_3,E_4)$ 
with $\rho (\mathcal S) = 0$ 
\begin{eqnarray*}
&&H = [e_1,\ldots, e_k,e_{k+1}, f_1,\ldots,f_k], \\
&&E_1 = [e_1,\ldots, e_k, e_{k+1}], \  E_2 = [f_1,\ldots,f_k],\\
&&E_3 = [e_1, (e_2+f_1), \ldots , (e_{k+1}+f_{k})],\\
&&E_4 = [(e_1 + f_1), \ldots, (e_k + f_k)].
\end{eqnarray*}

\noindent
(8)$\mathcal S (2k+1,-2) = (H;E_1,E_2,E_3,E_4)$ 
with $\rho (\mathcal S) = -2$ 
\begin{eqnarray*}
&&H = [e_1,\ldots, e_k,e_{k+1}, f_1,\ldots,f_k], \\
&&E_1 = [e_1,\ldots, e_k], \  E_2 = [f_1,\ldots,f_k],\\
&&E_3 = [(e_2+f_1), \ldots , (e_{k+1}+f_{k})],\\
&&E_4 = [(e_1 + f_2), \ldots, (e_{k-1} + f_k),(e_k+e_{k+1})].
\end{eqnarray*}

\noindent
(9)$\mathcal S (2k+1,2) = (H;E_1,E_2,E_3,E_4)$ 
with $\rho (\mathcal S) = 2$ 
\begin{eqnarray*}
&&H = [e_1,\ldots, e_k,e_{k+1}, f_1,\ldots,f_k], \\
&&E_1 = [e_1,\ldots, e_k,e_{k+1}], \ \ E_2 = [f_1,\ldots,f_k, e_{k+1}],\\
&&E_3 = [e_1, (e_2+f_1), \ldots ,(e_{k+1}+f_{k})], \\
&&E_4 = [f_1, (e_1 + f_2), \ldots, (e_{k-1} + f_k),(e_k+e_{k+1})].
\end{eqnarray*}

We put $\mathcal S _i(2k+1,-1)=\sigma_{1,i}\mathcal S _1(2k+1,-1)$,
\ $\mathcal S _i(2k+1,+1) = \sigma_{2,i}\mathcal S _2(2k+1,1)$,
\ $\mathcal S _{i,j}(2k+1,0) 
= \sigma_{1,i}\sigma_{3,j}\mathcal S _{1,3}(2k+1,0)$
for $i,j \in \{1,2,3,4 \}$.

\begin{thm}[Gelfand-Ponomarev \cite{GP}] 
If a system ${\mathcal S}$ of four subspaces in a finite-dimensional $H$
is indecomposable, then $\mathcal S$ is isomorphic to one of the 
following systems: 

$\mathcal S _{i,j}(m,0)$, $(i<j , i,j \in \{1,2,3,4\}, m = 1,2,...)$;
$\mathcal S (2k,0;\lambda)$, $(\lambda \in \mathbb C, \lambda \not=0, \lambda
\not=1,
k=1,2,...)$, 
 $\mathcal S_i(m,-1)$, $\mathcal S_i(m,1)$, $(i \in \{1,2,3,4\}, m =
1,2,...)$; $\mathcal S (2k+1,-2)$, $\mathcal S(2k+1,+2), (k=0,1,...)$.
\end{thm}

\noindent 
{\bf Remark.}It is known that if $\mathcal S$ is an indecomposable 
system of four subspaces in the above Theorem satisfying  $\rho (\mathcal S) \not= 0$, then $\mathcal S$ is transitive, for example, see \cite{B}.

\section{exotic indecomposable systems of four subspaces}

In this section we shall construct uncountably many, exotic, 
indecomposable systems of four subspaces, that is, 
indecomposable systems 
which are not isomorphic to any closed operator system 
under any permutaion of subspaces.   

\bigskip
\noindent
{\bf Exotic examples.} Let $L =\ell^{2}(\mathbb{N})$ with a 
standard basis $\{e_1,e_2, \dots  \}$.  Put $K = L \oplus L$ and 
$H = K \oplus K = L \oplus L \oplus L \oplus L$.  
Consider a unilateral shift
$S:L \rightarrow  L$ by $Se_n = e_{n+1}$ for $n = 1,2,\dots$ .
For a fixed paramater $\gamma \in {\mathbb C}$ with 
$|\gamma | \geq 1$, we consider an operator 
\[
T_{\gamma} = 
    \begin{pmatrix}
     \gamma S^* & I \\
     0         & S 
    \end{pmatrix} 
\in B(K) = B(L \oplus L). 
\] 
Let
$E_1=K\oplus 0$, 
$E_2= 0\oplus  K$, \\
$E_3=\{(x,T_{\gamma}x) \in K \oplus K ;  x \in K\} + {\mathbb C}(0,0,0,e_1) 
= \graph T_{\gamma} +  {\mathbb C}(0,0,0,e_1) $,
and $E_{4}=\{(x,x) \in K \oplus K ; x \in K\}$.
Consider a system $\mathcal{S}_{\gamma}=(H;E_1,E_2,E_3,E_4)$. 
We shall show that $\mathcal{S}_{\gamma}$ is indecomposable.  
If $|\gamma | > 1$, then $\mathcal{S}_{\gamma}$ is not
isomorphic to any closed operator systems under any permutation.  
We could  regard the system  $\mathcal{S}_{\gamma}$ is a 
one-dimensional \lq\lq deformation" of an operator system. 
First we start with an easy fact. 

\begin{lemma} Assume that  a bounded operator $A\in B(\ell^{2}(\mathbb{N}))$ 
is represented as an upper triangular matrix 
$A = (a_{ij})_{ij}$ by a standard basis $\{e_1,e_2, \dots  \}$. 
If the diagonal is constant $\lambda$, i.e., 
$a_{ii} = \lambda$ for $i=1, \dots$, and $A$ is an idempotent, 
then $A = 0$ or $A = I$. 
\label{lemma:Toeplitz-idempotent}
\end{lemma}
\begin{proof}
Put $N = A - \lambda I$.  Then $N$ is an upper triangular matrix 
with zero diagonal.  Comparing the diagonals for 
\[
\lambda I + N = A = A^2 = \lambda ^2I + 2\lambda N + N^2 , 
\]
we have  $\lambda ^2 = \lambda$. Hence $\lambda = 0$ or $1$. 
If $\lambda = 0$, then $N^2 = N$.  
Since $N$ is an idempotent and an upper triangular matrix 
with zero diagonal,  
$N = 0$, that is, $A = 0$. If $\lambda = 1$, then 
$(I - A)$ is an idempotent and an upper triangular matrix 
with zero diagonal, $I - A = 0$, that is, $A = I$.

\end{proof}

\begin{thm}
If $|\gamma | \geq 1$, then the above system  
$\mathcal{S}_{\gamma}=(H;E_1,E_2,E_3,E_4)$ is 
indecomposable.
\label{thm:exotic examples} 
\end{thm}
\begin{proof}We shall show that 
$\{V \in End (\mathcal{S}_{\gamma}) ; V^2 = V\} = \{0,I\}$. 
Let $V \in End (\mathcal{S}_{\gamma})$ satisfy $V^2 = V$. 
Since $V(E_i) \subset E_i$ for $i = 1,2,4$, we have 
\[
V = 
    \begin{pmatrix}
     U & 0 \\
     0 & U \\ 
    \end{pmatrix} 
  \in B(H) \ \ \text{ for some } \ U \in B(K)
\]
We write 
\[
U = 
    \begin{pmatrix}
     A & B \\
     C & D 
    \end{pmatrix} 
\in B(K), 
\]
for some $
A = (a_{ij})_{ij}, \ B = (b_{ij})_{ij},
          \ C = (c_{ij})_{ij}, \ D = (d_{ij})_{ij} \in B(K)$. 
We shall investigate the condition that $V(E_3) \subset E_3$. 
Since $E_3 = \graph T_{\gamma} +  {\mathbb C}(0,0,0,e_1) $, $E_3$ is spanned by 
\[
\{
\begin{pmatrix}
 e_1 \\
 0   \\
 0   \\
 0
 \end{pmatrix}, \ 
\begin{pmatrix}
 e_m \\
 0   \\
 \gamma e_{m-1}  \\
 0
 \end{pmatrix}, \ 
\begin{pmatrix}
 0 \\
 e_n   \\
 e_n  \\
 e_{n+1}
 \end{pmatrix}, 
\begin{pmatrix}
 0 \\
 0  \\
 0  \\
 e_1
 \end{pmatrix}
; m = 2,3,\dots   , n = 1,2,\dots
\ 
\}.
\]
We may write 
\[
E_3 = \{
\begin{pmatrix}
 (\lambda _n)_n \\
 (\mu _n)_n  \\
 (\gamma \lambda _{n+1} + \mu _n)_n \\
 (\alpha , (\mu _n)_n)
\end{pmatrix}
; \lambda _n, \mu _n, \alpha \in {\mathbb C}, 
 \ \sum _n |\lambda _n|^2 < \infty, \ 
  \sum _n |\mu _n|^2 < \infty
\}
\]
Since $(e_1,0,0,0) \in E_3$, we have 
\[
\begin{pmatrix}
     A & B & 0 & 0 \\
     C & D & 0 & 0 \\
     0 & 0 & A & B \\
     0 & 0 & C & D 
    \end{pmatrix}
\begin{pmatrix}
 e_1 \\
 0   \\
 0   \\
 0
 \end{pmatrix}
= 
 \begin{pmatrix}
 Ae_1 \\
 Ce_1   \\
 0   \\
 0
 \end{pmatrix}
=
\begin{pmatrix}
 (\lambda _m)_m \\
 (\mu _m)_m  \\
 (\gamma \lambda _{m+1} + \mu _m)_m \\
 (\alpha , (\mu _m)_m)
\end{pmatrix}
\in E_3.
\]
Then, for any $m = 1,2,\dots $, we have $c_{m1} = \mu_m = 0$. 
Moreover $0 = \gamma \lambda _{m+1} + \mu _m = \gamma \lambda _{m+1}$.  
Hence $\lambda _{m+1} = 0$ because $\gamma \not= 0$. 
Therefore $a_{m+1,1} = \lambda _{m+1} = 0$. Thus the first column of 
$C$ is zero and  the first column of $A$ is zero except $a_{11}$. 
We shall show that $C = 0$ and $A$ is an upper triangular 
Toeplitz matrix with by the induction of $n$-th columns.　
The case when $n =1$ is already shown. Assume that the assertion 
hold for $n$-th columns.    
Since $(e_{n+1},0,\gamma e_n,0) \in E_3$, we have 
\[
\begin{pmatrix}
     A & B & 0 & 0 \\
     C & D & 0 & 0 \\
     0 & 0 & A & B \\
     0 & 0 & C & D 
    \end{pmatrix}
\begin{pmatrix}
 e_{n+1} \\
 0   \\
 \gamma e_n   \\
 0
 \end{pmatrix}
=
 \begin{pmatrix}
 Ae_{n+1} \\
 Ce_{n+1}   \\
 \gamma Ae_n  \\
 \gamma Ce_n
 \end{pmatrix}
=
\begin{pmatrix}
 (\lambda _m)_m\\
 (\mu _m)_m  \\
 (\gamma \lambda _{m+1} + \mu _m)_m \\
 (\alpha , (\mu _m)_m)
\end{pmatrix}
\in E_3.
\]
Then  $c_{m,n+1} = \mu_m = \gamma c_{m+1,n} = 0$. 
And  $\gamma a_{m,n}= 
\gamma \lambda _{m+1} + \mu _m = \gamma \lambda _{m+1}$.  
Since $\gamma \not= 0$, 
$a_{m,n} = \lambda _{m+1} = a_{m+1,n+1}$. 
Thus  we have shown that 
$C = 0$ and $A$ is an upper triangular Toeplitz matrix. 
Since $V$ is an idempotent, so is 
\[
U = 
    \begin{pmatrix}
     A & B \\
     0 & D 
    \end{pmatrix} 
. 
\]
Hence $A$ is also an idempotent.  
By Lemma \ref{lemma:Toeplitz-idempotent}, 
we have two cases $A = 0$ or $A = I$.  

\noindent
(i)the case  $A = 0$: we shall show that $B = D = 0$. This 
 immediately implies  $U=0$, so that $V = 0$. 

\noindent
(ii)the case $A = I$: Since $I - V \in  End (\mathcal{S}_{\gamma})$ is 
is also an idempotent and it can be reduced to the case (i) and 
we have $V = I$.

Hence we may assume that $A = 0$. Since $U$ is an idempotent, 
$D$ is also an idempotent. Since $(0,0,0,e_1) \in E_3$, we have 
\[
\begin{pmatrix}
     0 & B & 0 & 0 \\
     0 & D & 0 & 0 \\
     0 & 0 & 0 & B \\
     0 & 0 & 0 & D 
    \end{pmatrix}
\begin{pmatrix}
 0   \\
 0   \\
 0   \\
 e_1
 \end{pmatrix}
= 
 \begin{pmatrix}
 0 \\
 0  \\
 Be_1   \\
 De_1
 \end{pmatrix}
=
\begin{pmatrix}
 (\lambda _m)_m \\
 (\mu _m)_m  \\
 (\gamma \lambda _{m+1} + \mu _m)_m \\
 (\alpha , (\mu _m)_m)
\end{pmatrix}
\in E_3.
\]
Then, for any $m = 1,2,\dots $, we have $\mu_m = \lambda _m = 0$. 
Hence $b_{m1} = \gamma \lambda _{m+1} + \mu _m = 0 $ and 
$d_{m+1,1} = \mu _m = 0$.  Thus the first column of 
$B$ is zero and  the first column of $D$ is zero except $d_{11}$. 
We shall show that $D$ is an upper triangular 
Toeplitz matrix by the induction of $n-$th columns.　
The case when $n =1$ is already shown. Assume that the assertion 
hold for $n-$th columns. Since $(0,e_n,e_n,e_{n+1}) \in E_3$, 
\[
\begin{pmatrix}
     0 & B & 0 & 0 \\
     0 & D & 0 & 0 \\
     0 & 0 & 0 & B \\
     0 & 0 & 0 & D 
    \end{pmatrix}
\begin{pmatrix}
 0 \\
 e_n \\
 e_n   \\
 e_{n+1}
 \end{pmatrix}
=
 \begin{pmatrix}
 Be_n \\
 De_n\\
 Be_{n+1}   \\
 De_{n+1} 
 \end{pmatrix}
=
\begin{pmatrix}
 (\lambda _m)_m\\
 (\mu _m)_m  \\
 (\gamma \lambda _{m+1} + \mu _m)_m \\
 (\alpha , (\mu _m)_m)
\end{pmatrix}
\in E_3. 
\]
We have $d_{m+1,n+1} = \mu_m = d_{mn}$. 
Hence $D$ is an upper triangular Toeplitz matrix. 
Since $D$ is also an idempotent,  $D = O$ or $D = I$ by Lemma 
\ref{lemma:Toeplitz-idempotent}. 

If $D = 0$, then $U = U^2 = 0$.  Thus $B = 0$, and the assertion 
is verified. We shall show that the case when $D = I$ will not 
occur.  On the contrary, suppose that $D = I$.      
 We have 
\[
V
\begin{pmatrix}
 0   \\
 0   \\
 0   \\
 e_1
 \end{pmatrix}
=
\begin{pmatrix}
     0 & B & 0 & 0 \\
     0 & I & 0 & 0 \\
     0 & 0 & 0 & B \\
     0 & 0 & 0 & I 
    \end{pmatrix}
\begin{pmatrix}
 0   \\
 0   \\
 0   \\
 e_1
 \end{pmatrix}
= 
 \begin{pmatrix}
 0 \\
 0  \\
 Be_1   \\
 e_1
 \end{pmatrix}
=
\begin{pmatrix}
 (\lambda _m)_m \\
 (\mu _m)_m  \\
 (\gamma \lambda _{m+1} + \mu _m)_m \\
 (\alpha , (\mu _m)_m)
\end{pmatrix}
\in E_3.
\]
Then, for any $m = 1,2,\dots $, we have $\mu_m = \lambda _m = 0$. 
Hence $b_{m1} = \gamma \lambda _{m+1} + \mu _m = 0 $ 
Thus the first column of $B$ is zero. 
We shall show that $B$ should be the following form 
by the induction of $n-$th columns:
\[
B =
\begin{pmatrix}
0&1&0&\gamma&0&\gamma^2&0&\gamma^3&0&\cdots \\
0&0&1&0&\gamma&0&\gamma^2&0&\gamma^3&\cdots \\
0&0&0&1&0&\gamma&0&\gamma^2&0&\ddots\\
0&0&0&0&1&0&\gamma&0&\gamma^2&\ddots\\
0&0&0&0&0&1&0&\gamma&0&\ddots\\
0&0&0&0&0&0&1&0&\gamma&\ddots\\
0&0&0&0&0&0&0&1&0&\ddots\\
0&0&0&0&0&0&0&0&1&\ddots\\
0&0&0&0&0&0&0&0&0&\ddots\\
\vdots&\vdots&\vdots&\vdots&\vdots&\vdots&\vdots&\ddots&\ddots&\ddots
\end{pmatrix}, 
\]
that is, $b_{ij} = \gamma^{k-1}$ if $j>i$ and $j-i = 2k -1$, and 
$b_{ij} = 0$ if otherwise. 

The case when $n =1$ is already shown. Assume that the assertion 
hold for $n$-th columns.  Since 
\[
\begin{pmatrix}
     0 & B & 0 & 0 \\
     0 & I & 0 & 0 \\
     0 & 0 & 0 & B \\
     0 & 0 & 0 & I 
    \end{pmatrix}
\begin{pmatrix}
 0 \\
 e_n \\
 e_n   \\
 e_{n+1}
 \end{pmatrix}
=
 \begin{pmatrix}
 Be_n \\
 e_n\\
 Be_{n+1}   \\
 e_{n+1} 
 \end{pmatrix}
=
\begin{pmatrix}
 (\lambda _m)_m\\
 (\mu _m)_m  \\
 (\gamma \lambda _{m+1} + \mu _m)_m \\
 (\alpha , (\mu _m)_m)
\end{pmatrix}
\in E_3, 
\]
for any $m = 1,2,\dots $, we have $\mu_m = \delta_{m,n}$. 
And  
\[
b_{m,n+1} = \gamma \lambda _{m+1} + \mu _m 
= \gamma \lambda _{m+1} + \delta_{m,n},
\]
that is, 
\[
((n+1)\text{-th column of }B) 
= \gamma S^*(n \text{ -th column of }B) + e_n. 
\]
By the induction  we have shown that 
$B$ is the above  form.  But then 
\[
\|B^*e_1\|^2 = \|(\text {the first row of }B)\|^2 
= \sum _{k=1}^{\infty} |\gamma|^{2(k-1)} = \infty, 
\] 
because $|\gamma| \geq 1$.  This contradicts to that 
$B$ is bounded. Therefore $D \not= I$.  This finishes the proof.
\end{proof}

\begin{thm}
If $|\beta | \geq 1$, $|\gamma | \geq 1$ and 
$|\beta | \not= |\gamma |$, 
then the above systems 
$\mathcal{S}_{\beta}=(H;E_1,E_2,E_3^{\beta},E_4)$ and  
$\mathcal{S}_{\gamma}=(H;E_1,E_2,E_3^{\gamma},E_4)$ are 
not isomorphic.
\end{thm}
\begin{proof}
 On the contrary, suppose that there were 
an isomorphism 
$\ V: \mathcal{S}_{\beta} \rightarrow  \mathcal{S}_{\gamma}$. 
We shall show a contradiction. We may and do assume that 
$|\beta | > |\gamma |$. 
Since $V(E_i) = E_i$ for $i = 1,2,4$, we have 
\[
V = 
    \begin{pmatrix}
     U & 0 \\
     0 & U \\ 
    \end{pmatrix} 
  \in B(H) \ \ \text{ for some invertible } \ U \in B(K)
\]
We write 
\[
U = 
    \begin{pmatrix}
     A & B \\
     C & D 
    \end{pmatrix} 
\in B(K), 
\]
for some $
A = (a_{ij})_{ij}, \ B = (b_{ij})_{ij},
          \ C = (c_{ij})_{ij}, \ D = (d_{ij})_{ij} \in B(K)$. 
We shall investigate the condition that 
$V(E_3^{\beta}) = E_3^{\gamma}$. 
Since $E_3^{\beta} = \graph T_{\beta} +  {\mathbb C}(0,0,0,e_1) $, 
$E_3^{\beta}$ is spanned by 
\[
\{
\begin{pmatrix}
 e_1 \\
 0   \\
 0   \\
 0
 \end{pmatrix}, \ 
\begin{pmatrix}
 e_m \\
 0   \\
 \beta e_{m-1}  \\
 0
 \end{pmatrix}, \ 
\begin{pmatrix}
 0 \\
 e_n   \\
 e_n  \\
 e_{n+1}
 \end{pmatrix}, 
\begin{pmatrix}
 0 \\
 0  \\
 0  \\
 e_1
 \end{pmatrix}
; m = 2,3,\dots   , n = 1,2,\dots
\ 
\}.
\]
We also  write 
\[
E_3^{\gamma} = \{
\begin{pmatrix}
 (\lambda _n)_n \\
 (\mu _n)_n  \\
 (\gamma \lambda _{n+1} + \mu _n)_n \\
 (\alpha , (\mu _n)_n)
\end{pmatrix}
; \lambda _n, \mu _n, \alpha \in {\mathbb C}, 
 \ \sum _n |\lambda _n|^2 < \infty, \ 
  \sum _n |\mu _n|^2 < \infty
\}. 
\]
Since $(e_1,0,0,0) \in E_3^{\beta}$, we have 
\[
0 \not= 
\begin{pmatrix}
     A & B & 0 & 0 \\
     C & D & 0 & 0 \\
     0 & 0 & A & B \\
     0 & 0 & C & D 
    \end{pmatrix}
\begin{pmatrix}
 e_1 \\
 0   \\
 0   \\
 0
 \end{pmatrix}
= 
 \begin{pmatrix}
 Ae_1 \\
 Ce_1   \\
 0   \\
 0
 \end{pmatrix}
=
\begin{pmatrix}
 (\lambda _m)_m \\
 (\mu _m)_m  \\
 (\gamma \lambda _{m+1} + \mu _m)_m \\
 (\alpha , (\mu _m)_m)
\end{pmatrix}
\in E_3^{\gamma}.
\]
Then, for any $m = 1,2,\dots $, we have $c_{m1} = \mu_m = 0$. 
Moreover $0 = \gamma \lambda _{m+1} + \mu _m = \gamma \lambda _{m+1}$.  
Hence $\lambda _{m+1} = 0$ because $\gamma \not= 0$. 
Therefore $a_{m+1,1} = \lambda _{m+1} = 0$. Thus the first column of 
$C$ is zero and  the first column of $A$ is zero except $a_{11}$. 
Since $Ae_1 \not= 0$, $a_{11} \not= 0$.
We shall show that $C = 0$ and $A$ is an upper triangular 
matrix satisfying 
\[
a_{i+1,j+1} = \frac{\beta}{\gamma} a_{ij}  \ \ \text{ if }
 i \leq j 
\]
and $a_{ij} = 0$ if $i>j$, by the induction of $n$-th columns.　
The case when $n =1$ is already shown. Assume that the assertion 
hold for $n$-th columns.    
Since $(e_{n+1},0,\beta e_n,0) \in E_3^{\beta}$, we have 
\[
\begin{pmatrix}
     A & B & 0 & 0 \\
     C & D & 0 & 0 \\
     0 & 0 & A & B \\
     0 & 0 & C & D 
    \end{pmatrix}
\begin{pmatrix}
 e_{n+1} \\
 0   \\
 \beta e_n   \\
 0
 \end{pmatrix}
=
 \begin{pmatrix}
 Ae_{n+1} \\
 Ce_{n+1}   \\
 \beta Ae_n  \\
 \beta Ce_n
 \end{pmatrix}
=
\begin{pmatrix}
 (\lambda _m)_m\\
 (\mu _m)_m  \\
 (\gamma \lambda _{m+1} + \mu _m)_m \\
 (\alpha , (\mu _m)_m)
\end{pmatrix}
\in E_3^{\gamma}.
\]
Then we have $c_{m,n+1} = \mu_m = \beta c_{m+1,n} = 0$. 
Moreover 
\[
\beta a_{m,n}= 
\gamma \lambda _{m+1} + \mu _m = \gamma \lambda _{m+1} 
= \gamma a_{m+1,n+1}.   
\]
Since $\gamma \not= 0$, 
$a_{m+1,n+1} = \frac{\beta}{\gamma} a_{m,n}.$
This completes the induction.  
Then we have 
\[
|a_{nn}| = |\frac{\beta}{\gamma}|^{n-1}|a_{11}| \rightarrow 
\infty ,
\]
because $a_{11} \not= 0$ and $|\frac{\beta}{\gamma}| > 1$. 
But This contradicts to that the operator $A$ is bounded. 
Therefore 
$\mathcal{S}_{\beta}$ and  $\mathcal{S}_{\gamma}$ are 
not isomorphic.
\end{proof}

Next we shall show that if $\gamma >1 $, then 
$\mathcal{S}_{\gamma}$ is not isomorphic to any 
closed operator system.  We introduce a necessary criterion 
for the purpose. 

\bigskip
\noindent
{\bf Definition}(intersection diagram) Let 
$\mathcal{S} =(H;E_1,E_2,E_3,E_4)$ be a system of fours 
subspaces. The {\it intersection diagram} for a system $\mathcal{S}$
is an undirected graph 
$\Gamma_{\mathcal{S}} 
= (\Gamma_{\mathcal{S}}^0,\Gamma_{\mathcal{S}}^1)$ 
with the set of vertices $\Gamma_{\mathcal{S}}^0$ and the 
set of edges $\Gamma_{\mathcal{S}}^1$ defined by 
$\Gamma_{\mathcal{S}}^0 = \{1,2,3,4\}$ 
and for $i \not= j \in \{1,2,3,4\}$ 
\[
\circ _i  \ ^{\line(1,0){20}}  \circ _j \ \text{ if and only if } 
E_i \cap E_j = 0 .
\]

\begin{lemma}
Let 
$\mathcal{S} = \mathcal{S}_{T,S} =(H;E_1,E_2,E_3,E_4)$ be a closed 
operator system. Then the intersection diagram $\Gamma_{\mathcal{S}}$ 
for the system $\mathcal{S}$ contains 
\[
\circ _4 \ ^{\line(1,0){20}} \circ _1 \ ^{\line(1,0){20}} 
\circ _2 \ ^{\line(1,0){20}} \circ _3 \ , 
\]
that is, $E_4 \cap E_1 = 0$, $E_1 \cap E_2 = 0$ and 
$E_2 \cap E_3 = 0$.  In particular, then the intersection diagram 
$\Gamma_{\mathcal{S}}$ is a connected graph.  
\end{lemma}
\begin{proof}
It follows form Proposition \ref{prop:operator-system}. 
\end{proof}  

\begin{prop}
If $\gamma > 1$, then the system $\mathcal{S}_{\gamma}$ 
is not isomorphic to any closed operator system under any 
permutation of subspaces. 
\end{prop}
\begin{proof}
It is clear that $E_4 \cap E_1 = 0$, $E_1 \cap E_2 = 0$ and 
$E_2 \cap E_4 = 0$.  Since $(e_1,0,0,0) \in E_1 \cap E_3$, 
we have $E_1 \cap E_3 \not= 0$.  Because 
$(0,0,0,e_4) \in E_2 \cap E_3$, 
we have $E_2 \cap E_3 \not= 0$. 
Since $|\gamma| > 1$, 
$a := (1,\gamma^{-1}, \gamma^{-2}, \gamma^{-3},..., ) 
\in \ell ^2(\mathbb N)$.  Then $(a,0,a,0) \in E_3 \cap E_4$, 
so that $E_3 \cap E_4 \not= 0$.  Therefore the vertex $3$ is 
not connected to any other vertices $1, 2, 4$.  Thus 
the intersection diagram 
$\Gamma_{\mathcal{S}}$ is not a connected graph. This implies 
that $\mathcal{S}_{\gamma}$ 
is not isomorphic to any closed operator system under any 
permutation of subspaces. 
\end{proof}

Combining the preceeding two propositions , 
we have the existence of uncountably many, exotic, 
indecomposable systems of four subspaces. 
 
\begin{thm}
There exists uncountably many, indecomposable systems 
of four subspaces which are not isomorphic to any closed 
operator system under any permutation of subspaces.
\end{thm}
\begin{proof}
A family  
$\{\mathcal{S}_{\gamma} ; \gamma > 1, \gamma \in {\mathbb R}\}$
of indecomposable systems above is a desired one.
\end{proof}

\section{Defects for systems of four subspaces.}

Gelfand and Ponomarev introduced an integer valued  invariant
$\rho ({\mathcal S})$, called {\it defect},  for a system 
${\mathcal S} = (H ; E_1, E_2, E_3, E_4)$ of four subspaces by 
\[
\rho ({\mathcal S}) = \sum _{i=1} ^4 \dim E_i - 2\dim H.  
\]
They showed that if a system of four subspaces is 
indecomposable, then the possible value of the defect 
$\rho ({\mathcal S})$ is one of 
five values $\{-2, -1, 0, 1, 2\}$ 
We shall extend their  notion of  defect for a certain class of 
systems relating with Fredholm index. 

Let ${\mathcal S} = (H ; E_1, E_2, E_3, E_4)$ be a system of 
four subspaces. We first introduce elementary numerical invariants 
\[
m_{ij} = \dim (E_i \cap E_j) \text{ and } 
m_{ijk} = \dim (E_i \cap E_j \cap E_k).
\]
Similarly put
\[
n_{ij} = \dim ((E_i + E_j)^{\perp}) \text { and } 
n_{ijk} = \dim ((E_i +  E_j + E_k)^{\perp}). 
\]
If ${\mathcal S}$ is indecomposable and $\dim H \geq 2$, 
then $m_{ijk} = 0$ and $n_{ijk} = 0$ by Proposition \ref{prop:n-1}. 

If $H$ is finite dimensional, then
\begin{align*}
& \dim E_i + \dim E_j -\dim H \\
& = \dim (E_i + E_j) + \dim (E_i \cap E_j) -
(\dim (E_i + E_j) + \dim ((E_i + E_j)^{\perp})) \\
& = \dim (E_i \cap E_j) - \dim ((E_i + E_j)^{\perp})
\end{align*}  

In order to make the  numerical invariant unchanged under 
any permutation of subspaces, counting $\ _4C_2 = 6$ pairs of subspaces 
\[
(E_1,E_2), (E_1,E_3), (E_1,E_4), (E_2,E_3), (E_2,E_4), (E_3,E_4), 
\]
we have the following expression of the defect:    

\begin{align*}
\rho ({\mathcal S}) & = \sum _{i=1} ^4 \dim E_i - 2\dim H \\
& = \frac{1}{3} \sum _{1\leq i<j\leq 4} 
(\dim E_i + \dim E_j -\dim H ) \\
& = \frac{1}{3} \sum _{1\leq i<j\leq 4} 
(\dim (E_i \cap E_j) - \dim ((E_i + E_j)^{\perp})).  
\end{align*}

\noindent
{\bf Definition} Let ${\mathcal S} = (H ; E_1, E_2, E_3, E_4)$ 
be a system of four subspaces. For any distinct $i,j = 1,2,3,4$, 
define an adding operator 
\[
A_{ij} : E_i \oplus E_j \ni (x,y) \rightarrow x+y \in H. 
\] 
Then 
\[
\Ker A_{ij} = \{(x,-x) \in E_i \oplus E_j ; x \in E_i \cap E_j \}
\]
and 
\[
\Im A_{ij} = E_i + E_j .
\]
We say ${\mathcal S} = (H ; E_1, E_2, E_3, E_4)$ is a
{\it Fredholm} system if  $A_{ij}$ is a Fredholm operator for any 
$i,j = 1,2,3,4$ with $i \not= j$. Then $\Im A_{ij} = E_i + E_j$ is closed 
and 
\[
\Index A_{ij} = \dim \Ker A_{ij}  - \dim \Ker A_{ij}^* 
= \dim (E_i \cap E_j) - \dim ((E_i + E_j)^{\perp}). 
\] 
T. Kato called the number 
 $dim (E_i \cap E_j) - \dim ((E_i + E_j)^{\perp})$ 
the index of the pair $E_i, E_j$ in (\cite{K};IV section 4).

\noindent
{\bf Definition} We say 
${\mathcal S} = (H ; E_1, E_2, E_3, E_4)$ is a
{\it quasi-Fredholm} system if $E_i \cap E_j$ and $(E_i + E_j)^{\perp}$ 
are finite-dimensional for any $i \not= j$. In the case we define 
the {\it defect} $\rho ({\mathcal S})$ of ${\mathcal S}$ by 

\begin{align*}
\rho ({\mathcal S}) & := 
\frac{1}{3} \sum _{1\leq i<j\leq 4} 
(\dim (E_i \cap E_j) - \dim (E_i + E_j)^{\perp})) \\
& = \frac{1}{3} \sum _{1\leq i<j\leq 4} 
(\dim (E_i \cap E_j) - \codim \overline{E_i + E_j}) 
\end{align*}
which coincides with the Gelfand-Ponomarev original defect if 
$H$ is finite-dimensional. 
Moreover, if ${\mathcal S}$ is a Fredholm system, then 
it is a quasi-Fredholm system and 
\[
\rho ({\mathcal S}) = \frac{1}{3} \sum _{1\leq i<j\leq 4} \Index A_{ij} . 
\]

\begin{prop}
Let $\mathcal{S}_T = (H;E_1,E_2,E_3,E_4)$ be a 
bounded operator system associated with a single operator
$T \in B(K)$. Then $\mathcal{S}_T$ is a Fredholm system 
if and only if $T$ and $T-I$ are Fredholm operators. If the 
condition is satisfied, then the defect is given by 
\[
\rho (\mathcal{S}_T) = \frac{1}{3}(\Index T + \Index (T-I))
\]
Similarly $\mathcal{S}_T$ is a quasi-Fredholm system 
if and only if $\Ker T$, $\Ker T^*$,  $\Ker (T-I)$ and 
$\Ker (T-I)^*$ are finte-dimensional.  If the condition is 
satisfied, 
then the defect is given by 
\[
\rho (\mathcal{S}_T) = \frac{1}{3}
(\dim \Ker T - \dim \Ker T^* + \dim \Ker (T-I) - \dim \Ker (T-I)^*)
\]
\label{proposition:single defect}
\end{prop}
\begin{proof}
It is clear that $E_i \cap E_j = 0$ and $E_i + E_j = H$ 
for $(i,j) = (1,2), (1,4),(2,4),(2,3)$. Since  
$\Ker A_{13} = E_1 \cap E_3 = \Ker T \oplus 0$  and 
$(\Im A_{13})^{\perp} = (E_1 + E_3)^{\perp} = (K \oplus \Im T)^{\perp}$, 
they are  finite-dimensional 
if and only if $\Ker T$ and $(\Im T)^{\perp} = \Ker T^*$ 
are finite-dimensional. 
And $\Im A_{13}$ is closed if and only if $\Im T$ is closed. 
We transform $E_3$ and $E_4$ by an invertible operator 
$
R = 
\begin{pmatrix}
     I  & 0 \\
     -I & I 
\end{pmatrix}
\in B(H) = B(K\oplus K)
$
, then $R(E_3) = \{(x,(T-I)x) \in  \ K \oplus K ; x \in K\}$ 
and $R(E_4) = K \oplus 0$. Hence 
$R(E_3 \cap E_4) = \Ker (T-I) \oplus 0$  and 
$R(E_3 + E_4) = K \oplus \Im (T-I)$. Then 

\begin{align*}
\dim ((E_3 + E_4)^{\perp}) 
& = \codim \ \overline{E_3 + E_4} \\
& = \codim \ \overline{R(E_3) + R(E_4)}) 
=dim ((R(E_3 + E_4))^{\perp}) 
\end{align*}

Thus $E_3 \cap E_4$ and $(E_3 + E_4)^{\perp}$ 
 are  finite-dimensional 
if and only if $\Ker (T-I)$ and $(\Im (T-I))^{\perp} = \Ker (T-I)^*$ 
are finite-dimensional. 
And $\Im A_{13} = E_3 + E_4$ is closed if and only if $Im (T-I)$ is closed.
It follows the desired conclusion.
\end{proof}

We shall show that the defect could  have a fractional value.
 
\noindent
{\bf Example. }Let $S$ be a 
unilateral shift on $K = \ell^2(\mathbb N)$. Then 
the operator system $\mathcal{S}_S$ is an indecomposable.  
It is not a Fredholm system but a quasi-Fredholm system 
and  $\rho (\mathcal{S}_S) = -\frac{1}{3}$. 
The operator system $\mathcal{S}_{S + \frac{1}{2}I}$ is a 
Fredholm system 
and  $\rho (\mathcal{S}_{S + \frac{1}{2}I}) = - \frac{2}{3}$. Moreover 
$(\mathcal S_{T+\alpha I})_{\alpha\in \mathbb C}$
is  uncountable family of indecomposable ,
quasi-Fredholm systems.  Fredholm systems among them and 
their defect are given by 
\[
 \ \rho(\mathcal S_{S+\alpha I}) 
= 
\begin{cases}
- \frac{2}{3}, & \ 
               (\vert\alpha\vert<1 \text{ and } \vert\alpha-1\vert<1) \\
-\frac{1}{3},  & \ (\vert\alpha\vert<1 \text{ and } \vert\alpha-1\vert>1)
           \text{ or } (\vert\alpha\vert>1 \text{ and } \vert\alpha-1\vert<1) \\

0, & \ (\vert\alpha\vert>1 \text{ and } \vert\alpha-1\vert>1). 
\end{cases}
\]

\begin{cor}
Let $\mathcal{S}_T = (H;E_1,E_2,E_3,E_4)$ be a 
bounded operator system associated with a single operator
$T \in B(K)$. If $\mathcal{S}_T$ is a Fredholm system, then 
$\mathcal{S}_{T^*}$ is a Fredholm system and 
$\rho (\mathcal{S}_{T^*}) =  - \rho (\mathcal{S}_T)$. 
Similarly If  $\mathcal{S}_T$ is a quasi-Fredholm system 
then $\mathcal{S}_{T^*}$ is a quasi-Fredholm system and 
$\rho (\mathcal{S}_{T^*}) =  - \rho (\mathcal{S}_T)$. 
\end{cor}
\begin{proof}
Use the fact that $T$ is Fredholm if and only if $T^*$ is a 
Fredholm, and then $\Index T^* = - \Index T$. 
\end{proof}

\begin{prop}
Let $\mathcal{S}=  (H;E_1,E_2,E_3,E_4)$ be a system 
of four subspaces.  If $\mathcal{S}$ is a Fredholm system, then 
the orthogonal complement $\mathcal{S} ^{\perp} 
= (H;E_1^{\perp},E_2^{\perp},E_3^{\perp},E_4^{\perp}) $ 
is a Fredholm system and 
$\rho (\mathcal{S} ^{\perp}) =  - \rho (\mathcal{S})$. 
Similarly if  $\mathcal{S}$ is a quasi-Fredholm system 
then $\mathcal{S} ^{\perp}$ is a quasi-Fredholm system and 
$\rho (\mathcal{S} ^{\perp}) =  - \rho (\mathcal{S})$. 
\end{prop}
\begin{proof}Recall elementary facts that 
$E_i^{\perp} \cap E_j^{\perp} = (E_i + E_j)^{\perp}$ and 
$(E_i^{\perp} + E_j^{\perp})^{\perp} = E_i \cap E_j$. 
The only non-trivial thing is to know  that 
$E_i + E_j$ is closed if and only if 
$E_i^{\perp} + E_j^{\perp}$ is closed, see,  for example, 
(\cite{K};IV Theorem 4.8).  
\end{proof}

\noindent
{\bf Example}. For $\gamma \in \mathbb C $ with 
$|\gamma | \geq 1$, let  
$\mathcal{S}_{\gamma}=  (H;E_1,E_2,E_3,E_4)$ 
be an exotic system of four subspaces in 
Theorem \ref{thm:exotic examples}.  
Then $\mathcal{S}_{\gamma}$ is a quasi-Fredholm system 
and 
\[
\rho (\mathcal{S}_{\gamma}) 
= \frac{1}{3}(\Index A_{13} + \Index A_{23} + \Index A_{34})
= \frac{1}{3}(1+1+1) = 1.
\] 
In fact, $E_1 \cap E_3 = {\mathbb C}(e_1,0,0,0)$, 
$E_2 \cap E_3 = {\mathbb C}(0,0,0,e_1)$ and 
$E_4 \cap E_3 = {\mathbb C}(a,0,a,0)$, where 
$a = (\gamma ^{n-1})_n \in L =\ell^{2}(\mathbb{N})$. 
All the other terms are  zeros. 

\medskip
\noindent
{\bf Definition}. Let $\mathcal{S}=  (H;E_1,E_2,E_3,E_4)$ be a system 
of four subspaces. We say that $\mathcal{S}$ is {\it non-degenerate} 
if $E_i + E_j = H$ and $E_i \cap E_j = 0$ for  $i \not= j$.  Then 
$\mathcal{S}$ is  clearly a Fredholm system with the defect 
$\rho (\mathcal{S}) = 0$. Thus the defect measures the failure from 
being non-degenerate.

\begin{prop}
Let $\mathcal{S}=  (H;E_1,E_2,E_3,E_4)$ be a system 
of four subspaces.  Then $\mathcal{S}$ is non-degenerate 
if and only if $\mathcal{S}^{\perp}$ is non-degenerate.
\end{prop} 
\begin{proof}
It follows from the fact that $E_i + E_j = H$  if and only if 
$E_i^{\perp}  \cap E_j^{\perp} = 0$.  
\end{proof}

\begin{prop}
 Let $\mathcal S_{T,S}$ be a bounded operator system. 
Then $\mathcal S_{T,S}$ is a Fredholm system if and only if 
$S,T$ and $ST-I$ are Fredholm operators. And if the condition is satisfied, 
then
\[
\rho(\mathcal S_{T,S}) 
= \frac{1}{3}(\Index T + \Index S + \Index (ST-I)).
\]
\end{prop}
\begin{proof}
It is clear that $E_i \cap E_j = 0$ and $E_i + E_j = H$ 
for $(i,j) = (1,2), (1,4),(2,3)$. Since  
$\Ker A_{13} = E_1 \cap E_3 = \Ker T \oplus 0$  and 
$(\Im A_{13})^{\perp} = (E_1 + E_3)^{\perp} = (K_1 \oplus \Im T)^{\perp}$, 
they are  finite-dimensional 
if and only if $\Ker T$ and $(\Im T)^{\perp} = \Ker T^*$ 
are finite-dimensional. 
And $\Im A_{13}$ is closed if and only if $\Im T$ is closed. 
Similarly 
$\Ker A_{24} = E_2 \cap E_4= 0 \oplus \Ker S$  and 
$(\Im A_{24})^{\perp} = (E_2 + E_4)^{\perp} = (\Im S \oplus K_2)^{\perp}$.
Hence they are  finite-dimensional 
if and only if $\Ker S$ and $(\Im S)^{\perp} = \Ker S^*$ 
are finite-dimensional. 
And $\Im A_{24}$ is closed if and only if $\Im S$ is closed. 
Nextly, 
\[
\Ker A_{34} = E_3 \cap E_4 = \{(x,Tx) \in K_1 \oplus K_2 ; x \in \Ker (ST-I) \}.\]
\[
\Im A_{34} = \{
\begin{pmatrix}
 x + Sy \\
 Tx + y   
\end{pmatrix}
; x \in K_1, y \in K_2 \} = 
\begin{pmatrix}
     I & S\\
     T & I
\end{pmatrix} 
\begin{pmatrix}
 x  \\
 y  
\end{pmatrix}
; x \in K_1, y \in K_2 \}.
\]
Multiplying invertible operator matrices from both sides, we have 
\[ 
\begin{pmatrix}
     I & -S\\
     0 & I
\end{pmatrix} 
\begin{pmatrix}
     I & S\\
     T & I
\end{pmatrix} \begin{pmatrix}
     I & 0\\
     -T & I
\end{pmatrix} 
=
\begin{pmatrix}
     I-ST & 0\\
     0 & I
\end{pmatrix} 
.
\]
Hence $\Im A_{34}$ is closed if and only if $\Im (ST-I)$ is closed, 
and $(\Im A_{34})^{\perp}$ is finite-dimensional if and only if 
$(\Im (ST-I))^{\perp}$ is finite-dimensional.  
Now it is easy to see the desired conclusons.
\end{proof}

Let $\mathcal S$ and $\mathcal S^{\prime}$ be two quasi-Fredholm systems 
of four subspaces.  Then it is evident that 
$\mathcal S\oplus \mathcal S^{\prime}$ is also a quasi-Fredholm system 
and 
\[
\rho(\mathcal S\oplus \mathcal S^{\prime}) 
=  \rho(\mathcal S) + \rho(\mathcal S^{\prime}). 
\]
Therefore we should investigate the possible values of the defect 
for  indecomposable systems.

\begin{thm}
The set of the possible values of the defect of indecomposable systems of 
four subspaces is  exactly $\mathbb Z/3$
\end{thm}
\begin{proof}
Let $S$ be a unilateral shift on $L=\ell^{2}(\mathbb N)$.
Let $K = L \otimes  \mathbb C^n$ and $H = K \oplus K$. 
For a positive integer $n$, put 
\[
 V = \begin{pmatrix} 
S& 0 & 0 & \cdots & 0  \\ 
I & S & 0 & \cdots & 0 \\
0 & I & S & \cdots & 0 \\
\vdots & \vdots & \ddots & \ddots & \vdots \\
0 & 0 & \cdots & I & S 
 \end{pmatrix}
\in M_{n}(\mathbb C)\otimes B(L) = B(K). 
\]
Let ${\mathcal S}_{V}=(H;E_1,E_2,E_3,E_4)$ be the operator system 
associated with a single operator $V$. 
We shall show that ${\mathcal S}_{V}$ is indecomposable.
Let $T = (T_{ij})_{ij} \in B(K)$ be an idempotent which commutes with $V$. 
It is enough to show that $T =  0$ or $T = I$.  

Since $VT = TV$, we have 
\[
ST_{11} = T_{11}S + T_{12}, \ \dots \  ST_{1(n-1)} = T_{1(n-1)}S + T_{1n}, 
\ T_{1n}S = ST_{1n}. 
\]
By the Kleinecke-Shirokov theorem, $T_{1n}$ is a quasinilpotent. 
Since $T_{1n}$ commutes with a unilateral shift $S$, $T_{1n}$ is a Toeplitz 
operator.  Then $\|T_{1n}\| = r(T_{1n}) = 0$. Thus $T_{1n} = 0$ by 
\cite{Ha3}.  
Inductively we can show that $T_{12}=T_{13}=\cdots=T_{1n}=0$. 
Similar argument shows that  $T$ is a lower triangular operator matrix, 
i.e., $T_{ij}=0$ for $i<j$. Since $T^{2}=T,$ we have 
$T_{ii}^{2}=T_{ii}$ for $ i=1,\cdots,n$. The diagonal of  $VT = TV$ shows 
that each $T_{ii}$ commutes with a unilatral shift $S$. This implies that 
$T_{ii} = 0$ or $I$ as in Lemma \ref{lemma:Toeplitz-idempotent}. 

\noindent
(i)the case that $T_{11} = 0$: The $2$-$1$th component of  $VT = TV$ shows 
that $T_{22} = ST_{21}-T_{21}S$. Hence $T_{22}$ cannot be $I$. Thus 
$T_{22} = 0$.  Similarly we can show that $T_{ii} = 0$ for $i = 1, \dots, n$. 
Thus the diagonal of operator matrix $T$ is zero.  Furthermore 
$T$ is a lower triangular operator matrix and idempotent.  Hence $T = O$.  

\noindent
(ii) the case that $T_{11} = I$: Considering $I - T$ instead of $T$, we can 
use the case (i) and shows that $T = I$.  Therefore ${\mathcal S}_{V}$ is indecomposable. 

The defect is given by 

\begin{align*}
\rho ({\mathcal S}_{V}) 
& = \frac{1}{3}
(\dim \Ker V - \dim \Ker V^* + \dim \Ker (V-I) - \dim \Ker (V-I)^*) \\
& = \frac{1}{3}(0-n+0-0) = \frac{-n}{3}.
\end{align*} 

In fact, 
\[
\Ker V^* = \{(a,-S^*a,(-S^*)^2a,\dots,(-S^*)^{n-1}) 
\in (\ell^{2}(\mathbb N))^n ; 
a \in \Ker S^{*n} \}
\]
is $n$-dimensional.

Similarly  ${\mathcal S}_{V^*}$ is an indecomposable system 
with $\rho({\mathcal S}_{V^*}) = \frac{n}{3}$.  

For $n = 0$, consider an indecomposable system ${\mathcal S}_{S+3I}$ 
as in Example after Proposition \ref{proposition:single defect}. 
Then $\rho({\mathcal S}_{S+3I}) = 0$.

Therefore  the defect for indecomposable systems of 
four subspaces can take any value in  $\mathbb Z/3$. 
\end{proof}

\noindent
{\bf Remark}. Indecomposablity of the  system ${\mathcal S}_V$ can also be 
derived by Theorem 3.4 in \cite{JW}, although we give our direct proof.

\begin{cor}
For any $n\in \mathbb Z$ there exist uncountable family of indecomposable
systems ${\mathcal S}$ of four subspaces with the same defect 
$\rho({\mathcal S}) = \frac{n}{3}$.
\end{cor}

\begin{proof} For a positive integer $n$, 
consider a family $({\mathcal S}_{V+\alpha I})_{\alpha \in (0,1)}$ and 
$({\mathcal S}_{V^*+\alpha I})_{\alpha \in (0,1)}$
of bounded operator systems similarly as in the above theorem. 
Then any ${\mathcal S}_{V+\alpha I}$ is also indecomposable and 
\[
\rho({\mathcal S}_{V+\alpha I}) 
= \frac{1}{3}(0-n+0-0) = \frac{-n}{3} .
\]

If $\alpha \not= \beta$, then the spectrum 
$\sigma(V +\alpha I) \not= \sigma(V + \beta I)$ .
Since $V +\alpha I$ and $V + \beta I$ are not similar, 
${\mathcal S}_{V+\alpha I}$ and ${\mathcal S}_{V+\beta I}$ are 
not isomorphic each other.  

We also have $\rho({\mathcal S}_{V^*+\alpha I}) = \frac{n}{3}$.　
And they are not isomorphic each other.  

For $n = 0$, consider a family 
$({\mathcal S}_{S+3I + \alpha I})_{\alpha \in [0,1]}$ in Example after 
Proposition 8.5. 
They are indecomposable , not isomorphic each other and 
$\rho ({\mathcal S}_{S+3I + \alpha I}) = 0$.  
\end{proof}

\section{Coxeter functors}
In \cite{GP} Gelfand and Ponomarev introduced two functors 
$\Phi ^+$ and $\Phi ^-$ on the category of systems $\mathcal S$ 
of $n$ subspaces in finite-dimensional vector spaces. 
They used the functors $\Phi ^+$ and $\Phi ^-$ to give a complete 
classification of indecomposable systems of four subspaces with defect 
$\rho (\mathcal S) \not= 0$ in finite-dimensional vector spaces. If  
the defect $\rho (\mathcal S) < 0$, then there exists a positive integer 
$\ell$ such that $(\Phi ^+)^{{\ell}-1}(\mathcal S) \not= 0$ and 
$(\Phi ^+)^{\ell}(\mathcal S) = 0$.  Combining the facts that indecomposable 
systems $\mathcal T$ with $\Phi ^+(\mathcal T) = 0$ can be classified easily 
and that $\mathcal S$ is isomorphic to (and recovered as) 
$(\Phi ^-)^{{\ell}-1}(\Phi ^+)^{{\ell}-1}(\mathcal S)$, they 
provided a complete classification.  A similar argument 
holds for systems $\mathcal S$ with defect $\rho (\mathcal S) > 0$. 

In their argument the finiteness of dimension is used crucially. In fact 
if an indecomposable system ${\mathcal S}=(H;E_1,E_2,E_3,E_4)$ with 
$\dim H > 1$ satisfies that the defect $\rho (\mathcal S) < 0$, then 
$\Phi ^+(\mathcal S) = (H^+;E_1^+,E_2^+,E_3^+,E_4^+)$  has the property 
that $\dim H^+ < \dim H$. The property guarantees the existence of 
a positive integer 
$\ell$ such that $(\Phi ^+)^{\ell}(\mathcal S) = 0$. Although we can not expect
such an argument anymore in the case of infinite-dimensional space, these 
functors $\Phi ^+$ and $\Phi ^-$ are interesting on their own right. Therefore 
we shall extend these functors $\Phi ^+$ and $\Phi ^-$  on infinite-dimensional Hilbert spaces and  show that the Coxeter functors preserve the 
defect and indecomposability under certain conditions.  

\bigskip
\noindent
{\bf Definition}.(Coxeter functor $\Phi ^+$) Let ${\mathcal Sys}^n$ be the 
category of the systems of $n$ subspaces in Hilber spaces 
and homomorphisms. Let $\mathcal S = (H;E_1,\ldots ,E_n)$ be a system
of $n$ subspaces  in a Hilbert space $H$. 
Let $R := \oplus _{i=1}^n E_i$ and 
\[
\tau : R \ni x = (x_1,\dots,x_n) \longmapsto \tau(x) = \sum _{i=1}^n x_i \in H .\]
Define $\mathcal S^+ = (H^+;E_1^+,\dots ,E_n^+)$ by 
\[
H^+ := \Ker \tau  \text{ and }   E_k^+ := \{(x_1,\dots,x_n) \in H^+ ; x_k = 0 \}.\] 
Let ${\mathcal T} = (K;F_1, \dots , F_n)$  
be  another system of $n$ subspaces in a Hilbert space $K$ and  
$\varphi : {\mathcal S} \rightarrow {\mathcal T}$ be a 
homomorphism. Since $\varphi : H \rightarrow K$ is a bounded linear operator 
with $\varphi(E_i) \subset F_i$, 
we can define a bounded linear operator 
$\varphi^+: H^+ \rightarrow K^+$ by 
$\varphi^+(x_1,\dots,x_n) = (\varphi(x_1),\dots,\varphi(x_n))$. 
Since $\varphi^+(E_i^+) \subset F_i^+$, $\varphi^+$ define a 
homomorphism $\varphi ^+ : {\mathcal S}^+ \rightarrow {\mathcal T}^+$. 
Thus we can introduce a covariant 
functor $\Phi ^+: {\mathcal Sys}^n \rightarrow 
{\mathcal Sys}^n$ by 
\[
\Phi ^+({\mathcal S}) = {\mathcal S}^+
\text{ and } \Phi ^+(\varphi) =  \varphi^+. 
\]

\bigskip
\noindent
{\bf Example}.If ${\mathcal S} = (\mathbb C; \mathbb C , \mathbb C,\mathbb C)$, then ${\mathcal S}^+ \cong (\mathbb C^2; \mathbb C (1,0), 
\mathbb C (0,1),\mathbb C (1,1))$.

\begin{lemma}
Let ${\mathcal S} = (H;E_1,E_2,E_3,E_4)$ be 
a system of four subspaces and  consider 
${\mathcal S}^+  = (H^+;E_1^+,E_2^+,E_3^+,E_4^+)$ . 
Then 
\[
E_1^+ \cap E_2^+ = \{(0,0,a,-a) \in \oplus _{i=1}^4 E_i ; a \in E_3 \cap E_4 \} .
\]
In particular, we have $\dim E_1^+ \cap E_2^+ = \dim E_3 \cap E_4$. 
Same formulae hold under permutation of subspaces.
\label{lemma:nondegenerate1}
\end{lemma}
\begin{proof}
Let $x = (x_1,x_2,x_3,x_4) \in E_1^+ \cap E_2^+$, then $x_1 = x_2 = 0$. 
Since $x \in H^+$, $\tau (x) = x_3 + x_4 = 0$. 
Thus $a := x_3 = - x_4  \in E_3 \cap E_4$ and $x = (0,0,a,-a)$. The converse 
inclusion is clear.    
\end{proof}

\begin{lemma}
Let ${\mathcal S} = (H;E_1,E_2,E_3,E_4)$ be 
a system of four subspaces and  consider 
${\mathcal S}^+  = (H^+;E_1^+,E_2^+,E_3^+,E_4^+)$ . 
If $E_3 \cap E_4 = 0$ and $E_3 + E_4 = H$, then $E_1^+ + E_2^+ = H^+$. 
Same formulae  hold under permutation of subspaces.
\label{lemma:nondegenerate2}
\end{lemma}
\begin{proof}
Let $z = (z_1,z_2,z_3,z_4) \in H^+$. Put $y_1 := z_1$ and 
$x_2 := z_2$.  Since  $E_3 + E_4 = H$, there exist $y_3 \in E_3$ and 
$y_4 \in E_4$ such that $-y_1 = y_3 + y_4$.  Since $y_1 + y_3 + y_4 = 0$, 
$y := (y_1,0,y_3,y_4) \in H^+$.  Similarly 
there exist $x_3 \in E_3$ and 
$x_4 \in E_4$ such that $-x_2 = x_3 + x_4$, so that 
$x := (0,x_2,x_3,x_4) \in H^+$. 

Since $z \in H^+$, $z_1 + z_2 + z_3 + z_4 = 0$. Hence 
\begin{align*}
z_3 + z_4  & = -z_1 -z_2 = -y_1 -y_2 \\
           & = (y_3 + y_4) + (x_3 + x_4) 
           = (x_3 + y_3) + (x_4 + y_4) \in E_3 + E_4. 
\end{align*} 
Because $E_3 \cap E_4 = 0$, we have $z_3 = x_3 + y_3$ and $z_4 = x_4 + y_4$. 
Therefore $z = x + y \in E_3^+ + E_4^+$. 
\end{proof}

\bigskip
\noindent
{\bf Example}. Let ${\mathcal S}_{S,T} = (H;E_1,E_2,E_3,E_4)$ be 
a bounded operator system.  Combining the preceding two lemmas 
Lemma \ref{lemma:nondegenerate1} and Lemma \ref{lemma:nondegenerate2} 
with a characterization of bounded operator systems in Corollary 
\ref{cor:bounded operator system}, we have that  
${\mathcal S}^+  = (H^+;E_1^+,E_2^+,E_3^+,E_4^+)$ is a bounded 
operator system up to permutation of subspaces. More precisely,  
$(H^+;E_3^+,E_4^+,$ $E_1^+,E_2^+)$ is a bounded operator system.

\bigskip
Let $0\oplus E_i \oplus 0
:= 0 \oplus \dots \oplus 0 \oplus E_i 
\oplus 0 \oplus \dots \oplus 0 \subset R$ 
and $q_i \in B(R)$ be the projection onto 
$0\oplus E_i \oplus 0$.  Let $\imath _+ : H^+ \rightarrow R$ be a 
canonical embedding. Then we have an exact sequence: 
\[
0 \longrightarrow H^+ \overset{\imath _+}{\longrightarrow} R
\overset{\tau}{\longrightarrow} H
\]
Furthermore we have 
\[
\Ker \tau q_i = \Ker q_i, \ 
E_i = \Im \tau q_i = \overline{\Im \tau q_i} \ 
\text{ and } \ E_i^+ = \Ker q_i \imath _+ .
\] 
These properties characterize 
${\mathcal S}^+  = (H^+;E_1^+,E_2^+,E_3^+,E_4^+)$ . 

\begin{prop}Let $X, Y$ and $Z$ be Hilbert spaces and 
$T:X \rightarrow Y$ and $S:Y \rightarrow Z$ be bounded 
linear maps. Suppose that a sequence
\[
0 \longrightarrow X \overset{T}{\longrightarrow} Y
\overset{S}{\longrightarrow} Z.
\]
is exact. Let $p_1,...,p_n \in B(Y)$ be projections  with 
$\sum _i p_i = I$ and $p_ip_j = 0$ for $i \not= j$.  
Furthermore we assume that 
\[
\Ker Sp_i = \Ker p_i  \text{ and  } 
\Im Sp_i \text { is closed in } Z .
\]
Let $E_i := \Im Sp_i \subset Z $
and 
$E_i' := \Ker p_iT \subset X$.  Define 
$\mathcal S = (Z;E_1,\dots ,E_n)$ and 
$\mathcal S' = (X;E_1',\dots ,E_n')$. 
Then  $\mathcal S' \cong \Phi ^+({\mathcal S})$
\label{prop:characterize Phi+}
\end{prop}
\begin{proof}
Consider the restriction 
$S_i:= S|_{\Im p_i} : \Im p_i \rightarrow \Im Sp_i$. 
Since  $\Ker Sp_i = \Ker p_i$, $S_i$ is one to one. 
Because  $\Im Sp_i$ is closed, $\Im Sp_i$ is complete. 
Therefore $S_i$ is an invertible operator by open 
mapping theorem. 
Define $\varphi : Y = \oplus _{i=1}^n \Im p_i \rightarrow \oplus _{i=1}^n E_i$
by $\varphi((y_i)_i ) = (S_i(y_i))_i$ for $(y_i)_i \in \oplus _{i=1}^n \Im p_i$.Then $\varphi$ is an invertible operator. 
Consider  $\tau : \oplus _{i=1}^n E_i \rightarrow Z$ given 
$\tau ((z_i)_i) = \sum _{i=1}^n z_i$.  Let $Z^+ = \Ker \tau$ and 
$\imath _+ : Z^+ \rightarrow \oplus _{i=1}^n E_i$ be a canonical embedding. 
Then $\tau \varphi = S$. Define $\psi : X \rightarrow Z^+$  by 
$\psi (x) = \varphi T(x)$ for $x \in X$. The map $\psi$ is well-defined, 
because $\tau(\psi (x)) = \tau (\varphi T(x)) = ST(x) = 0$.  Then the following diagram 
\[
\begin{CD}
0 @>>> X @>T>> Y @>S>> Z \\
@.  @V \psi VV  @V \varphi VV @V id_Z VV \\
0 @>>> Z^+  @>\imath _+>> \oplus _{i=1}^n E_i 
@>\tau>>Z
\end{CD}
\]
is commutative. Furthermore maps $\psi$  and $\varphi$ are invertible 
operators. Let $q_i \in B(\oplus _{i=1}^n E_i)$ be a projection onto 
$0\oplus E_i \oplus 0$.  Then $q_i = \varphi p_i \varphi^{-1}$, 
$E_i^+ = \Ker (q_i \imath _+)$ and $E_i' = \Ker (p_i T)$.  Therefore 
$\psi(E_i') = E_i^+$.  Thus 
$\psi : \mathcal S' \rightarrow \Phi ^+({\mathcal S})$ 
is a desired isomorphism. 
\end{proof}

\bigskip
\noindent
{\bf Definition}.(Coxeter functor $\Phi ^-$) In \cite{GP} Gelfand and 
Ponomarev introduced a dual functor  $\Phi ^-$ using quotients of 
vector spaces.  If $H$ is a Hilbert space and $K$ a subspace of $H$, 
then it is convenient to identify the quotient space $H/K$  with the 
orthogonal complement $K^{\perp}$.  Therefore we shall generalize their 
functor $\Phi ^-$ in terms of orthogonal complements instead of quotients 
in our case of Hilbert spaces. 
Let $\mathcal S = (H;E_1,\ldots ,E_n)$ be a system
of $n$ subspaces  in a Hilbert space $H$. Let $e_i^{\perp} \in B(H)$ 
be the projection onto $E_i^{\perp} \subset H$.  
Let $Q := \oplus _{i=1}^n E_i^{\perp}$ and 
\[
\mu: H \ni x \longmapsto \mu(x) = (e_1^{\perp}x,\dots,e_n^{\perp}x) \in Q .
\]
Then $\mu ^* : Q \rightarrow H$ is given by 
$\mu ^*(y_1,\dots,y_n) = \sum _{i=1}^n y_i$. Define  
$H^- := \Ker \mu ^* \subset Q$.  Let  $\imath _- : H^- \rightarrow Q$ 
be a canonical embedding.  Then 
$q_- := \imath _-^*: Q \rightarrow H^-$ is the projection. 
Let $0\oplus E_i^{\perp} \oplus 0
:= 0 \oplus \dots 0 \oplus E_i^{\perp} \oplus 0 \dots \oplus 0 \subset Q$ 
and $r_i \in B(Q)$ be the  projection onto 
$0\oplus E_i^{\perp} \oplus 0$.
Define $\mathcal S^- = (H^-;E_1^-,\dots ,E_n^-)$ by 
\[
E_i^- := \overline{q_- (0\oplus E_i^{\perp} \oplus 0)} 
= \overline{\Im q_- r_i} \subset H^- . 
\] 
We note that 
\[
H^- := \Ker \mu ^* = Q \cap (\Im \mu)^{\perp} 
\cong Q/\overline{\Im \mu}. 
\]
We have an exact sequence 
\[
0 \longrightarrow H^- \overset{\imath _-}{\longrightarrow} Q
\overset{\mu^*}{\longrightarrow} H
\]
and a sequence 
\[
H \overset{\mu}{\longrightarrow} Q
\overset{q_-}{\longrightarrow} H^- \longrightarrow 0 , 
\]
satisfying that $\overline{\Im \mu} = \Ker q_-$ and $q_-$ is onto. 
Thus it is easy to see that our definition of 
$\mathcal S^- = (H^-;E_1^-,\dots ,E_n^-)$ coincides with 
the original one by Gelfand and Ponomarev up to isomorphism in 
the case of finite-dimensional spaces. 

Define $\Phi ^-(\mathcal S) := \mathcal S^- = (H^-;E_1^-,\dots ,E_n^-)$. 
Then there is a relation  between ${\mathcal S}^+$ and ${\mathcal S}^-$.  
We recall some elementary facts first.  

\begin{lemma}
Let $H$ and $K$ be Hilbert spaces and $M$ a closed subspace of $H$. Let 
$T: H \rightarrow K$ be a bounded operator. Consider $T^*: K \rightarrow H$. 
Then $\overline{T(M^{\perp})} = ((T^*)^{-1}(M))^{\perp} \subset K$. 
\end{lemma}

\begin{lemma}
Let $L$ be a Hilbert space and $M$, $K$  closed subspaces of $L$. 
Let $P_K \in B(L)$ be the projection onto $K$. Then 
$\overline{P_K(M^{\perp})} = K \cap (K \cap M)^{\perp}$. 
\end{lemma}
\begin{proof}
By the preceding lemma, 
\[
(\overline{P_K(M^{\perp})})^{\perp} = P_K^{-1}(M) = \{x \in L ; P_Kx \in M \}. 
\]
Decompose $x \in L$ such that  
$ x = x_1 + x _2$ with $x_1 \in K$,  $x_2 \in K^{\perp}$. 
Then $P_Kx \in M$ if and only if $x_1 \in M$.  Therefore  
$(\overline{P_K(M^{\perp})})^{\perp} = (K \cap M)  + K^{\perp}$. Thus 
$\overline{P_K(M^{\perp})} = K \cap (K \cap M)^{\perp}$. 
\end{proof}

\begin{prop}
Let $\mathcal S = (H;E_1,\ldots ,E_n)$ be a system
of $n$ subspaces  in a Hilbert space $H$. Then we have 
\[
\Phi ^-(\mathcal S) = \Phi ^{\perp}\Phi ^+\Phi ^{\perp}(\mathcal S) .
\]
\end{prop}
\begin{proof}
Since $\Phi ^{\perp}(\mathcal S) = (H;E_1^{\perp},\ldots ,E_n^{\perp})$, 
we have 
\[
\Phi ^+\Phi ^{\perp}(\mathcal S) 
= (H';(E_1^{\perp })^+,\ldots ,(E_n^{\perp})^+), 
\]
where $H' = \{(y_1,\dots,y_n) \in \oplus _{i=1}^n E_i^{\perp} ; 
y_1 + \dots + y_n = 0 \}$. Therefore we have  $H' = H^-$. 

Applying the preceding Lemma by putting 
$L = \oplus _{i=1}^n E_i^{\perp}$, 
$M = \{(y_1,\dots,y_n) \in L ; y_k = 0 \}$ and $K = H^- \subset L$, we have 
\[
E_k^- = \overline{q_- (0\oplus E_k^{\perp} \oplus 0)} 
= \overline{P_K(M^{\perp})} = K \cap (K \cap M)^{\perp} 
= H^- \cap ((E_k^{\perp})^+)^{\perp} .
\]
Therefore $(E_k^-)^{\perp} = (E_k^{\perp})^+$ in $H^-$.  
Hence $\Phi ^{\perp}\Phi ^-(\mathcal S)  = \Phi ^+\Phi ^{\perp}(\mathcal S)$. 
This implies the conclusion. 
\end{proof}  

Let $\mathcal S = (H;E_1,\ldots ,E_n)$ be a system
of $n$ subspaces  in a Hilbert space $H$ and 
${\mathcal T} = (K;F_1, \dots , F_n)$  
be  another system of $n$ subspaces in a Hilbert space $K$.  
Let $\varphi : {\mathcal S} \rightarrow {\mathcal T}$ be a 
homomorphism, i.e., $\varphi : H \rightarrow K$ is a bounded linear operator 
with $\varphi(E_i) \subset F_i$. 
Define $\varphi ^- : \Phi ^-({\mathcal S}) \rightarrow 
\Phi ^-({\mathcal T})$ by 
\[
\varphi ^- := \Phi ^{\perp}\Phi ^+\Phi ^{\perp} (\varphi) .
\]
Thus we can introduce a covariant 
functor $\Phi ^-: {\mathcal Sys}^n \rightarrow 
{\mathcal Sys}^n$ by 
\[
\Phi ^-({\mathcal S}) = {\mathcal S}^-
\text{ and } \Phi ^-(\varphi) =  \varphi^-. 
\]

\bigskip
\noindent
{\bf Remark}. 
Let $\mathcal S = (H;E_1,\ldots ,E_n)$ be a system
of $n$ subspaces  in a Hilbert space $H$. 
Let $R := \oplus _{i=1}^n E_i$ and 
$\tau : R \rightarrow H$ is given by 
$\tau (x) = \sum _{i=1}^n x_i$. Let $H^0 := \Ker \tau$ 
and $q_0 : R \rightarrow H^0$ be the canonical projection. 
Define $E_k^0 :=  \overline{q_0 (0\oplus E_k \oplus 0)}$. 
Let $\mathcal S^0 := (H^0;E_1^0,\dots ,E_n^0)$ and 
$\Phi ^0({\mathcal S}) = \mathcal S^0$.  
Then we have 
\[
\Phi ^+(\mathcal S) = \Phi ^{\perp}\Phi ^0(\mathcal S)  
\text{ and } 
\Phi ^- (\mathcal S)= \Phi ^0 \Phi ^{\perp}(\mathcal S) .
\]
Furthermore 
\[
\Phi ^- \Phi ^+(\mathcal S) = (\Phi ^0)^2(\mathcal S)
\text{ and }
\Phi ^+ \Phi ^-(\mathcal S) 
= \Phi ^{\perp}(\Phi ^0)^2\Phi ^{\perp}(\mathcal S). 
\]
Suppose that $H$ is finite-dimensional. Then 
\[
\dim H^0 = \dim \Ker \tau = \dim  R -\dim \Im \tau 
= \sum _i \dim E_i - \dim (\sum _i E_i)
\]
In particular, if  $\mathcal{S}=(H;E_1,E_2,E_3,E_4)$ is an 
indecomposable system of  four subspaces with $\dim H \geq 2$, 
then  $\dim H^0 = \sum _i \dim E_i - \dim H$ and the defect 
\[
\rho (\mathcal{S}) = \sum _i \dim E_i - 2\dim H = \dim H^0 - \dim H .
\]

We shall characterize $\Phi ^- (\mathcal S)$. The 
following fact is useful: 
Let $H$ and $K$ be Hilbert spaces and  $T:H \rightarrow K$ be 
a bounded linear operator. Then $\Im T$ is closed in $K$ 
if and only if $\Im T^*$ is closed in $H$. 

\begin{prop}Let $U, V$ and $W$ be Hilbert spaces and 
$A:U \rightarrow V$ and $B:V \rightarrow W$ be bounded 
linear operators. Suppose that a sequence
\[
U \overset{A}{\longrightarrow} V
\overset{B}{\longrightarrow} W \longrightarrow 0
\]
is exact. Let $p_1,...,p_n \in B(V)$ be projections  with 
$\sum _i p_i = I$ and $p_ip_j = 0$ for $i \not= j$.  
Furthermore we assume that 
\[
\Im p_i A \text { is closed in } V \text{ and  } 
\Im p_i A  = \Im p_i . 
\]
Let $L_i' := \overline{\Im Bp_i}  \subset W $
and 
$L_i := \Ker p_iA \subset U$.  Define 
$\mathcal S = (U;L_1,\dots ,L_n)$ and 
$\mathcal S' = (W;L_1',\dots ,L_n')$. 
Then  $\mathcal S' \cong \Phi ^-({\mathcal S})$
\label{prop:characterize Phi-}
\end{prop}
\begin{proof}Since $\Im B = W$ is closed, $\Im B^* \subset V$ 
is also closed. Then 
\[
\Im B^* = (\Ker B)^{\perp} = (\Im A)^{\perp} = \Ker A^* 
\]
and $\Ker B^* = (\Im B)^{\perp} = W^{\perp} = 0$. 
Hence the dual sequence 
\[
0 \longrightarrow W \overset{B^*}{\longrightarrow} V
\overset{A^*}{\longrightarrow} U
\]
is exact. We shall apply Proposition \ref{prop:characterize Phi+}
by putting $X= W$, $Y = V$, $Z = U$, $T = B^*$ and $S = A^*$. 
We can check the assumption of the Proposition. In fact, 
\[
\Ker Sp_i = \Ker A^*p_i = (\Im p_iA)^{\perp} = (\Im p_i)^{\perp} 
= \Ker p_i ,
\]
and $\Im Sp_i = \Im A^*p_i = \Im(p_iA)^*$ is closed, because 
$\Im(p_iA)$ is closed.   
Let \[
E_i := \Im Sp_i = \Im(p_iA)^* 
= (\Ker p_iA)^{\perp} = (L_i)^{\perp} \subset U
\]
and 
\[
E_i' := \Ker p_iT = \Ker p_iB^* 
=(\Im Bp_i)^{\perp} =(L_i')^{\perp} \subset W .
\] 
Then 
$(X;E_1',\dots ,E_n') \cong \Phi ^+(Z;E_1,\dots ,E_n)$, that is, 
we have 
\[
(W;(L_1')^{\perp}, \dots, (L_n')^{\perp}) \cong 
\Phi ^+(U;(L_1)^{\perp}, \dots,  (L_n)^{\perp}).
\]
Thus  $({\mathcal S}')^{\perp} \cong \Phi ^+({\mathcal S}^{\perp})$. 
Hence 
\[
{\mathcal S}' \cong \Phi ^{\perp}\Phi ^+\Phi ^{\perp}({\mathcal S}) 
= \Phi ^- ({\mathcal S}) .
\]
\end{proof}

\begin{prop}
Let $\mathcal S$ and $\mathcal T$ be systems
of $n$ subspaces  in a Hilbert space $H$. Then we have 
$\Phi ^+({\mathcal S} \oplus {\mathcal T}) 
\cong \Phi ^+({\mathcal S}) \oplus \Phi ^+({\mathcal T})$,  

$\Phi ^-({\mathcal S} \oplus {\mathcal T}) 
\cong \Phi ^-({\mathcal S}) \oplus \Phi ^-({\mathcal T})$, 
and  
$\Phi ^{\perp}({\mathcal S} \oplus {\mathcal T}) 
\cong \Phi ^{\perp}({\mathcal S}) \oplus \Phi ^{\perp}({\mathcal T}).
$
\end{prop}
\begin{proof}It is straightforward to prove them.
\end{proof} 

\bigskip
\noindent
{\bf Definition}. 
Let $\mathcal S = (H;E_1,\ldots ,E_n)$ be a system
of $n$ subspaces  in a Hilbert space $H$. Then 
$\mathcal S$ is said to be {\it reduced from above} 
if for any $k = 1, \dots, n$ 
\[
\sum _{i \not= k} E_i = H .
\]
In particular we have $E_k \subset \sum _{i \not= k} E_i$. 
Similarly $\mathcal S$ is said to be {\it reduced from below} 
if for any $k = 1, \dots, n$ 
\[
\sum _{i \not= k} E_i^{\perp}= H .
\] 
In particular we have 
$E_k^{\perp} \subset \sum _{i \not= k} E_i^{\perp}$ and 
$ \cap _{i \not= k} E_i = 0$

It is evident taht  $\mathcal S \oplus \mathcal T$ is 
reduced from above if and only if both $\mathcal S$ and $\mathcal T$ are
reduced from above. Similarly $\mathcal S \oplus \mathcal T$ is 
reduced from below if and only if both $\mathcal S$ and $\mathcal T$ are
reduced from below.

\bigskip
\noindent
{\bf Example}.(1) Any bounded operator system is reduced 
from above and reduced from below. In fact $E_1 + E_2 = H$, 
$E_1 + E_4 = H$, $E_2 + E_4 = H$ and  
$E_1^{\perp} + E_2^{\perp} = H$, 
$E_1^{\perp} + E_4^{\perp}= H$, 
$E_2^{\perp} + E_4^{\perp} = H$. 

\noindent
(2)The exotic examples in section 10 are reduced 
from above and reduced from below.   

\medskip
We shall show a duality theorem between Coxeter functors 
$\Phi ^+$ and  $\Phi ^-$. 

\begin{thm}
Let $\mathcal S = (H;E_1,\ldots ,E_n)$ be a system
of $n$ subspaces  in a Hilbert space $H$. Suppose that  
$\mathcal S$ is reduced from above. Then we have 
\[
\Phi ^- \Phi ^+({\mathcal S}) \cong {\mathcal S} .
\]
\label{thm:duality1}
\end{thm}
\begin{proof}Let $R = \oplus _{i=1}^n E_i$. 
Consider a sequence 
\[
H^+ \overset{\imath _+}{\longrightarrow} R
\overset{\tau}{\longrightarrow} H {\longrightarrow} 0 .
\]
Since $\mathcal S$ is reduced from above, 
$\Im \tau  = \sum _{i=1}^n  E_i = H$. Thus the 
above sequence is exact.  Let $p_i \in B(R)$ be the projection 
onto $0\oplus E_i \oplus 0$. We shall apply Proposition 
\ref{prop:characterize Phi-} by putting 
$U = H^+$, $V = R$, $W = H$, $A = \imath _+$ and $B = \tau$. 
We can check the assumption of the proposition. In fact, 
since $\mathcal S$ is reduced from above, for any 
$x_k \in  E_k$, there exist 
$x_i \in  E_i$ for $i \not= k$ such that 
$x_k = \sum_{i \not=k} -x_i$.  Then $\sum_{i=1}^n  x_i = 0$, 
that is, $x := (x_i)_i \in H^+$. Then 
\[
p_kA(x) = 0\oplus x_k \oplus 0 \in  0\oplus E_k \oplus 0 .
\]
Thus  $\Im p_kA = 0\oplus E_k \oplus 0 = \Im p_k$ and 
$\Im p_kA$ is closed. 
Therefore 
$ (W;L_1',\dots ,L_n') \cong  \Phi ^-(U;L_1,\dots ,L_n)$ . 
Since 
\[
L_k' = \overline{\Im Bp_k} = \overline{\Im \tau p_k} = E_k
\]
and 
\[
L_k = \Ker p_kA = \Ker p_k \imath _+ = E_k^+ , 
\]
we have 
\[
{\mathcal S} = (H;E_1,\dots,E_n) \cong  \Phi ^-(H^+;E_1^+,\dots, E_n^+) 
= \Phi ^- \Phi ^+ ({\mathcal S}) .
\]
\end{proof}

Similarly we have the follwoing: 
\begin{thm}
Let $\mathcal S = (H;E_1,\ldots ,E_n)$ be a system
of $n$ subspaces  in a Hilbert space $H$. Suppose that  
$\mathcal S$ is reduced from below. Then we have 
\[
\Phi ^+ \Phi ^-({\mathcal S}) \cong {\mathcal S} .
\]
\label{thm:duality2}
\end{thm}
\begin{proof}
If $\mathcal S$ is reduced from below, 
then ${\mathcal S}^{\perp}$ is reduced from above.  
Hence 
$\Phi ^- \Phi ^+({\mathcal S}^{\perp}) 
\cong {\mathcal S}^{\perp}$.  
Then \[
{\mathcal S} 
\cong \Phi ^{\perp}\Phi ^- \Phi ^+\Phi ^{\perp}({\mathcal S})
=\Phi ^{\perp} \Phi ^- \Phi ^{\perp} \Phi ^{\perp} \Phi ^+ \Phi ^{\perp}
({\mathcal S})
= \Phi ^+ \Phi ^- ({\mathcal S}) .
\]
\end{proof}　

\begin{prop}
Let $\mathcal S = (H;E_1,\ldots ,E_n)$ be a system
of $n$ subspaces  in a Hilbert space $H$.
Then $\Phi ^+({\mathcal S}) = 0$ if and only if 
for any $k = 1, \dots, n$ 
\[
E_k \cap (\sum _{i \not= k} E_i) = 0.
\] 
\end{prop}
\begin{proof}It is easy to see that 
$\Phi ^+({\mathcal S}) = 0$ if and only if
for any $x_i \in E_i$ with  $i = 1, \dots, n$
$\sum_i x_i = 0$ imples $x_1 = \dots = x_n = 0$. 
The latter condition is equal to that 
$E_k \cap (\sum _{i \not= k} E_i) = 0$ 
for any $k = 1, \dots, n$. 
\end{proof}
 
The above conditon $E_k \cap (\sum _{i \not= k} E_i) = 0$ 
for any $k = 1, \dots, n$ is something like an opposite of that 
$\mathcal S$ is reduced from above.

\begin{prop}
Let $\mathcal S = (H;E_1,\ldots ,E_n)$ be a system
of $n$ subspaces  in a Hilbert space $H$. 
Then $\Phi ^+({\mathcal S}) = 0$ and $\sum _{i=1}^n E_i$ is closed in $H$ 
if and only if  $(H;E_1,\ldots ,E_n,(\sum _{i=1}^n  E_i)^{\perp})$ is  
isomorphic to a system of direct sum decomposition, 
that is, there is an orthogonal direct sum decomposition 
$K = \oplus _{i=1}^{n+1} K_i$ 
of a Hilbert space $K$ and $(H;E_1,\ldots ,E_n,$ $(\sum _{i=1}^n  E_i)^{\perp})$ is isomorphic to a system $(K;K_1,\dots,K_{n+1})$, in particular ${\mathcal S}$ is  isomorphic to a commutative system. 
\end{prop}
\begin{proof} Assume that $\Phi ^+({\mathcal S}) = 0$ and 
$\sum _{i=1}^n E_i$ is closed in $H$. 
Let $E_{n+1} = (\sum _{i=1}^n  E_i)^{\perp}$. 
Let $R := \oplus _{i=1}^{n+1} E_i$ and 
$K_i := 0 \oplus \dots \oplus 0 \oplus E_i \oplus 0 \oplus \dots \oplus 0 
\subset R$. Define $\varphi: K \rightarrow H$ by $\varphi ((x_i)_i) = \sum _i x_i$. Then the  bounded operator $\varphi$ is onto, because  $\sum _{i=1}^n E_i$
is closed in $H$. Since $\Phi ^+({\mathcal S}) = 0$, $\varphi$ is one to one 
by the preceding proposition.  It is clear that $\varphi (K_i) = E_i$. 
Hence $(H;E_1,\dots,E_{n+1})$ is isomorphic to $(K;K_1,\dots,K_{n+1})$. 
The converse and the rest are trivial.
\end{proof}

\bigskip
\noindent
{\bf Example}. Let $T \in B(K)$ be a positive operator 
with dense range and $\Im T \not= K$. Let $H = K \oplus K$, 
$E_1 = K \oplus 0$ and $E_2 = \graph T$. Put 
${\mathcal S} = (H;E_1,E_2)$.  Then $\Phi ^+({\mathcal S}) = 0$ and 
$(E_1 + E_2)^{\perp} = 0$. But 
$(H;E_1,E_2,0)$ is not isomorphic to a system of direct sum decomposition. 
In fact $E_1 + E_2 = K \oplus \Im T$ is not closed. 

We also have the following:
\begin{prop}
Let $\mathcal S = (H;E_1,\ldots ,E_n)$ be a system
of $n$ subspaces  in a Hilbert space $H$
Then $\Phi ^-({\mathcal S}) = 0$ if and only if 
for any $k = 1, \dots, n$ 
\[
E_k^{\perp}  \cap (\sum _{i \not= k} E_i^{\perp} ) = 0.
\] 
\end{prop}

\begin{prop}
Let $\mathcal S = (H;E_1,\ldots ,E_n)$ be a system
of $n$ subspaces  in a Hilbert space $H$. If 
$\mathcal S$ is reduced from  above and $\mathcal S \not= 0$,  
then $\Phi ^+({\mathcal S}) \not= 0$. Similarly 
if $\mathcal S$ is reduced from below and $\mathcal S \not= 0$,  
then $\Phi ^-({\mathcal S}) \not= 0$.
\end{prop}
\begin{proof}Suppose that $E_i = 0$ for any $i = 1, \dots, n$. 
Then $H = \sum _{i=1}^{n-1}  E_i =0.$  This contradicts to that 
$\mathcal S \not= 0$.  Therefore $E_k \not= 0$ for some $k$. 
Since $\sum _{i \not= k} E_i = H$, for a non-zero $x_k \in E_k$, 
there exist $x_i \in E_k$ for $i \not= k$ such that 
$-x_k = \sum _{i \not= 0} x_i $. Therefore 
$x := (x_1,\dots,x_n) \in H^+ $ is non-zero, that is, 
$\Phi ^+({\mathcal S}) \not= 0$. The other is similarly proved. 
\label{prop:coxeter nonzero}
\end{proof}

\bigskip
\noindent
{\bf Remark}. By Proposition \ref{prop:n-1}, if a system of 
$n$ subspaces ${\mathcal S} = (H;E_1, \dots ,$ $E_n)$  
is indecomposable and $\dim H \geq 2$, then for any distinct $n$-$1$ subspaces 
$E_{i_1}, \dots, E_{i_{n-1}}$, we have that
\[
\bigcap _{k = 1}^{n-1} E_{i_k} = 0  \text{ and  } 
\bigvee _{k = 1}^{n-1} E_{i_k} = H, 
\]
that is, 
\[
\overline{\sum _{k = 1}^{n-1} E_{i_k}^{\perp} }= H   
\text{ and  } 
\overline{\sum_{k = 1}^{n-1} E_{i_k} } = H, 
\]
Unless $H$ is finite-dimensional, these conditions seems to be 
weaker than that $\mathcal S$ is reduced from  below and above. 

\bigskip
\noindent
{\bf Remark}. Let $\mathcal S = (H;E_1,\ldots ,E_n)$ be a system
of $n$ subspaces  in a Hilbert space $H$ and consider 
$\mathcal S^+ = (H^+;E_1^+,\dots ,E_n^+)$. Then 
for any distinct $n$-$1$ subspaces 
$E_{i_1}^+, \dots, E_{i_{n-1}}^+$, we have that
\[
\bigcap _{k = 1}^{n-1} E_{i_k}^+ = 0 .
\]
In fact, for example, let 
$(x_1,\dots,x_n) \in \cap _{k = 1}^{n-1} E_{k}^+ $. 
Then $x_1= x_2= \dots =x_{n-1} = 0$. Since 
$(x_1,\dots,x_n) \in H^+$, we have $\sum _{i=1}^n x_k = 0$. 
Hence $x_n = 0$. Thus  $\cap _{k = 1}^{n-1} E_{k}^+ = 0$. 

On the other hand  the above condition implies that 
\[
\overline{\sum _{k = 1}^{n-1} (E_{i_k}^+)^{\perp}} = H^+.  
\]
This condition is a little weaker than that $\mathcal S ^+$ is 
reduced from  below unless $H$ is finite dimensional.

Conider $\mathcal S ^- = \Phi ^{\perp}\Phi ^+\Phi ^{\perp}({\mathcal S})$ 
similarly.  Then we have 
\[
\overline{\sum _{k = 1}^{n-1} E_{i_k}^-} = H^-.  
\]
The condition is a little weaker than that $\mathcal S ^-$ is 
reduced from  above unless $H$ is finite dimensional.

\begin{thm}
Let $\mathcal S = (H;E_1,\ldots ,E_n)$ be a system
of $n$ subspaces  in a Hilbert space $H$. 
Suppose that  $\mathcal S$ is reduced from above  and 
$\mathcal S ^+ = \Phi ^+({\mathcal S})$ is reduced from below. 
If $\mathcal S$ is indecomposable, then $\Phi ^+({\mathcal S}) $ 
is also indecomposable.
\end{thm}
\begin{proof}On the contrary suppose  that $\mathcal S ^+$ were 
decomposable. Then there exist non-zero systems ${\mathcal T}_1$ and 
${\mathcal T}_2$ of $n$ subspaces such that 
$\mathcal S ^+ = {\mathcal T}_1 \oplus {\mathcal T}_2$. Since 
$\mathcal S$ is reduced from above, 
\[
{\mathcal S} \cong \Phi ^-\Phi ^+({\mathcal S}) 
= \Phi ^-({\mathcal T}_1) \oplus \Phi ^-({\mathcal T}_1) ,
\]
by a duality Theorem \ref{thm:duality1}. Since 
$\mathcal S ^+ = \Phi ^+({\mathcal S})$ is reduced from below, 
${\mathcal T}_1$ and ${\mathcal T}_2$ are also 
reduced from below. By another duality Theorem \ref{thm:duality2}, 
$\Phi ^+ \Phi ^- ({\mathcal T}_i) \cong {\mathcal T}_i$ for $i = 1,2$. 
Since  ${\mathcal T}_i \not= 0$, we have 
$\Phi ^- ({\mathcal T}_i) \not= 0$. (We could  use 
Propsition \ref{prop:coxeter nonzero} instead.)
This implies that $\mathcal S$ is decomposable.  This is a contradiction. 
Therefore $\mathcal S ^+$ is indecomposable.
\end{proof}

\bigskip
\noindent
{\bf Example}. Let ${\mathcal S}_{\gamma} = (H;E_1,E_2,E_3,E_4)$ 
be an exotic example in section 10. Since 
$E_i + E_j = H$ and $E_i \cap E_j = 0$ for distinct 
$i,j \in \{1,2,4\}$, we have 
$E_k^+ + E_m^+ = H$ and $E_k^+ \cap E_m^+ = 0$ for distinct 
$k,m \in \{3,4\}$  or $k,m \in \{1,3\}$ or $k,m \in \{2,3\}$  by 
Lemma \ref{lemma:nondegenerate1} and Lemma \ref{lemma:nondegenerate2}. 
Since $E_k^+ + E_m^+ = H$ is closed,  
$(E_k^+)^{\perp} + (E_m^+ )^{\perp}$ is closed. Hence 
$(E_k^+)^{\perp} + (E_m^+ )^{\perp} = H$
Therefore  ${\mathcal S}_{\gamma}$ is reduced from above  and 
$\Phi ^+({\mathcal S}_{\gamma})$ is reduced from below. 
Since ${\mathcal S}_{\gamma}$ is indecomposable, 
$\Phi ^+({\mathcal S}_{\gamma})$ 
is also indecomposable.

Similarly  we have the following:

\begin{thm}
Let $\mathcal S = (H;E_1,\ldots ,E_n)$ be a system
of $n$ subspaces in a Hilbert space $H$. 
Suppose that  $\mathcal S$ is reduced from below and 
$\mathcal S ^- = \Phi ^-({\mathcal S})$ is reduced from above.
If $\mathcal S$ is indecomposable, then $\Phi ^-({\mathcal S})$ 
is also indecomposable.
\end{thm}

\bigskip

We shall show that the Coxeter functors $\Phi ^+$ and $\Phi ^-$ 
preserve the defect under certain conditions.

Let $\mathcal S = (H;E_1,\ldots ,E_n)$ be a system
of $n$ subspaces in a Hilbert space $H$. Consider 
$\mathcal S^+ = (H^+;E_1^+,\dots ,E_n^+)$. Let 
$R = \oplus _{i=1}^n E_i$ and $p_0 \in B(R)$ be the 
projection of $R$ onto $H^+$.  Let $e_i \in B(H)$ be 
the projection of $H$ onto $E_i$.  Recall that 
$\tau : R \rightarrow H$ is given by 
$\tau(a) = \sum _{i=1}^n a_i$ for $a = (a_1,\dots,a_n) \in R$.   

\begin{lemma} Suppose that 
$\sum _{i=1}^n e_i$ is invertible. Then for 
$a = (a_1,\dots,a_n) \in R$ we have 
\[
p_0(a) = (a_k -e_k(\sum _{i=1}^n e_i)^{-1}(\tau (a)))_k \in H^+
\]
\label{lemma:p0}
\end{lemma}
\begin{proof}
Recall that $\tau^*: H \rightarrow R$ is given by 
$\tau^*(y) = (e_1y, \dots, e_ny)$ for $y \in H$. 
Consider the orthogonal decomposition $R = H^+ \oplus (H^+)^{\perp}$. 
Since $H^+ = \Ker \tau$, $(H^+)^{\perp} = \overline{\Im \tau^*}$ 
in $R$. Define 
\[
x = (x_k)_k := (a_k -e_k(\sum _{i=1}^n e_i)^{-1}(\tau (a)))_k \in R.
\] 
Then 
\[
\tau (x) = \sum _{k=1}^n (a_k -e_k(\sum _{i=1}^n e_i)^{-1}(\tau (a)))
         = \tau (a) - (\sum _{k=1}^n e_k)(\sum _{i=1}^n e_i)^{-1}(\tau (a)) 
         = 0 .
\]
Therefore $x \in H^+$.  Put $y := (\sum _{i=1}^n e_i)^{-1}(\tau (a)) \in H$. 
Then $\tau^*(y) = (e_1y, \dots, e_ny) \in (H^+)^{\perp}$. 
Since $a = x + \tau^*(y) \in H^+ \oplus (H^+)^{\perp}$,
we have $p_0(a) = x$.  
\end{proof}

\begin{cor} Suppose that 
$\sum _{i=1}^n e_i$ is invertible. 
Then $\Im \tau^*$ is closed and 
\[
(H^+)^{\perp} = \Im \tau^* =\{(e_1y, \dots, e_ny) \in R ; y \in H\} .
\]
\end{cor}
\begin{proof}
By the above lemma, we have 
\[
(H^+)^{\perp} = \Im (I-p_0) = \{(e_1y, \dots, e_ny) \in R ; y \in H\} 
= \Im \tau^* .
\] 
\end{proof}

\bigskip
\begin{lemma}Let $\mathcal S = (H;E_1,\ldots ,E_n)$ be a system
of $n$ subspaces in a Hilbert space $H$. Let $e_i \in B(H)$ be 
the projection of $H$ onto $E_i$. Then 
$$
\sum _{i=1}^n E_i = \Im ((\sum _{i=1}^n e_i)^{1/2}). 
$$
Moreover $\sum _{i=1}^n E_i$ is closed if and only if 
$\sum _{i=1}^n e_i$ has a closed range.  
\end{lemma}
\begin{proof} See Filmore and Williams \cite{FW} for several facts on 
operator ranges. Let $T = (T_{ij})_{ij} \in B(H^n)$ be an operator matrix 
defined by $T_{1j} = e_j$ and $T_{ij} = 0$  for $i \not= 1$. Recall that
$\Im T = \Im ((TT^*)^{1/2})$ for any operator $T$.  Since 
$\Im T = (\sum _{i=1}^n E_i) \oplus 0 \oplus 0 \oplus 0$ and 
$\Im ((TT^*)^{1/2}) = (\Im ((\sum _{i=1}^n e_i)^{1/2})) \oplus 0 \oplus 0 \oplus 0$, 
we have $\sum _{i=1}^n E_i = \Im ((\sum _{i=1}^n e_i)^{1/2}).$   

It is a known fact that $\Im A$ is closed if and only if $\Im A^{1/2}$ is 
closed for any positive operator $A \in B(H)$. This implies the rest.
\end{proof}

\bigskip
\begin{cor}Let $\mathcal S = (H;E_1,\ldots ,E_n)$ be a system
of $n$ subspaces in a Hilbert space $H$.
If $\mathcal S$ is reduced from above, then 
$f:= \sum _{i=1}^n e_i$ is invertible. 
\end{cor}
\begin{proof}
Let $x \in \Ker f$. Then $(e_ix|x) = 0$ so that $e_ix = 0$.  
Since $\mathcal S$ is reduced from above, $x \in \cap _i E_i^{\perp} = 0$ 
Thus $\Ker f = 0$. Then $\overline{\Im f} = (\Ker f)^{\perp} = H$. Since 
$\mathcal S$ is reduced from above,  
$\sum _{i=1}^n E_i = H$ is clearly closed. By the preceding lemma, 
$f$ has a closed range.  Thus $\Im f = H$.  Therefore $f$ is invertible. 
\end{proof}

\begin{lemma}
Suppose that $\mathcal S$ is reduced from above. Then for $k = 1, \dots ,n$
\[
 (E_k^+)^{\perp} = \{(\delta _{jk}a_j -e_j(\sum _{i=1}^n e_i)^{-1}(a_k))_j 
\in H^+ ; a_k \in E_k \} .
\]
\end{lemma}
\begin{proof}
Since $\mathcal S$ is reduced from above, we have 
$\Im p_kp_0 = 0 \oplus E_k \oplus$. In fact, for any $a_k \in E_k$, 
there exist $a_i \in E_i, (i \not= k)$ 
such that $-a_k = \sum _{i \not= k} a_i$.  
Then $(a_1, \dots, a_n) \in H^+$ and 
\[
p_kp_0(a_1, \dots, a_n) = (0,\dots ,0,a_k,0,\dots , 0) \in  0 \oplus E_k 
\oplus 0  .
\]
The converse inclusion is trivial. Since $\Im p_kp_0 = 0 \oplus E_k \oplus$ is 
closed, $(\Im p_kp_0)^* = \Im p_0p_k$ is also closed. Hence 
\[
(E_k^+)^{\perp} = E_k^0 = \Im p_0p_k = \{p_0(0,\dots ,0,a_k,0,\dots, 0) 
; a_k \in E_k \}
\]
Therefore the conclusion follows from Lemma \ref{lemma:p0} . 
\end{proof}

\begin{prop}
Let ${\mathcal S} = (H;E_1,E_2,E_3,E_4)$ be 
a system of four subspaces and  
${\mathcal S}^+  = (H^+;E_1^+,E_2^+,E_3^+,E_4^+)$ . 
Suppose that $\mathcal S$ is reduced from above. Then 
$f:= e_1 + e_2 + e_3 + e_4$ is invertible and 
\begin{align*}
 & (E_1^+)^{\perp} \cap (E_2^+)^{\perp} \\
& = \{(e_1u - e_1f^{-1}e_1u,-e_2f^{-1}e_1u,-e_3f^{-1}e_1u,-e_4f^{-1}e_1u) ; 
u \in E_3^{\perp} \cap E_4^{\perp}  \} .
\end{align*}
Moreover we have 
\[
\dim ((E_1^+)^{\perp} \cap (E_2^+)^{\perp})  
= \dim (E_3^{\perp} \cap E_4^{\perp}) .
\] 
The same formulae hold under permutation of subspaces. 
\label{prop:preserve perp}
\end{prop}

\begin{proof}
Let $x = (x_1,x_2,x_3,x_4) \in (E_1^+)^{\perp} \cap (E_2^+)^{\perp}$. Then 
by the preceding lemma,  there exist $a_1 \in E_1$ and $a_2 \in E_2$ such that 
\begin{align*}
x & = (x_1,x_2,x_3,x_4) \\
  & = (a_1 - e_1f^{-1}a_1, - e_2f^{-1}a_1, - e_3f^{-1}a_1 - e_4f^{-1}a_1) \\
  & = (- e_1f^{-1}a_2 , a_2 - e_2f^{-1}a_2,- e_3f^{-1}a_2 - e_4f^{-1}a_2) .
\end{align*}
Put $u := f^{-1}(a_1 - a_2) \in H$.  Then $a_1 = e_1u$, $a_2 = -e_2u$, 
$e_3u = 0$ and $e_4u = 0$.  Therefore $u \in E_3^{\perp} \cap E_4^{\perp}$ and 
\[
x = (e_1u - e_1f^{-1}e_1u,-e_2f^{-1}e_1u,-e_3f^{-1}e_1u,-e_4f^{-1}e_1u) . 
\]
Conversely suppose that 
\[
x = (e_1u - e_1f^{-1}e_1u,-e_2f^{-1}e_1u,-e_3f^{-1}e_1u,-e_4f^{-1}e_1u) , 
\]
for some $u \in  E_3^{\perp} \cap E_4^{\perp}$. Put $a_1:=e_1u \in E_1$ and 
$a_2:=-e_2u \in E_2$. Since $e_3u = 0$ and $e_4u = 0$, we have 
\[
a_1 -a_2 = e_1u + e_2u = e_1u  + e_2u + e_3u + e_4u = fu .
\]
Because $f$ is invertible, $u = f^{-1}(a_1 -a_2)$. Therefore 
\[
x = (a_1 - e_1f^{-1}a_1, - e_2f^{-1}a_1, - e_3f^{-1}a_1 - e_4f^{-1}a_1) 
\in (E_1^+)^{\perp} . 
\]
On the other hand, $a_1 = e_1u = e_1f^{-1}(a_1 -a_2)$. Hence 
\[
a_1 - e_1f^{-1}a_1 = - e_1f^{-1}a_2 .
\]
Since $a_2 = -e_2u = -e_2f^{-1}(a_1 -a_2)$, we have 
\[
- e_2f^{-1}a_1 = a_2 - e_2f^{-1}a_2 .
\]
Since 
$e_3f^{-1}(a_1 -a_2) = e_3u = 0$, we have $e_3f^{-1}a_1 = e_3f^{-1}a_2$. 
Similarly $e_4f^{-1}a_1 = e_4f^{-1}a_2$. Therefore 
\[
x = (- e_1f^{-1}a_2 , a_2 - e_2f^{-1}a_2,- e_3f^{-1}a_2 - e_4f^{-1}a_2)
\in (E_2^+)^{\perp} . 
\]
Thus $x \in (E_1^+)^{\perp} \cap (E_2^+)^{\perp}$.  

Moreover define 
$T: E_3^{\perp} \cap E_4^{\perp} \rightarrow 
(E_1^+)^{\perp} \cap (E_2^+)^{\perp}$ by 
\[
Tu = (e_1u - e_1f^{-1}e_1u,-e_2f^{-1}e_1u,-e_3f^{-1}e_1u,-e_4f^{-1}e_1u) 
\]
for $u \in E_3^{\perp} \cap E_4^{\perp}$. Then $T$ is a bounded, 
surjective operator. We shall show that $T$ is one to one.  Suppose that 
$Tu = 0$.  Since $e_2f^{-1}e_1u = 0$, $f^{-1}e_1u \in E_2^{\perp}$. 
Similarly $f^{-1}e_1u \in E_3^{\perp}$ and $f^{-1}e_1u \in E_4^{\perp}$. 
Since $\mathcal S$ is reduced from above, 
\[
f^{-1}e_1u \in E_2^{\perp} \cap E_3^{\perp} \cap  E_4^{\perp} 
= (E_2 + E_3 + E_4)^{\perp} = H^{\perp} = 0 .
\]
Hence $e_1u = 0$. Similary we have $e_2u = 0$.  Therefore 
$fu =e_1u + e_2u + e_3u + e_4u = 0$. Since $f$ is invertible, $u = 0$. 
Thus $T$ is an invertible operator. Therefore 
$\dim ((E_1^+)^{\perp} \cap (E_2^+)^{\perp})  
= \dim (E_3^{\perp} \cap E_4^{\perp})$. 
\end{proof}

\begin{thm}
Let ${\mathcal S} = (H;E_1,E_2,E_3,E_4)$ be 
a system of four subspaces.  
Suppose that $\mathcal S$ is reduced from above. 
If ${\mathcal S}$ is a quasi-Fredholm system, then 
$\Phi ^+({\mathcal S})$ is also a quasi-Fredholm system 
and 
\[\rho(\Phi ^+({\mathcal S})) = \rho ({\mathcal S}) .
\]
\end{thm}
\begin{proof}
It follows from Lemma \ref{lemma:nondegenerate1} and 
  Proposition \ref{prop:preserve perp} . 
\end{proof}

\begin{thm}
Let ${\mathcal S} = (H;E_1,E_2,E_3,E_4)$ be 
a system of four subspaces.  
Suppose that $\mathcal S$ is reduced from below. 
If ${\mathcal S}$ is a quasi-Fredholm system, then 
$\Phi ^-({\mathcal S})$ is also a quasi-Fredholm system 
and 
\[
\rho(\Phi ^-({\mathcal S})) = \rho ({\mathcal S}) .
\]
\end{thm}
\begin{proof}
Recall that $\mathcal S$ is reduced from below if and only if 
$\Phi ^{\perp}({\mathcal S})$ is reduced from above, and 
${\mathcal S}$ is a quasi-Fredholm system if and only if 
$\Phi ^{\perp}({\mathcal S})$ is a quasi-Fredholm system. 
Applying the preceding theorem, 
$\Phi ^-({\mathcal S}) = \Phi ^{\perp}\Phi ^+\Phi ^{\perp}(\mathcal S)$
is a quasi-Fredholm system and 
\[
\rho(\Phi ^-({\mathcal S})) = -\rho(\Phi ^+\Phi ^{\perp}(\mathcal S)) 
= - \rho(\Phi ^{\perp}(\mathcal S)) = \rho ({\mathcal S}) .
\]
\end{proof}

\bigskip
\noindent
{\bf Example}. 
Let $\mathcal S$  be an operator system. 
Since $E_{1}=K\oplus 0,E_{2}=0\oplus K$, we have that 
$f = \sum_{i=1}^{4}e_{i}  \geq I$ is invertible. 
Moreover if $\mathcal S =  {\mathcal S}_T$ is associated with 
a single bounded operator $T$, then $E_4 = \{(x,x) \in H ; x \in K \}$.
Thus  $E_i + E_j = H$ for $(i,j) = (1,2), (1,4),(2,4)$ and 
$\mathcal S$ is reduced from above. Therefore, 
if ${\mathcal S}_T$ is a quasi-Fredholm system, then 
$\Phi ^+({\mathcal S}_T)$ is also a quasi-Fredholm system 
and $\rho(\Phi ^+({\mathcal S}_T)) = \rho ({\mathcal S}_T)$. 
Similarly, let  ${\mathcal S}_{\gamma}$ 
be an exotic example in section 10. Then 
${\mathcal S}_{\gamma}$ is reduced from above and $f$
is invertible. Since ${\mathcal S}_{\gamma}$ is a quasi-Fredholm system, 
$\Phi ^+({\mathcal S}_{\gamma})$ is also a quasi-Fredholm system 
and $\rho(\Phi ^+({\mathcal S}_{\gamma})) = \rho ({\mathcal S}_{\gamma})$.

\end{document}